\documentclass[11pt]{amsart}
\usepackage[margin= 1in]{geometry}
\usepackage{graphicx,dsfont}
\usepackage{amssymb}
\usepackage{mathtools}
\usepackage{tikz-cd}
\usepackage{color}
\usepackage[title]{appendix}
\usepackage{amsmath, mathrsfs, stmaryrd, color,slashed,bbm}

\usepackage{float}
\usepackage{verbatim}
\usepackage{comment}
\usepackage[utf8]{inputenc}
\usepackage[english]{babel}
\usepackage{enumitem}
\usepackage{booktabs}

\newtheorem{theorem}{Theorem}
\newtheorem{remark}[theorem]{Remark}
\newtheorem{claim}[theorem]{Claim}
\newtheorem{definition}[theorem]{Definition}
\newtheorem{proposition}[theorem]{Proposition}
\newtheorem{corollary}[theorem]{Corollary}
\newtheorem{lemma}[theorem]{Lemma}

\newcommand{\free}{\textup{free}}

\newcommand{\N}{\mathbb{N}}
\newcommand{\R}{\mathbb{R}}
\newcommand{\Z}{\mathbb{Z}}
\newcommand{\E}{\mathbb{E}}
\newcommand{\bP}{\mathbb{P}}

\newcommand{\whitenoise}{\ensuremath{\mathscr{\dot{W}}}}
\newcommand{\Ham}{\mathbf{H}}
\newcommand{\cL}{\mathcal{L}}
\newcommand{\Fext}{\mathcal{F}_{\textrm{ext}}}

\numberwithin{theorem}{section}
\numberwithin{equation}{section}

\begin{document}
\title[]{Convergence of the KPZ line ensemble}
\author{Xuan Wu}
\address{Department of Mathematics, University of Chicago, 
Chicago IL, 60637}
\email{xuanw@uchicago.edu}

\begin{abstract}
In this paper we study the KPZ line ensembles under the KPZ scaling. Based on their Gibbs property, we derive quantitative local fluctuation estimates for the scaled KPZ line ensembles. This allows us to show that the family of scaled KPZ line ensembles is tight. Together with the recent progresses in \cite{QS20}, \cite{Vir} and \cite{DM}, the tightness result yields the conjectural convergence of the scaled KPZ line ensembles to the Airy line ensemble.
\end{abstract}

\vspace{1.5cm}

\maketitle

\section{Introduction}
\subsection{Kardar-Parisi-Zhang equation}
The main object of study in this paper, the KPZ line ensemble, can be viewed as a multi-layer extension to the well-known Kardar-Parisi-Zhang (KPZ) equation. The KPZ equation was introduced in 1986 by Kardar, Parisi and Zhang \cite{KPZ} as a model for random interface growth. In one-spatial dimension (sometimes also called (1+1)-dimension to emphasize that time is one dimension too), it describes the evolution of a function $\mathcal{H}(t, x)$ recording the height of an interface at time $t > 0$ above position $x\in\R$. The KPZ equation is written formally as a stochastic partial differential equation (sPDE), 
\begin{equation}\label{eq:KPZ}
\partial_t\mathcal{H} = \frac{1}{2}\partial_x^2 \mathcal{H} + \frac{1}{2}(\partial_x\mathcal{H})^2 + \whitenoise,
\end{equation}
where $\whitenoise$ is a space-time white noise (for mathematical background or literature review, see \cite{Cor, QS15} for instance). 

The KPZ equation \eqref{eq:KPZ} is associated with a famous universality class, the KPZ universality class, which bears the same name. The KPZ equation is a canonical member of the associated KPZ universality class and a model belongs to the KPZ universality class if it bears the same long-time, large-scale behavior as the KPZ equation. The KPZ universality class hosts a large class of models, which covers a wide range
of mathematical and physical systems of distinct origins, including interacting particle systems, random
matrices, traffic models, directed polymers in random media and non-linear stochastic PDEs. 

All models in the KPZ universality class can be transformed to a kinetically growing interface reflecting the competition between growth in a direction normal to the surface, a surface tension smoothing force, and a stochastic term which tends to roughen the interface. These features may be illustrated by the Laplacian $\frac{1}{2}\partial_x^2 \mathcal{H}$, non-linear term $\frac{1}{2}(\partial_x\mathcal{H})^2$ and white noise $\whitenoise$ in the KPZ equation \eqref{eq:KPZ}. Numerical simulations along with some theoretical results have confirmed that in the long time $t$ scaling limit, fluctuations in the height of such evolving interfaces scale like $t^{1/3}$ and display non-trivial spatial correlations in the scale $t^{2/3}$ (known as the $3:2:1$ KPZ scaling).

The KPZ equation is related to the stochastic heat equation (SHE) with multiplicative noise through the Hopf–Cole transformation. Denote $\mathcal{Z}(t, x)$ as the solution to the following SHE, 
\begin{equation}\label{eq:SHE}
\partial_t \mathcal{Z} = \frac{1}{2}\partial_x^2 \mathcal{Z} + \mathcal{Z} \whitenoise.
\end{equation} 
The Hopf-Cole solution to the KPZ equation \eqref{eq:KPZ} is defined through taking $$\mathcal{H}(t,x)=\log \mathcal{Z}(t,x).$$  It was first proved in \cite{Mue} that $\mathcal{Z}(t,x)$ is almost surely strictly positive (with positive initial conditions), which justifies the validity of the transform. The fundamental solution $\mathcal{Z}^{nw}(t, x)$ to SHE \eqref{eq:SHE} is of great importance. It solves \eqref{eq:SHE} with a delta mass initial value at origin, i.e. $\mathcal{Z}(0, x)=\delta_{x=0}$. Meanwhile $\mathcal{H}^{nw}(t, x) = \log \mathcal{Z}^{nw}(t, x)$ is known as the narrow wedge solution to the KPZ equation. The initial condition of $\mathcal{H}^{nw}(0, x)$ is not well-defined; however $\mathcal{H}^{nw}(t, x)$ is stationary around a parabola $\frac{-x^2}{2t}$, which resembles a sharp wedge for small $t$, hence known as the narrow wedge initial condition.

Using the Feynman-Kac representation, $\mathcal{Z}^{nw}(t, x)$ formally takes the following expression,
\begin{equation}\label{eq:Z^nw}
\mathcal{Z}^{nw}(t, x)=p(t,x)\E\left[:\exp :\left(\int_0^t\whitenoise(s,B_s)ds\right)\right],
\end{equation}
where $p(t,x)=(2\pi t)^{-1/2} \exp(-x^2/2t)$ is the heat kernel, the expectation is taken with respect to a Brownian bridge $(B_s,s\leq t)$ which starts at origin at time $0$ and arrives at $x$ at time $t$. The $:{\rm exp}:$ is the {\it Wick exponential}, see \cite{Jan} for instance. This bridge representation arises because of the $\delta_{x=0}$ initial condition, and hence the factor $p(t,x)$. This Feynman-Kac representation is mostly formal since the integral of white noise over a Brownian path is not well-defined pathwise or to exponentiate the integral.

\begin{figure}[H]
\begin{center}
\begin{tabular}{lr}
      \raisebox{-10 ex}{\includegraphics[width=5cm]{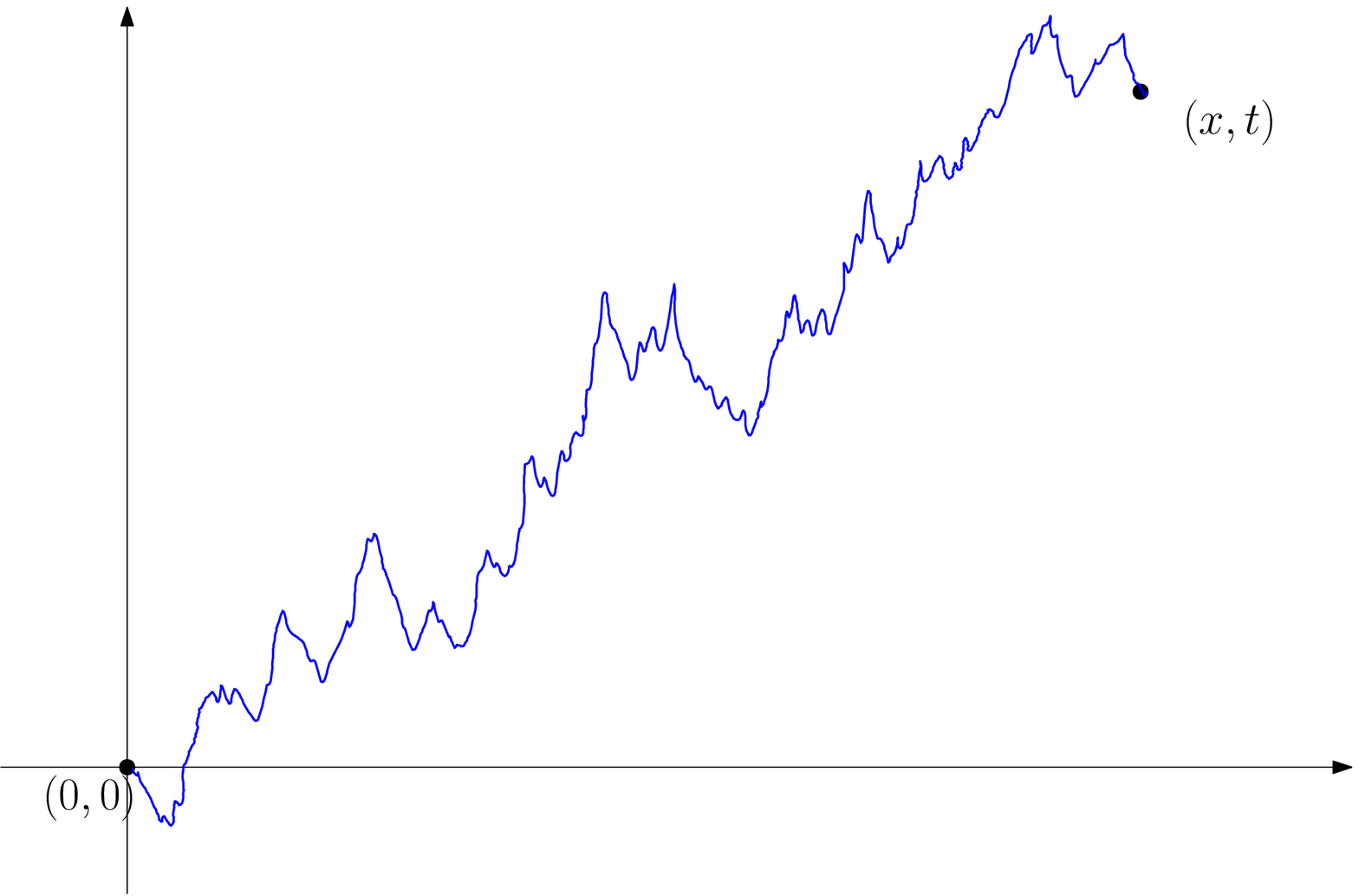}}
        &
           \raisebox{-10 ex}{\includegraphics[width=5cm]{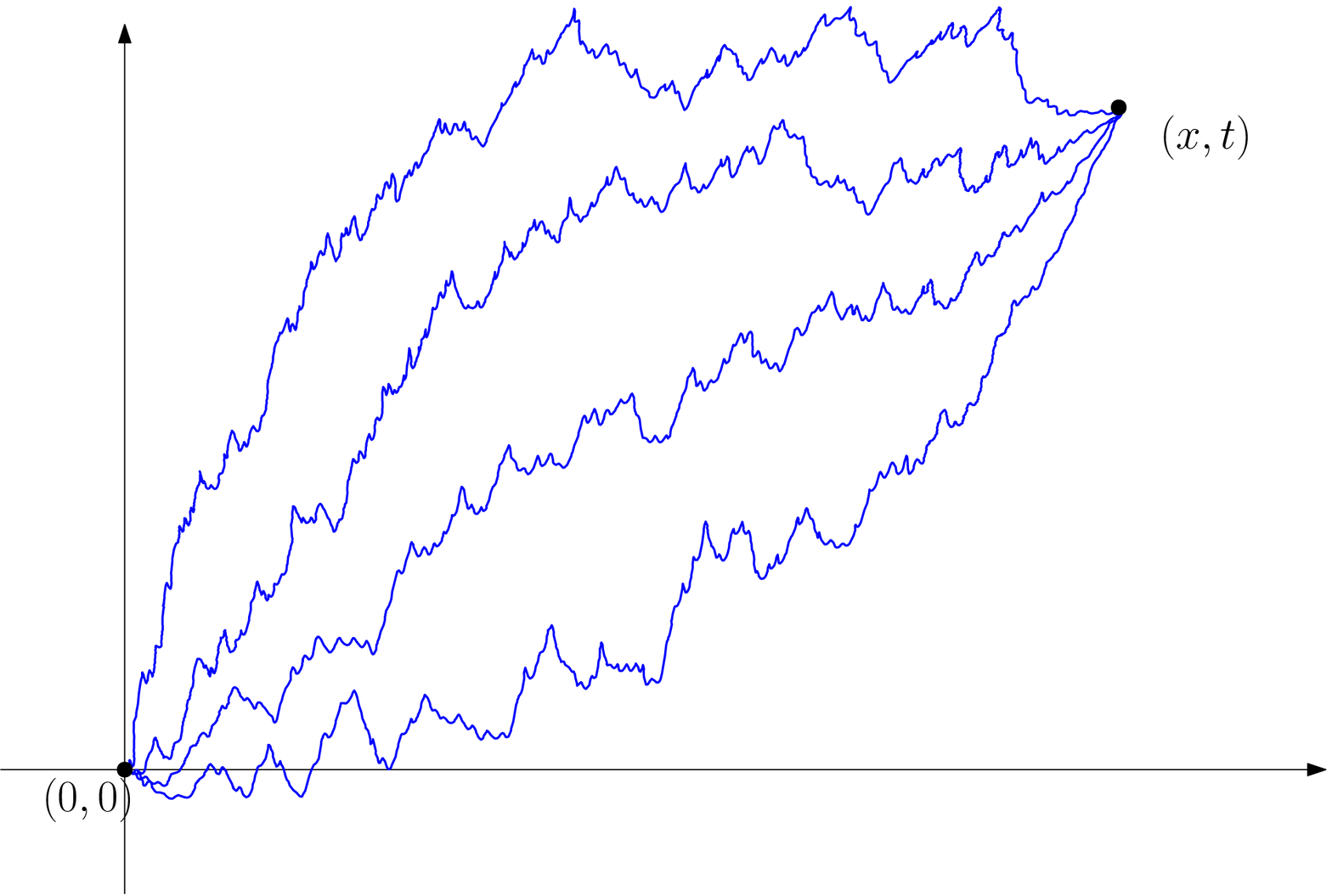}}  
        \end{tabular}
\end{center}
\caption{On the left, the CDRP can be formulated in terms of the general directed polymer framework, i.e. a random path interacting with a random environment. More precisely, the underlying path measure is that of a Brownian bridge which starts at $0$ at time $0$ and ends at $x$ at time $t$; and the random environment is the white noise $\whitenoise$. The multi-layer extension on the right involves non-intersecting Brownian bridges and is defined similarly.}\label{figure:polymer}
\end{figure}

We adopt this representation to emphasize on its interpretation as being the partition function of a continuum directed random polymer (CDRP) that is modeled by a continuous path interacting with a space-time white noise. We emphasize that this approach is very useful for a generalization to the multi-layer scenario which involves multiple Brownian bridges, see an illustration in Figure \ref{figure:polymer}. We will make sense of the expression \eqref{eq:Z^nw} through a chaos expansion after we introduce its multi-layer extension below.

\subsection{The KPZ line ensemble and CDRP partition functions}\label{sec:KPZLE}
Motivated by recent developments on solvable directed polymer models, \cite{OW} generalizes the above Feynman-Kac representation / continuum polymer partition function $\mathcal{Z}(t,x)$ to accommodate multiple non-intersecting Brownian bridges. More precisely, they defined a hierarchy of partition functions $\mathcal{Z}_n(t,x)$, formally written as 
\begin{equation}\label{wick}
\mathcal{Z}_{n}(t,x) = p(t,x)^n \E\left[:{\rm exp}:\, \left\{\sum_{i=1}^{n} \int_{0}^{t} \whitenoise(s,B_i(s)) ds \right\} \right],
\end{equation}
where the expectation is taken with respect to the law on $n$ independent Brownian bridges $\{B_i\}_{i=1}^{n}$ starting at $0$ at time $0$ and ending at $x$ at time $t$. Intuitively these path integrals represent energies of non-intersecting paths, and thus the expectation of their exponential represents the partition function for this path (or directed polymer) model. It is worth noting that the first layer, $\mathcal{Z}_1$, is the sames as the fundamental solution $\mathcal{Z}^{nw}$ to the SHE \eqref{eq:SHE}. These partition functions $\mathcal{Z}_n(t,x)$ also solve a multi-layer extension of SHE \eqref{eq:SHE}, see \cite{OW}.

For $n\in \N$, $t\geq 0$ and $x\in \R$, $\mathcal{Z}_n(t,x)$ is rigorously defined via the following chaos expansion,
\begin{equation}\label{Zpartfunc}
\mathcal{Z}_{n}(t,x) = p(t,x)^n \sum_{k=0}^{\infty} \int_{\Delta_k(t)}\int_{\R^k} R_k^{(n)}\left((t_1,x_1),\cdots, (t_k,x_k)\right) \whitenoise(dt_1 dx_1)\cdots \whitenoise(dt_k dx_k),
\end{equation}
where $\Delta_k(t) = \{0<t_1<\cdots <t_k<t\}$, and $R_k^{(n)}$ is the $k$-point correlation function for a collection of $n$ non-intersecting Brownian bridges each of which starts at $0$ at time $0$ and ends at $x$ at time $t$. For notational simplicity, set $\mathcal{Z}_0(t,x)\equiv 1$. For details about integration with respect to a white noise, we refer to \cite{Jan}.

\cite{LW} show that for any $t>0$, with probability 1, for all $x\in \R$ and all $n\in \N$, $\mathcal{Z}_n(t,x)>0$. The positivity result permits the following important definition of the KPZ line ensemble $\mathcal{H}_n^t(x)$, a process given by the logarithm of ratios of partition functions $\mathcal{Z}_n$. 
\begin{definition}\label{def:KPZlineensemble}
For $t>0$ fixed, the KPZ line ensemble is a continuous $\N$-indexed line ensemble  $\mathcal{H}^t = \{\mathcal{H}^t_{n}\}_{n\in\N}$  on $\R$ given by
\begin{equation}
\mathcal{H}^t_{n}(x) := \log \left(\frac{\mathcal{Z}_n(t,x)}{\mathcal{Z}_{n-1}(t,x)}\right),
\end{equation}
with convention $\mathcal{Z}_0\equiv 1$. Note that $\mathcal{H}^{t}_1(\cdot)$ equals $\mathcal{H}^{nw}(t,\cdot)$, the time $t$ spatial process of the Hopf-Cole solution to the KPZ equation (\ref{eq:KPZ}) with narrow-wedge initial data. Sometimes we will omit $t$ and just write the KPZ line ensemble as $\mathcal{H}$.
\end{definition}

\begin{figure}
     \begin{center}
\includegraphics[width=0.5\textwidth]{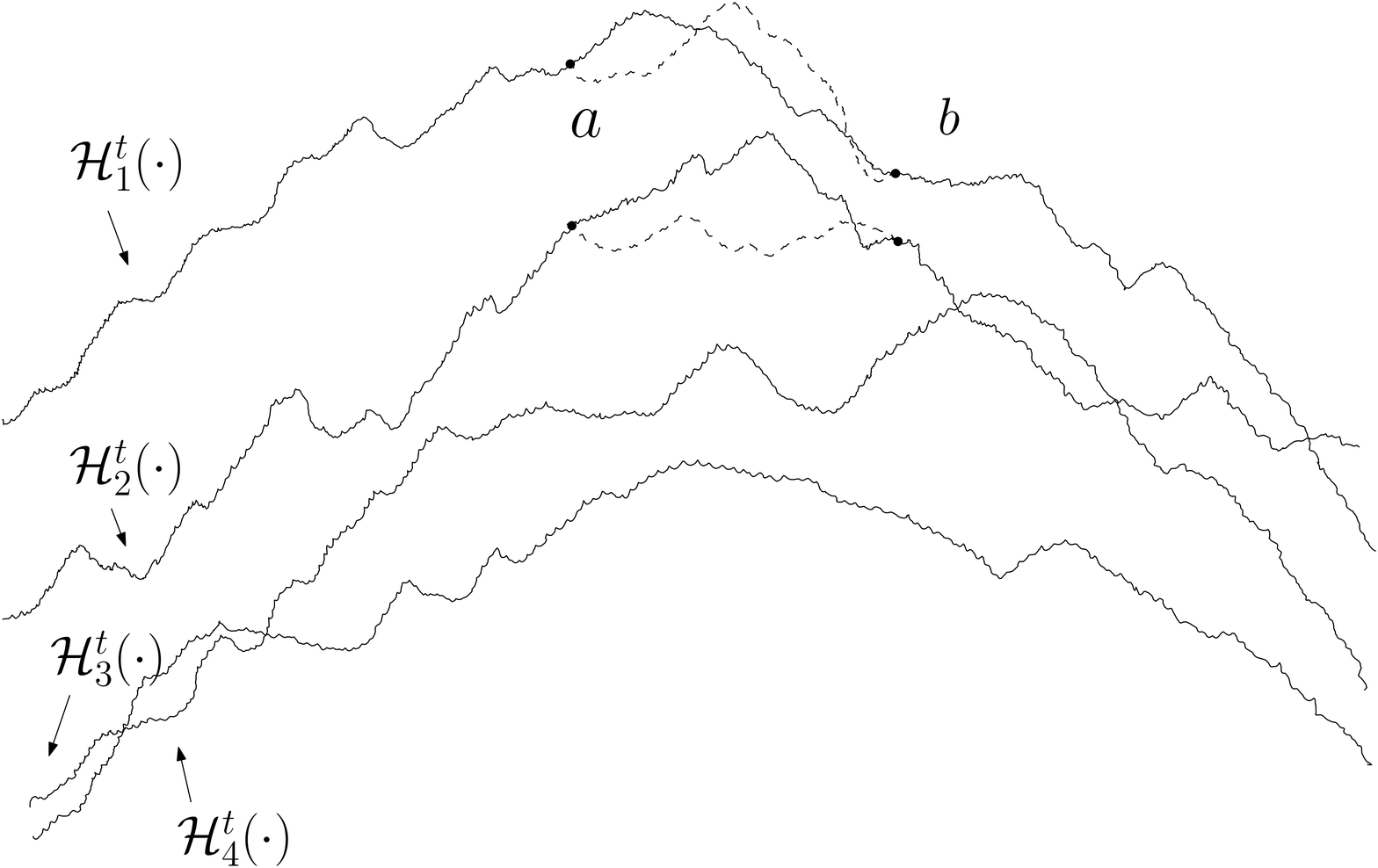}
     \caption{An overview of $\textrm{KPZ}_t$ line ensemble $\mathcal{H}^t$ for fixed time $t$. Curves $\mathcal{H}^t_1(\cdot)$ through $\mathcal{H}^t_4(\cdot)$ are sampled. The lowest indexed curve $\mathcal{H}^t_1$ is distributed according to the time $t$ solution to KPZ equation with narrow wedge initial data. The dotted first two curves $\mathcal{H}^t_1$ and $\mathcal{H}^t_2$ between $a$ and $b$ indicate a possible resampling according to two independent Brownian bridges which subject to an exponential penalization for going out of order (Definition \ref{def:H-BGP}), as a demonstration for the $\Ham$-Brownian Gibbs property.}
      \label{KPZLE}
      \end{center}
      \end{figure} 
 
The KPZ line ensemble arises as scaling limits of O’Connell-Yor semi-discrete direct polymers \cite{OY} and log-gamma discrete directed polymers \cite{Sep}. For O’Connell-Yor semi-discrete direct polymers, the convergence follows from the subsequential compactness result in \cite{CH16} and finite dimensional convergence in \cite{Nic}. The Gibbs property of $\mathcal{H}$ was introduced in \cite{CH16}, which plays a key role in their investigation. The analogue result about sequential compactness and limiting Gibbs property for log-gamma line ensemble was carried out in \cite{Wu19}.

\subsection{$\mathbf{H}$-Brownian Gibbs property}\label{sec:gibbs}
In this subsection, we describe the key tool we exploit in this paper, the $\mathbf{H}$-Brownian Gibbs property of the KPZ line ensemble.

An $\mathbf{H}$-Brownian Gibbs property can be viewed as a resampling invariance. Let $\mathcal{L}=\{\mathcal{L}_1,\mathcal{L}_2,\dots\}$ be a line ensemble, i.e., countably many random functions indexed by $\mathbb{N}$. We may also view $\cL$ as a random variable taking values in $C(\mathbb{N}\times\mathbb{R},\mathbb{R})$. Then $\mathcal{L}$ satisfies the $\mathbf{H}$-Brownian Gibbs property if for any integers $k_1\leq k_2$ and $[a,b]\subset\mathbb{R}$, the law of $\mathcal{L}$ is unchanged if one replaces the restriction of $\mathcal{L}$ on $[k_1,k_2]_{\Z}\times [a,b]$ by independent Brownian bridges, re-weighted by an interaction factor (determined by \textbf{H}) between adjacently indexed curves. See Definition~\ref{def:H_Brownian} for the precise definition.

There is a large class of stochastic integrable models from random matrix theory, interacting particle systems, last passage percolation and directed polymers that enjoy certain $\mathbf{H}$-Brownian Gibbs properties. Among them, the parabolic Airy line ensemble $\tilde{\mathcal{A}}$ (defined in Definition~\ref{def:Airy} in the next subsection) is the most celebrated one and it satisfies the $\mathbf{H}_{\textup{ord}}$-Brownian Gibbs property. See Definition~\ref{def:todo} for the expression of $\mathbf{H}_{\textup{ord}}$. We remark that in the literature, the $\mathbf{H}_{\textup{ord}}$-Brownian Gibbs property is referred to as the Brownian Gibbs property. 

After the discovery of the $\mathbf{H}_{\textup{ord}}$-Brownian Gibbs property of $\tilde{\mathcal{A}}$ in \cite{CH14}, this Gibbs property has served as a powerful probability tool in studying the Airy line ensemble and other line ensembles enjoying same Gibbs property. Recently, one of the authors of \cite{CH14}, developed a more delicate treatment in \cite{Ham1} to estimate the modulus of continuity for line ensembles with the $\mathbf{H}_{\textup{ord}}$-Brownian Gibbs property (e.g. the Airy line ensemble and the line ensemble associated with Brownian last passage percolation). \cite{Ham1} also established $L^{\infty-}$ bounds on the Radon-Nikodym derivative of the line ensemble curves with respect to Brownian bridges and other refined regularity properties. Furthermore in the subsequent papers \cite{Ham2,Ham3,Ham4}, the work in \cite{Ham1} was applied to understanding the geometry of last passage paths in Brownian last passage percolation with more general initial data. Another breakthrough is the construction of the directed landscape \cite{DOV}, based on the $\mathbf{H}_{\textup{ord}}$-Brownian Gibbs property of the Brownian last passage percolation.

The KPZ line ensemble (defined in Definition~\ref{def:KPZlineensemble}) is an central integrable model in the KPZ universality class. It satisfies the $\mathbf{H}$-Brownian Gibbs property with $\mathbf{H}(x)=e^x$. In view of the successful applications of the $\mathbf{H}_{\textup{ord}}$-Brownian Gibbs property, it is natural to exploit this Gibbs property  to study the KPZ line ensemble. In this paper, we employ this strategy and obtain regularity and convergence results of the scaled KPZ line ensembles. We discuss these results in the next subsection Section~\ref{sec:main results} and ideas in Section~\ref{sec:proofideas}.

\subsection{Main results}\label{sec:main results}
This paper aims to investigate the KPZ line ensemble under the KPZ scaling. The scaled KPZ line ensemble $\mathfrak{H}^t$ is defined as follows.
\begin{definition}\label{def:scaledKPZLE}
For $t>0$ fixed, the scaled KPZ line ensemble is a continuous $\N$-indexed line ensemble $\mathfrak{H}^t = \{\mathfrak{H}^t_{n}\}_{n\in\N}$  on $\R$ given by
\begin{equation}\label{eq:scaledKPZLE}
\mathfrak{H}^t_n(x):=\frac{\mathcal{H}^t_n(t^{2/3}x)+\frac{t}{24}}{t^{1/3}}+ (n-1)t^{-1/3}\log t^{2/3}. 
\end{equation}
Sometimes we will omit $t$ and just write $\mathfrak{H}$.
\end{definition}

\begin{remark}
Our definition of $\mathfrak{H}^t$ is slightly different from the one given in \cite[Theorem 2.15]{CH16}. The reason is to obtain a simpler Gibbs property. See Corollary~\ref{cor:KPZGibbs}. 
\end{remark}

The {\em main results} in this paper are Theorem~\ref{thm:main1} and Theorem~\ref{thm:main2}. We established a quantitative local fluctuation estimate in Theorem \ref{thm:main1}. The analogue of Theorem~\ref{thm:main1} for the Airy line ensemble is proved in \cite[Theorem 2.14]{Ham1}. With the local fluctuation estimates, we are able to show the tightness of the scaled KPZ line ensembles $\mathfrak{H}^t$ as $t$ varies (Theorem \ref{thm:main2}(i)). Furthermore, as $t$ increases, curves in $\mathfrak{H}^t$ becomes more and more ordered and the limit satisfies the $\mathbf{H}_{\textup{ord}}$-Brownian Gibbs property, see Theorem \ref{thm:main2}(ii). Theorem \ref{cor:main} is the {\em main application} of Theorem~\ref{thm:main2}, which proves that the scaled KPZ line ensemble converges to the parabolic Airy line ensemble. The local fluctuation result Theorem~\ref{thm:main1} is applied in \cite{Wu21} to study the Brownian regularity of the curves in $\mathfrak{H}^t$ (with an affine shift) with respect to Brownian bridges.

\begin{theorem}[Local fluctuation estimates]\label{thm:main1}
Fix $k\in\mathbb{N}$. There exists a constant $D_0=D_0(k)$ depending only on $k$ such that the following statement holds. For all $t\geq 1$,  $d\in (0,1]$ and $K\geq 0$, we have
\begin{align}\label{osc}
\mathbb{P}\left(\sup_{u,v\in [0,d]}\left|\mathfrak{H}_k^t(u)-\mathfrak{H}_k^t(v) \right|\geq Kd^{1/2}\right)\leq  D_0^{-1} e^{-D_0 K^{3/2}}.
\end{align}
\end{theorem}

\begin{theorem}[Tightness of the scaled KPZ line ensemble]\label{thm:main2}
\hfill
\begin{enumerate}[label=(\roman*)]
\item As $t$ varies, $\mathfrak{H}^t$ is tight in the following sense. Given an increasing sequence $t_N$ converging to infinity, there exists a subsequence, denoted by $t_{N_j}$, such that $\mathfrak{H}^{t_{N_j}}$ converges weakly (see Definition \ref{def:weakconvergence}) as $\mathbb{N}\times \mathbb{R}$-indexed line ensembles.
\item Any subsequential limit $\mathfrak{H}^{\infty}$ is a non-intersecting line ensemble and enjoys the $\mathbf{H}_{\textup{ord}}$-Brownian Gibbs property.
\end{enumerate}
\end{theorem}

\begin{remark}
The fluctuation bound (Theorem~\ref{thm:main1}) is weaker compared to the one for Brownian motions (Lemma~\ref{lem:BB-osc}) because $e^{-K^2}$ decays faster than $e^{-K^{3/2}}$. This may be viewed as the analogue of the \textit{high jump difficulty} discussed in \cite{Ham1}. In the follow-up paper \cite{Wu21}, we further improve \eqref{osc} by developing a soft jump ensemble method, inspired by \cite{Ham1}.  
\end{remark}

A longstanding conjecture about the KPZ equation is that its solution converges to the KPZ fixed point under the KPZ scaling. Recently, \cite{QS20} and \cite{Vir} independently made the breakthrough and gave an affirmative answer to this conjecture. Before that, one point convergence was independently proved in \cite{SS, ACQ} for the narrow wedge initial condition. 

As a major application of Theorem~\ref{thm:main2}, the convergence can be extended to the level of line ensembles. That is, the scaled KPZ line ensemble (defined in Definition~\ref{def:scaledKPZLE}) converges to the Airy line ensemble. The line ensemble convergence was conjectured in \cite{CH16} where the (scaled) KPZ line ensembles were introduced. To state this result, let us give the definition of the Airy line ensemble, which was first introduced as a multilayer Airy process in \cite{PS}. 
 
\begin{definition}\label{def:Airy}
The Airy line ensemble $\mathcal{A} = \{\mathcal{A}_1 > \mathcal{A}_2 > \cdots \}$ is countable many random functions indexed by $\mathbb{N}$. The law of $\mathcal{A}$ is uniquely determined by its determinantal structure. More precisely, for any finite set $I=\{u_1, \cdots, u_k\}\subset\mathbb{R}$, the point process on $I\times\mathbb{R}$ given by $\{(s,\mathcal{A}_n(s)):n\in\mathbb{N},s\in I\}$ is a determinantal point process with kernel given by  
\begin{equation*}
K(s_1,x_1;s_2,x_2)=\left\{ \begin{array}{cc}
\int_0^{\infty} e^{-z(s_1-s_2)}\textrm{Ai}(x_1+z)\textrm{Ai}(x_2+z)dz & \textrm{if} \quad s_1\geq s_2,\\
-\int_{-\infty}^0 e^{-z(s_1-s_2)}\textrm{Ai}(x_1+z)\textrm{Ai}(x_2+z)dz  & \textrm{if} \quad s_1< s_2,
\end{array} \right.
\end{equation*}
where $\textrm{Ai}$ is the Airy function. The parabolic Airy line ensemble $\tilde{\mathcal{A}}=\{\tilde{\mathcal{A}}_1>\tilde{\mathcal{A}}_2>\dots\}$ is defined by
\begin{equation*}
\tilde{\mathcal{A}}_n(x):=2^{-1/3}\mathcal{A}_n(2^{-1/3}x)-2^{-1}x^2.
\end{equation*}
\end{definition}

It is proved in \cite{CH14} that the parabolic Airy line ensemble $\tilde{\mathcal{A}}$ enjoys the $\mathbf{H}_{\textup{ord}}$-Brownian Gibbs property. See Definition~\ref{def:todo} for the definition of the $\mathbf{H}_{\textup{ord}}$-Brownian Gibbs property. 

We are ready to state the line ensemble convergence result.
\begin{theorem}\label{cor:main}
When $t$ goes to infinity, the scaled KPZ line ensemble $\mathfrak{H}^t$ (defined in Definition~\ref{def:scaledKPZLE}) converges weakly to the parabolic Airy line ensemble $\tilde{\mathcal{A}}$ (defined in Definition~\ref{def:Airy}). See Definition \ref{def:weakconvergence} for the definition of the weak convergence of line ensembles.
\end{theorem}

\begin{proof}
By Theorem \ref{thm:main2}, the scaled KPZ line ensemble $\mathfrak{H}^t$ is tight for $t\geq 1$ and any subsequential limit satisfies the $\mathbf{H}_{\textup{ord}}$-Brownian Gibbs property. It is proved in \cite{QS20} and \cite{Vir} that $\mathfrak{H}^t_1$ converges weakly to $\tilde{\mathcal{A}}_1$. Moreover, \cite[Theorem 1.1]{DM} shows that a $\mathbf{H}_{\textup{ord}}$-Brownian Gibbsian line ensemble is completely characterized by the finite dimensional distributions of its top curve. In view of the $\mathbf{H}_{\textup{ord}}$-Brownian Gibbs property of $\tilde{\mathcal{A}}$, this implies that any subsequential limit of $\mathfrak{H}^t$ agrees with $\tilde{\mathcal{A}}$. Hence $\mathfrak{H}^t$ converges weakly to $\tilde{\mathcal{A}}$ as $t$ goes to infinity.
\end{proof}

\subsection{Framework to analyze Gibbsian line ensembles}\label{sec:proofideas}
In this subsection, we present a general strategy to analyze Gibbsian line ensembles. As discussed in Subsection~\ref{sec:gibbs}, this strategy has been very successful in studying the parabolic Airy line ensemble or other $\mathbf{H}_{\textup{ord}}$-Brownian Gibbs line ensembles. However, new difficulties arise when trying to fit the KPZ line ensembles into this framework. The novelty of this paper is to resolve these difficulties through a two-step inductive resampling procedure, see more details in Section~\ref{sec:Z}.

To describe such a strategy, we need to introduce some notations. We focus on line ensembles consisting of three random continuous curves. See Definition~\ref{def:H_Brownian} for the general setup. For any $\vec{x}=(x_1,x_2,x_3), \vec{y}=(y_1,y_2,y_3)\in\mathbb{R}^3$ and $(a,b)\subset\mathbb{R}$, we write $\bP^{(a,b),\vec{x},\vec{y}}_{\free}$ for the probability measure on $C([1,3]_{\Z}\times [a,b],\R)$ which is induced by independent Brownian bridges $B=(B_1,B_2,B_3)$ with $B(a)=\vec{x}$ and $B(b)=\vec{y}$. Consider a line ensemble $\mathcal{L}=(\mathcal{L}_1,\cL_2,\cL_3)$ and a Hamiltonian function $\mathbf{H}:\R\to [0,\infty]$. We write $\bP_{\cL}$ for the law of $\cL$. Let $\bP^{(a,b),\vec{x},\vec{y}}_{\cL}$ be the probability measure on $C([1,3]_{\Z}\times [a,b],\R)$ induced by $\cL$ conditioned on $\mathcal{L}(a)=\vec{x}$ and $\mathcal{L}(b)=\vec{y}$. Then $\cL$ enjoys the $\mathbf{H}$-Brownian Gibbs property if for any $\vec{x},\vec{y},(a,b)$ as above, it holds that
\begin{equation}\label{equ:RN_intro}
\frac{\textup{d} \bP^{(a,b),\vec{x},\vec{y}}_{\cL}}{\textup{d} \bP^{(a,b),\vec{x},\vec{y}}_{\free}}(\tilde{\cL})=\frac{W_{\mathbf{H}}^{(a,b),\vec{x},\vec{y}}(\tilde{\cL})}{Z_{\mathbf{H}}^{(a,b),\vec{x},\vec{y}}}. 
\end{equation} 
Here the Boltzmann weight $W_{\mathbf{H}}^{(a,b),\vec{x},\vec{y}}$ is a functional on $C([1,3]_{\Z}\times [a,b],\R)$ defined through
\begin{align*}
W_{\mathbf{H}}^{(a,b),\vec{x},\vec{y}}(\tilde{\cL}):=\exp\left( - \sum_{i=1}^2 \int_a^b\mathbf{H}(\tilde{\cL}_{i+1}(u)-\tilde{\cL}_i(u))\, du \right),
\end{align*}
and $Z_{\mathbf{H}}^{(a,b),\vec{x},\vec{y}}$ is the normalizing constant which equals the expectation of $W_{\mathbf{H}}^{(a,b),\vec{x},\vec{y}}$ with respect to $\bP^{(a,b),\vec{x},\vec{y}}_{\free}$. Let $\mu^{(a,b)}$ be the probability measure on $\R^3\times\R^3$ induced by $(\cL(a),\cL(b))$. 

Let $\mathsf{A}\subset  C([1,3]_{\Z}\times [-1,1],\R)$ be a Borel subset. We abuse the notation and denote by $\mathsf{A}$ the event that (the restriction of) random curves are contained in $\mathsf{A}$. Then from \eqref{equ:RN_intro}, we have
\begin{align*}
\mathbb{P}_{\cL}(\mathsf{A})=&\int_{\R^3\times\R^3} d\mu^{(-1,1)}(\vec{x},\vec{y})\ \bP_{\cL}(\mathsf{A}\, |\, (\cL(-1),\cL(1))=(\vec{x},\vec{y})) \\
=& \int_{\R^3\times\R^3} d\mu^{(-1,1)}(\vec{x},\vec{y})\ \left(Z^{(-1,1),\vec{x},\vec{y}}_{\mathbf{H}}\right)^{-1} \mathbb{E}^{(-1,1),\vec{x},\vec{y}}_{\free}\left[ \mathbbm{1}_{\mathsf{A}}\cdot W_{\mathbf{H}}^{(-1,1),\vec{x},\vec{y}} \right]. 
\end{align*}
Let $\mathsf{GB}$ be a Borel subset of $\R^3\times\R^3$. We again abuse the notation and denote by $\mathsf{GB}$ the event that $\{(\cL(-1),\cL(1))\in \mathsf{GB}\}$ and we think of $\mathsf{GB}$ as the collection of good boundary data $(\cL(-1),\cL(1))$. Multiplying the above equality by $1=\mathbbm{1}_{\mathsf{GB}}+\mathbbm{1}_{\mathsf{GB}^{\textup{c}}}$ and using the fact $W^{(-1,1),\vec{x},\vec{y}}_{\mathbf{H}}\leq 1$, we obtain
\begin{align}\label{equ:key}
\mathbb{P}_{\cL}(\mathsf{A})\leq \mathbb{P}_{\cL}(\mathsf{GB}^{\textup{c}})+ \sup_{(\vec{x},\vec{y})\in\mathsf{GB}} \left(Z^{(-1,1),\vec{x},\vec{y}}_{\mathbf{H}}\right)^{-1} \mathbb{P}^{(-1,1),\vec{x},\vec{y}}_{\free}(\mathsf{A}) .
\end{align}
The inequality \eqref{equ:key} provides a setup to bound $\mathbb{P}_{\cL}(\mathsf{A})$ for general $\mathsf{A}$. Now we consider the concrete case $\mathsf{A}=\mathsf{BigFluc}$, where
\begin{align*}
\mathsf{BigFluc}:=\left\{\displaystyle\sup_{0\leq u<v\leq d} |\tilde{\cL}_1(u)-\tilde{\cL}_1(v)|\geq K\sqrt{d}\right\} 
\end{align*}
for some fixed $K>0$ and $d\in (0,1)$. Suppose that $\mathsf{GB}$, as a subset of $\R^3\times\R^3$, is contained in $[-K,K]^3\times [-K,K]^3$. Then $\mathbb{P}^{(-1,1),\vec{x},\vec{y}}_{\free}(\mathsf{BigFluc})\leq e^{-C^{-1}K^2}$ for some $C>0$ (see Lemma~\ref{lem:BB-osc}). Set 
\begin{equation}\label{def:GB_intro}
\mathsf{GB}:=\{(\vec{x},\vec{y})\, |\, Z^{(-1,1),\vec{x},\vec{y}}_{\mathbf{H}}\geq e^{-2^{-1}C^{-1}K^2},\ |x_i|\leq K, |y_i|\leq K \}.
\end{equation}
Then we have
\begin{align}\label{equ:middle_intro}
\mathbb{P}_{\cL}(\mathsf{BigFluc})\leq \mathbb{P}_{\cL}(\mathsf{GB}^{\textup{c}})+ e^{-2^{-1}C^{-1}K^2}. 
\end{align}
Therefore, the task of bounding $\mathbb{P}_{\cL}(\mathsf{BigFluc})$ reduces to find an upper bound of $\mathbb{P}_{\cL}(\mathsf{GB}^{\textup{c}})$.
 
The above argument applies for any $\mathbf{H}$-Brownian Gibbs line ensemble. Next, we specialize to the scaled KPZ line ensemble $\mathfrak{H}^t$. Strictly speaking, $\mathfrak{H}^t$ has infinite many random curves and the Gibbs property of $(\mathfrak{H}^t_1,\mathfrak{H}^t_2,\mathfrak{H}^t_3)$ is more involved than \eqref{equ:RN_intro}. We will ignore this issue because it suffices for the purpose of illustrating the main difficulties of analyzing the KPZ line ensemble and the main novelty of our arguments.

In view of the tail estimate of $\mathfrak{H}^t_1$ in \cite{CG1, CG2}, the optimal bound we can hope for is
\begin{equation}\label{equ:optbound}
 \mathbb{P}_{\mathfrak{H}^t}(\mathsf{GB}^{\textup{c}})\leq e^{-\tilde{C}^{-1}K^{3/2}} 
\end{equation}
for some $\tilde{C}>0$. Once we have \eqref{equ:optbound}, together with \eqref{equ:middle_intro},  we conclude that
\begin{align*}
\mathbb{P}_{\mathfrak{H}^t}(\mathsf{BigFluc})\leq  e^{-\tilde{C}^{-1}K^{3/2}}+ e^{-2^{-1}C^{-1}K^2}. 
\end{align*}
This provides local fluctuation estimates of $\mathfrak{H}^t$ and is the content of Theorem~\ref{thm:main1}. 

Deriving the bound \eqref{equ:optbound} for the scaled KPZ line ensembles is considerably more challenging compared to similar estimates for the parabolic Airy line ensemble. In the following, we explain this difficulty in more detail and discuss how we resolve such an issue. This resolution is the main technical contribution of the paper.

For each $t>0$, $\mathfrak{H}^t$ satisfies the $\mathbf{H}_t$-Brownian Gibbs property with $\mathbf{H}_t(x)=e^{t^{1/3}x}$. We aim to show \eqref{equ:optbound} holds uniformly in $t\geq 1$. When $t$ goes to infinity, the Boltzmann weights converge to the indicator function on the set $\mathsf{Ord}\subset C([1,3]_{\Z}\times [-1,1],\mathbb{R}) $, where
$$\mathsf{Ord}:=\{ \tilde{\cL}_1(u) >\tilde{\cL}_2(u) >\tilde{\cL}_3(u)\, \textup{for all}\ u\in [-1,1]  \}.$$  
This implies that if $\vec{x}$ and $\vec{y}$ are not strictly decreasing, then $Z^{(-1,1),\vec{x},\vec{y}}_{\mathbf{H}_t}$ converges to zero as $t$ goes to infinity. In particular, $(\vec{x},\vec{y})$ is not in $\mathsf{GB}$ for $t$ large enough. Therefore, we have to show that
\begin{equation}\label{equ:order_intro}
\mathfrak{H}^t(\pm 1)\ \textup{are strictly decreasing with high probability.}
\end{equation} 
We remark that for the parabolic Airy line ensemble, \eqref{equ:order_intro} follows directly from its definition.

To show \eqref{equ:order_intro}, we apply the Gibbs resampling on a larger interval $[-2,2]$. Intuitively, the $\mathbf{H}_t$-Brownian Gibbs property should guarantee \eqref{equ:order_intro} for $t$ large enough. The subtlety here is that we need to derive \eqref{equ:order_intro} \textbf{without} any control on the normalizing constant, because the very reason we seek for \eqref{equ:order_intro} is to estimate the normalizing constant.

\begin{figure}
\begin{center}
\begin{tabular}{lr}
\raisebox{0 ex}{\includegraphics[width=6cm]{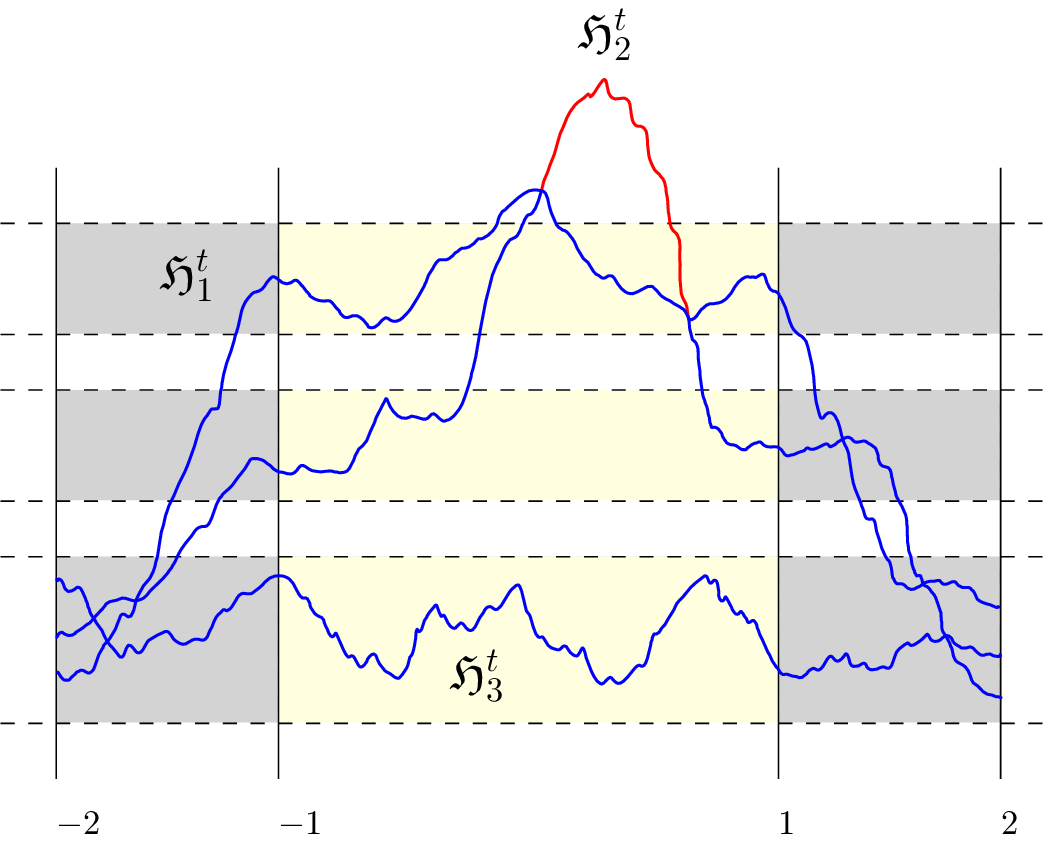}}
        &
\raisebox{0 ex}{\includegraphics[width=6cm]{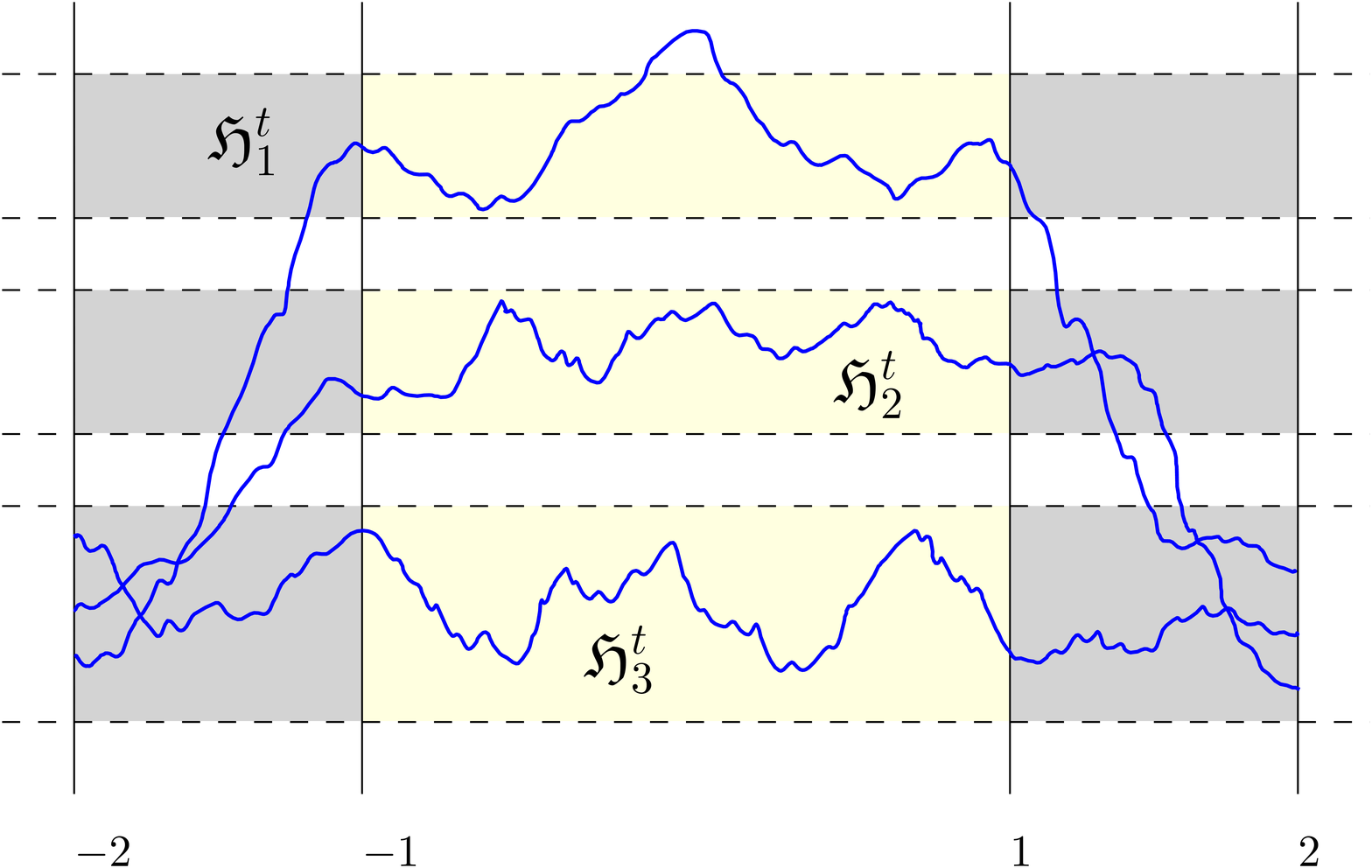}}
        \end{tabular}
\end{center}
\caption{An illustration of the two-step inductive argument to separate curves. In the first step, we ensure that $\mathfrak{H}^t_1>\mathfrak{H}^t_3$ and $\mathfrak{H}^t_2>\mathfrak{H}^t_3$ on the interval $[-1,1]$. However, $\mathfrak{H}^t_2\geq \mathfrak{H}^t_1$ may happen. See the picture on the left. In the second step, we push $\mathfrak{H}^t_3, \mathfrak{H}^t_2$ down to ensure $\mathfrak{H}^t_1>\mathfrak{H}^t_2>\mathfrak{H}^t_3$ on $[-1,1]$. See the picture on the right. }\label{fig:airy-KPZ}
\end{figure}

Now we are ready to explain the heuristics of our two-step inductive argument, introduced to resolve the difficulties caused by the intersecting nature of the KPZ line ensemble. To achieve \eqref{equ:order_intro}, we rely on the stochastic monotonicity, Lemma~\ref{monotonicity}. The stochastic monotonicity captures the idea that the interaction between $\mathfrak{H}^t_1$ and $\mathfrak{H}^t_2$ would push $\mathfrak{H}^t_1$ upward so $\mathfrak{H}^t_1$ is more likely to go up. This implies $\mathfrak{H}^t_1>\mathfrak{H}^t_2$ and $\mathfrak{H}^t_1>\mathfrak{H}^t_3$ on $[-1,1]$ are plausible to happen. However, because $\mathfrak{H}^t_2$ is interacting with $\mathfrak{H}^t_1$ and $\mathfrak{H}^t_3$ simultaneously, $\mathfrak{H}^t_2$ does not have a clear tendency. We overcome this difficulty through a novel inductive two-step resampling procedure. The idea is that if $\mathfrak{H}^t_1$ is very high, $\mathfrak{H}^t_2$ will almost not be affected by $\mathfrak{H}^t_1$. Hence $\mathfrak{H}^t_3$ can still push $\mathfrak{H}^t_2$ upward. We can lift $\mathfrak{H}^t_1$ first and then push $\mathfrak{H}^t_2$ up. After this step, both $\mathfrak{H}^t_1$ and $\mathfrak{H}^t_2$ are larger than $\mathfrak{H}^t_3$ on $[-1,1]$, but $\mathfrak{H}^t_2\geq \mathfrak{H}^t_1$ could occur. In the second step, we may apply a similar argument in the opposite direction to lower $\mathfrak{H}^t_2$ because $\mathfrak{H}^t_3$ is far below. At the end, we can achieve the ordering $\mathfrak{H}^t_1>\mathfrak{H}^t_2>\mathfrak{H}^t_3$ on $[-1,1]$. This in particular implies \eqref{equ:order_intro}. See Figure~\ref{fig:airy-KPZ} for an illustration.

\subsection{Outline} Section \ref{sec:basics of line ensembles} contains relevant definitions about Gibbsian line ensembles. Section \ref{sec:tools} contains lemmas about the stochastic monotonicity and the strong Markov property. Main results are stated in Section \ref{sec:main results} and are proved in Section \ref{sec:proofs}. The main technical result, Proposition \ref{pro:Z_lowerbound_k}, is proved in Section \ref{sec:Z}. We also record basic estimates on Brownian bridges in Appendix \ref{sec:Brownian} and deduce tail estimates for the scaled KPZ line ensembles in Appendix \ref{sec:tail}.
\subsection{Notations}
We would like to explain some notation here. The natural numbers are defined to be $\N = \{1, 2, . . .\}.$ For integers $k_1\leq k_2$, let $[k_1,k_2]_{\Z} := \{k_1,k_1+1,\ldots,k_2\}$. Events are denoted in a special font $\mathsf{E}$, their indicator function is written as $\mathbbm{1}_{\mathsf{E}}$  and their complement is written as $\mathsf{E}^c$.

Universal constants will be generally denoted by $C$ and constants depending only on $k$ will be denoted as $D=D(k)$. We label the ones in statements of theorems, lemmas and propositions with subscripts (e.g. $D_1(k), D_2(k), \cdots$ based on their order of occurrence) but the constants in proofs may change value from line to line.

\subsection{Acknowledgment} The author thanks Ivan Corwin for comments concerning a draft version of this paper. The author is grateful for the three Minerva lectures that Alan Hammond has given at Columbia in Spring 2019 on Gibbsian line ensembles.

\section{KPZ line ensemble and $\Ham$-Brownian Gibbs property}
In this section we define $\Ham$-Brownian Gibbsian line ensembles and discuss their properties in Section~\ref{sec:basics of line ensembles} and Section~\ref{sec:tools}. The Gibbs property for the KPZ line ensemble is formulated in Section~\ref{sec:KPZLEproperties}

\subsection{Line ensembles and the $\Ham$-Brownian Gibbs property}\label{sec:basics of line ensembles}
\begin{definition}[Line ensembles]\label{def:line-ensemble}
Let $\Sigma$ be an interval of $\Z$ and let $\Lambda$ be a subset of $\R$. Consider the set $C(\Sigma\times \Lambda,\R)$ of continuous functions from $\Sigma\times\Lambda$ to $\R$ endowed with the topology of uniform convergence on compact subsets of $\Sigma\times \Lambda$, and let $\mathcal{C}(\Sigma\times \Lambda,\R)$ denote the sigma-field generated by Borel sets in $C(\Sigma\times \Lambda,\R)$. A $\Sigma\times \Lambda$-indexed line ensemble $\mathcal{L}$ is a random variable on a probability space $(\Omega,\mathcal{B},\mathbb{P})$, taking values in $C(\Sigma\times \Lambda,\R)$ such that $\mathcal{L}$ is a measurable function from $\mathcal{B}$ to $\mathcal{C}(\Sigma\times \Lambda,\R)$.
\end{definition}

We think of such line ensembles as multi-layer random curves. We will generally write $\mathcal{L}:\Sigma\times \Lambda\to \R$ even though it is not $\mathcal{L}$, but rather $\mathcal{L}(\omega)$ for each $\omega\in \Omega$ which is such a function. We will also sometimes specify a line ensemble by only giving its law without reference to the underlying probability space.  We write $\mathcal{L}_i(\cdot):= \big(\mathcal{L}(\omega)\big)(i,\cdot)$ for the label $i\in \Sigma$ curve of the ensemble $\mathcal{L}$. 

\begin{definition}[Convergence of line ensembles]\label{def:weakconvergence} 
Given a $\Sigma\times \Lambda$-indexed line ensemble $\mathcal{L}$ and a sequence of such ensembles $\left\{\mathcal{L}^N\right\}_{N\geq 1}$, we will say that $\mathcal{L}^N$ converges to $\mathcal{L}$ weakly
as a line ensemble if for all bounded continuous functions $F:C(\Sigma\times \Lambda,\R)\to \R$, it holds that,
	\begin{equation*}
		\lim_{N\to\infty} \int F\big(\mathcal{L}^N(\omega)\big)d\mathbb{P}^N(\omega)= \int F\big(\mathcal{L}(\omega)\big) d\mathbb{P}(\omega).
	\end{equation*}
This is equivalent to weak-$*$ convergence in $C(\Sigma\times \Lambda,\R)$ endowed with the topology of uniform convergence on compact subsets of $\Sigma\times \Lambda$.
\end{definition}

We now start to formulate the $\Ham$-Brownian Gibbs property. We adopt the convention that all Brownian motions and bridges have diffusion parameter one.

\begin{definition}[$\Ham$-Brownian bridge line ensemble]\label{def:H_Brownian}
Fix $k_1\leq k_2$ with $k_1,k_2 \in \mathbb{Z}$, an interval $[a,b]\subset \mathbb{R}$ and two vectors $\vec{x},\vec{y}\in \mathbb{R}^{k_2-k_1+1}$. A $[k_1,k_2]_{\Z}\times [a,b]$-indexed line ensemble $\mathcal{L} = (\mathcal{L}_{k_1},\ldots,\mathcal{L}_{k_2})$ is called a free Brownian bridge line ensemble with entrance data $\vec{x}$ and exit data $\vec{y}$ if its law $\mathbb{P}^{k_1,k_2,[a,b],\vec{x},\vec{y}}_{\free}$ is that of $k_2-k_1+1$ independent standard Brownian bridges starting at time $a$ at the points $\vec{x}$ and ending at time $b$ at the points $\vec{y}$.

A Hamiltonian $\Ham$ is defined to be a measurable function $\Ham:\mathbb{R}\to [0,\infty]$. Throughout, we will make use of the special Hamiltonian
\begin{align}
\Ham_t(x) = e^{t^{1/3}x}.
\end{align} 
Fix a Hamiltonian $\Ham$. Let $k_1\leq k_2$ be two integers, $\vec{x},\vec{y}\in\mathbb{R}^{k_2-k_1-1}$, $f,g:(a,b)\to \mathbb{R}\cup\{\pm\infty\}$ be two measurable functions. The Boltzmann weight of $\mathcal{L}=(\mathcal{L}_{k_1},\mathcal{L}_{k_1+1},\dots, \mathcal{L}_{k_2})\in C([k_1,k_2]_{\Z}\times [a,b],\mathbb{R})$ is defined by 
\begin{equation}\label{def:Boltzmann_Brownian}
W_{\Ham}^{k_1,k_2,(a,b),\vec{x},\vec{y},f,g}(\mathcal{L}):= \exp\Bigg\{-\sum_{i=k_1-1}^{k_2}\int_a^b \Ham\Big(\mathcal{L}_{i+1}(u)-\mathcal{L}_{i}(u)\Big)du\Bigg\},
\end{equation}
where we adopt conventions that $\mathcal{L}_{k_1-1}=f$, $\mathcal{L}_{k_2+1}=g$. The normalizing constant is defined by
\begin{equation}\label{def:normalcont_Brownian}
Z_{\Ham}^{k_1,k_2,(a,b),\vec{x},\vec{y},f,g} :=\mathbb{E}^{k_1,k_2,(a,b),\vec{x},\vec{y}}_{\free}\Big[W_{\Ham}^{k_1,k_2,(a,b),\vec{x},\vec{y},f,g}(\mathcal{L})\Big],
\end{equation}
where $\mathcal{L}$ in the above expectation is distributed according to the measure $\mathbb{P}^{k_1,k_2,(a,b),\vec{x},\vec{y}}_{\free}$.

Suppose  $Z_{\Ham}^{k_1,k_2,(a,b),\vec{x},\vec{y},f,g}$ is positive. We define the $\Ham$-Brownian bridge line ensemble with entrance data $\vec{x}$, exit data $\vec{y}$ and boundary data $(f,g)$ to be a $[k_1,k_2]_{\Z}\times [a,b]$-indexed line ensemble $\mathcal{L} = (\mathcal{L}_{k_1},\ldots, \mathcal{L}_{k_2})$ with law $\mathbb{P}^{k_1,k_2,(a,b),\vec{x},\vec{y},f,g}_{\Ham}$ given according to the following Radon-Nikodym derivative relation:
 \begin{equation*}
 \frac{\textup{d}\mathbb{P}_{\Ham}^{k_1,k_2,(a,b),\vec{x},\vec{y},f,g}}{\textup{d}\mathbb{P}_{\free}^{k_1,k_2,(a,b),\vec{x},\vec{y}}}(\mathcal{L}) := \frac{W_{\Ham}^{k_1,k_2,(a,b),\vec{x},\vec{y},f,g}(\mathcal{L})}{Z_{\Ham}^{k_1,k_2,(a,b),\vec{x},\vec{y},f,g}}.
 \end{equation*}

Moreover, given $a'<b'$ contained in $(a,b)$, we define
\begin{equation}\label{def:Boltzmann_off}
W_{\Ham,(a',b')}^{k_1,k_2,(a,b),\vec{x},\vec{y},f,g}(\mathcal{L}):= \exp\Bigg\{-\sum_{i=k_1-1}^{k_2} \int_{[a,a']\cup [b',b]} \Ham\Big(\mathcal{L}_{i+1}(u)-\mathcal{L}_{i}(u)\Big)du\Bigg\}.
\end{equation}
That is, the weight is calculated on $[a,a']\cup [b',b]$ but not on $(a',b')$. We similarly define
\begin{equation}\label{def:normalcont_Brownian_t}
Z_{\Ham,(a',b')}^{k_1,k_2,(a,b),\vec{x},\vec{y},f,g} :=\mathbb{E}^{k_1,k_2,(a,b),\vec{x},\vec{y}}_{\free}\Big[W_{\Ham,(a',b')}^{k_1,k_2,(a,b),\vec{x},\vec{y},f,g}(\mathcal{L})\Big],
\end{equation}
and
\begin{equation*}
\frac{\textup{d}\mathbb{P}_{\Ham,(a',b')}^{k_1,k_2,(a,b),\vec{x},\vec{y},f,g}}{\textup{d}\mathbb{P}_{\free}^{k_1,k_2,(a,b),\vec{x},\vec{y}}}(\mathcal{L}) := \frac{W_{\Ham,(a',b')}^{k_1,k_2,(a,b),\vec{x},\vec{y},f,g}(\mathcal{L})}{Z_{\Ham,(a',b')}^{k_1,k_2,(a,b),\vec{x},\vec{y},f,g}}.
\end{equation*}
\end{definition}

$\Ham$-Brownian Gibbs property could be viewed as a spatial Markov property, more specifically, it provides a description of the conditional law inside a compact set.
\begin{definition}[$\Ham$-Brownian Gibbs property]\label{def:H-BGP}
Let $\Ham$ be a Hamiltonian, $\Sigma$ be an interval of $\Z$ and $\Lambda$ be a subset of $\R$. A $\Sigma \times \Lambda$-indexed line ensemble $\mathcal{L}$ satisfies the $\Ham$-Brownian Gibbs property if the following holds. For all $[k_1,k_2]_{\Z} \subset \Sigma$ and $[a,b]\subset \Lambda$, set $\vec{x}=\big(\mathcal{L}_{k_1}(a),\cdots ,\mathcal{L}_{k_2}(a)\big)$, $\vec{y}=\big(\mathcal{L}_{k_1}(b),\ldots ,\mathcal{L}_{k_2}(b)\big)$, $f=\mathcal{L}_{k_1-1}$ and $g=\mathcal{L}_{k_2+1}$ with the convention that if $k_1-1\notin\Sigma$ then $f\equiv +\infty$ and likewise if $k_2+1\notin \Sigma$ then $g\equiv -\infty$. Then it holds almost surely that $Z_{\Ham}^{k_1,k_2,(a,b),\vec{x},\vec{y},f,g}>0$. Moreover, 
\begin{equation}\label{def:conditioned}
\textrm{Law}\left(\mathcal{L} \left\vert_{[k_1,k_2]_{\Z} \times [a,b]}\ \textrm{conditional on } \mathcal{L}\right\vert_{(\Sigma \times \Lambda) \setminus ( [k_1,k_2]_{\Z} \times (a,b) )} \right) =\mathbb{P}_{\Ham}^{k_1,k_2,(a,b),\vec{x},\vec{y},f,g}
\end{equation}

The following equality \eqref{equ:expectation} is an equivalent formulation of \eqref{def:conditioned} and it is convenient for computations. 
For any bounded measurable function $F$ from $C\left([k_1,k_2]_{\mathbb{Z}}\times [a,b] , \R\right)$ to $\R$, it holds almost surely that
\begin{equation}\label{equ:expectation}
\E\left[F(\mathcal{L}\big|_{[k_1,k_2]_{\Z}\times [a,b]} ) \right\vert \Fext\big([k_1,k_2]_{\Z}\times (a,b) \big)\Big] = \mathbb{E}_{\Ham}^{k_1,k_2,(a,b),\vec{x},\vec{y},f,g}\left[F(\tilde{\mathcal{L}})\right],
\end{equation}
where \glossary{$\Fext\left(K\times (a,b)\right)$, Sigma-field generated by a line ensemble outside $[k_1,k_2]_{\Z}\times (a,b)$}
\begin{equation*}
\Fext\left([k_1,k_2]_{\Z}\times (a,b)\right) := \sigma\left(\mathcal{L}_{i}(s): (i,s)\in \Sigma\times \Lambda \setminus [k_1,k_2]_{\Z}\times (a,b)\right)
\end{equation*}
is the exterior sigma-field generated by the line ensemble outside $[k_1,k_2]_{\Z}\times (a,b)$. On the right-hand side $\tilde{\mathcal{L}}$ is distributed according to $\mathbb{P}_{\Ham}^{k_1,k_2,(a,b),\vec{x},\vec{y},f,g}$. 
\end{definition}

\begin{definition}\label{def:todo}
Consider the following special Hamiltonian, $$\mathbf{H}_{\textup{ord}}(x)=\left\{ \begin{array}{cc}
+\infty, & x>0,\\
0, & x\leq 0,
\end{array} \right. $$
Its corresponding $\mathbf{H}_{\textup{ord}}$-Brownian Gibbs property was also referred to as Brownian Gibbs property \cite{CH14}. The reason is that $\Ham_{\textup{ord}}$-Brownian bridge line ensemble (defined in Definition~\ref{def:H_Brownian}) is non-intersecting Brownian bridges conditioned on avoiding the upper and lower boundaries. We note that a Brownian Gibbsian line ensemble must be ordered with probability one.
\end{definition}

\subsection{Gibbs property of the KPZ line ensemble}\label{sec:KPZLEproperties}
For $t>0$, recall that the KPZ line ensemble $\mathcal{H}^t$ is a $\mathbb{N}\times\mathbb{R}$-indexed line ensemble defined in Definition~\ref{def:KPZlineensemble}. The Gibbs property of $\mathcal{H}^t$ is proved in \cite[Theorem 2.15]{CH16} and we record it as the theorem below.
\begin{theorem}\label{thm:KPZGibbs}
For any $t>0$, the KPZ line ensemble $\mathcal{H}^t$ satisfies the $\mathbf{H}$-Brownian Gibbs property with $\mathbf{H}(x)=e^x$.
\end{theorem} 
We are mostly interested in the scaled KPZ line ensemble $\mathfrak{H}^t$ defined in Definition~\ref{def:scaledKPZLE}. Recall $\Ham_t(x) = e^{t^{1/3}x}.$ Through a direct computation we have the following corollary.
\begin{corollary}\label{cor:KPZGibbs}
For any $t>0$, the line ensemble ${\mathfrak{H}}^t$ satisfies the $\mathbf{H}_t$-Brownian Gibbs property.
\end{corollary}

\subsection{Strong Gibbs property and stochastic monotonicity}\label{sec:tools}
We record some important properties, developed in \cite[Section 2]{CH16}, about $\Ham$-Gibbsian line ensembles in this section. We begin with the strong Gibbs property, which enables us to resample the trajectory within a stopping domain as opposed to a deterministic interval. 

\begin{definition}\label{def:stopdm}
Let $\Sigma$ be an interval of $\Z$, and $\Lambda$ be an interval of $\R$. Consider a $\Sigma\times\Lambda$-indexed line ensemble $\mathcal{L}$ which has the $\Ham$-Brownian Gibbs property for some Hamiltonian $\Ham$. For $[k_1,k_2]_{\mathbb{Z}} \subseteq \Sigma$ and $(\ell,r)\subseteq \Lambda$, $\Fext\big([k_1,k_2]_{\mathbb{Z}}\times (\ell,r)\big)$ denotes the sigma-field generated by the data outside $[k_1,k_2]_{\mathbb{Z}}\times (\ell,r)$. The random variable $(\mathfrak{l},\mathfrak{r})$ \glossary{$(\mathfrak{l},\mathfrak{r})$, Stopping domain} is called a {\it $[k_1,k_2]_{\mathbb{Z}}$-stopping domain} if for all $\ell<r$,
\begin{equation*}
\big\{\mathfrak{l} \leq \ell , \mathfrak{r}\geq r\big\} \in \Fext\big([k_1,k_2]_{\mathbb{Z}}\times (\ell,r)\big).
\end{equation*}
\end{definition}

For $[k_1,k_2]_{\mathbb{Z}}\subset \mathbb{Z}$, define
\begin{align*}
C^{k_1,k_2}:=\left\{ (\ell,r,f)\,:\, \ell<r,\ f=(f_{k_1},\dots,f_{k_2})\in C([k_1,k_2]_{\mathbb{Z}}\times [\ell,r]) \right\}.
\end{align*}

\begin{lemma}\label{lem:stronggibbs}
Consider a $\Sigma\times\Lambda$-indexed line ensemble $\mathcal{L}$ which has the $\Ham$-Brownian Gibbs property. Fix $[k_1,k_2]_{\mathbb{Z}}\subseteq \Sigma$. For all random variables $(\mathfrak{l},\mathfrak{r})$ which are $[k_1,k_2]_{\mathbb{Z}}$-stopping domains for $\mathcal{L}$, the following strong $\Ham$-Brownian Gibbs property holds: for all Borel functions $F: C^{k_1,k_2}\to\mathbb{R}$, $\bP$-almost surely,
\begin{equation}\label{eqn:stronggibbs}
\E\bigg[ F\Big(\mathfrak{l},\mathfrak{r}, \mathcal{L}\big\vert_{[k_1,k_2]_{\Z}\times [\mathfrak{l},\mathfrak{r}]}\Big) \Big\vert \Fext\big( [k_1,k_2]_{\Z	} \times (\mathfrak{l},\mathfrak{r})\big) \bigg]
=\mathbb{E}^{k_1,k_2,(\mathfrak{l},\mathfrak{r}),\vec{x},\vec{y},f,g}_{\mathbf{H}}\Big[F\big(\mathfrak{l},\mathfrak{r}, \tilde{\mathcal{L}} ) \Big],
\end{equation}
where $\vec{x} = \{\mathcal{L}_i(\mathfrak{l})\}_{i=k_1}^{k_2}$, $\vec{y} = \{\mathcal{L}_i(\mathfrak{r})\}_{i=k_1}^{k_2}$, $f =\mathcal{L}_{k_{1}-1}$ (or $\infty$ if $k_1-1\notin \Sigma$), $g=\mathcal{L}_{k_2+1}$ (or $-\infty$ if $k_2+1\notin \Sigma$). On the left-hand side $\mathcal{L}\big\vert_{[k_1,k_2]_{\Z}\times [\mathfrak{l},\mathfrak{r}]}$ is the restriction of curves distributed according to the law of $\mathcal{L}$ and on the right-hand side $\tilde{\mathcal{L}}$ is distributed according to $\bP^{k_1,k_2,(\mathfrak{l},\mathfrak{r}),\vec{x},\vec{y},f,g}_{\mathbf{H}}$.
\end{lemma}

For a convex Hamiltonian $\Ham$ (such as $\Ham_t$), $\Ham$-Brownian bridge line ensembles enjoy stochastic monotonicity.

\begin{lemma}\label{monotonicity}
Fix $k_1\leq k_2$ with $k_1,k_2\in \Z$ and $[a,b]\subset \mathbb{R}$. Consider two pairs of vectors $\vec{x}^{(i)},\vec{y}^{(i)}\in \R^{k_2-k_1+1}$ and two pairs of measurable functions $(f^{(i)},g^{(i)})$ for $i\in \{1,2\}$ such that $x^{(1)}_{j}\leq x^{(2)}_{j}$ and $y^{(1)}_{j}\leq y^{(2)}_{j}$ for all $j \in [k_1,k_2]_{\Z} $ and $f^{(1)}(u)\leq f^{(2)}(u), g^{(1)}(u)\leq g^{(2)}(u)$ for all $u\in (a,b)$. Suppose that $Z_{\mathbf{H}}^{k_1,k_2,(a,b),\vec{x}^{(i)},\vec{y}^{(i)},f^{(i)},g^{(i)}}>0$ for $i\in \{1,2\}$. Let $\mathcal{Q}^{(i)}=\{\mathcal{Q}^{(i)}_k, \dots, \mathcal{Q}^{(i)}_{k_2} \} $ be a $[k_1,k_2]_{\Z}\times [a,b]$-indexed line ensemble with law $\mathbb{P}_{\mathbf{H}}^{k_1,k_2,(a,b),\vec{x}^{(i)},\vec{y}^{(i)},f^{(i)},g^{(i)}}$. See Definition~\ref{def:H_Brownian}.

If the Hamiltonian $\Ham$ is convex, then there exists a coupling of $\mathcal{Q}^{(1)}$ and $  \mathcal{Q}^{(2)}$ such that almost surely $\mathcal{Q}^{(1)}_j(u)\leq \mathcal{Q}^{(2)}_j(u)$ for all $(j,u)\in [k_1,k_2]_{\Z}\times [a,b] $.
\end{lemma}

\section{Proof of main results}\label{sec:proofs}
In this section we prove Theorem~\ref{thm:main1} and Theorem~\ref{thm:main2}. The proofs rely heavily on the key technical result, Proposition~\ref{pro:Z_lowerbound_k}, which provides a quantitative lower bound on the normalization constant $Z$.  We state Proposition~\ref{pro:Z_lowerbound_k} below and postpone its proof to the next section.
\begin{proposition}\label{pro:Z_lowerbound_k}
Fix $k\in\mathbb{N} $ and $L\geq 1$. There exists a constant $D_1=D_1(k)$ depending only on $k$ such that the following statement holds. For all $t\geq 1$ and $K\geq L^2$, we have
\begin{align}\label{eq:Zsmall}
\mathbb{P} \left( Z^{1,k,(-L,L), \vec{x}, \vec{y},+\infty,\mathfrak{H}^t_{k+1}}_{\mathbf{H}_t}<  D_1^{-1}  e^{- D_1L^{-1}K^2}  \right)< e^{-K^{3/2}},
\end{align} 
where $\vec{x}=\big(\mathfrak{H}^t_i(-L)\big)_{i=1}^k$ and $\vec{y}=\big(\mathfrak{H}^t_i(L)\big)_{i=1}^k$.
\end{proposition}

\subsection{Proof of Theorem \ref{thm:main1} and Theorem \ref{thm:main2}(i)} We begin with proving Theorem \ref{thm:main1}.
\begin{proof}[Proof of Theorem \ref{thm:main1}]
Fix $L=1$ and $\varepsilon\in (0,e^{-1}]$. We will exploit Gibbs resampling on $[1,k]_{\mathbb{Z}}\times[-1,1]$. We start with controlling the boundary values $\mathfrak{H}^t_k(\pm1)$. Let $D'=\max\{D_9,D_{10}\}$ with $D_9,D_{10}$ from Proposition \ref{pro:lowertail} and Proposition \ref{pro:uppertail}. Set
$$r:=\left(D' \log\varepsilon^{-1}+D'\log(2D')\right)^{2/3}$$ and denote
$$\mathsf{E}_1:= \left\{ \sup_{u\in [-1,1]}\left| \mathfrak{H}^t_k(u)+2^{-1}u^2 \right| <r \right\}.$$ By the tail estimates in Proposition \ref{pro:lowertail} and Proposition \ref{pro:uppertail}, we have $\mathbb{P}(\mathsf{E}_1^c)<\varepsilon$.\\[0.1cm]
\indent Recall the convention that $\vec{x}=\big(\mathfrak{H}^t_i(-L)\big)_{i=1}^k$ and $\vec{y}=\big(\mathfrak{H}^t_i(L)\big)_{i=1}^k$ and denote \begin{align*}
\mathsf{E}_2:=\left\{   Z^{1,k,(-1,1),  \vec{x},  \vec{y},+\infty ,\mathfrak{H}^t_{k+1}}_{\mathbf{H}_t}\geq \delta    \right\}.
\end{align*}
Here $\delta$ is a small constant to be determined soon. Let $D_1$ be the constant in Proposition \ref{pro:Z_lowerbound_k}. Apply Proposition \ref{pro:Z_lowerbound_k} with $L=1$ and the right hand side being $\varepsilon\in (0,e^{-1}]$. It is easy to check that $\mathbb{P}(\mathsf{E}_2^c)<\varepsilon$ by picking
\begin{align}\label{eq:delta}
\delta:=&D_1^{-1}e^{-D_1(\log\varepsilon^{-1})^{4/3}}.
\end{align}

Let $B:[-1,1]\to\mathbb{R}$ be a free Brownian bridge with $B(\pm 1)=0$ and write $\mathbb{P}_{\free}$ to represent the law of $B$. Let $C_0$ be the universal constant in Lemma \ref{lem:BB-osc}. For all $K\geq 0$ and all $d\in (0,1]$, 
\begin{align*}
\mathbb{P}_{\textup{free}}\left( \sup_{u,v\in [0,d]}\left| B(u)-B(v) \right|\geq Kd^{1/2} \right)\leq C_0 e^{-C_0^{-1} K^2}.
\end{align*}
Define $K=K(\varepsilon)$ as
\begin{align}\label{def:K}
K(\varepsilon):=\max\{2 r, (8C_0D_1)^{1/2}\left(\log\varepsilon^{-1} \right)^{2/3}  \}.
\end{align}
Let 
\begin{align*}
\mathfrak{H}^{t,[-1,1]}_k(u):=\mathfrak{H}^{t}_k(u)-\left(2^{-1}(1-u)\mathfrak{H}^{t }_k(-1)+2^{-1}(u-1)\mathfrak{H}^{t }_k(1)\right).
\end{align*}
Define the events
\begin{align*}
\mathsf{G}\coloneqq &\left\{\sup_{u,v\in [0,d]}\left|\mathfrak{H}_k^t(u)-\mathfrak{H}_k^t(v) \right|\geq Kd^{1/2}\right\}\\ 
\mathsf{G}'\coloneqq&\left\{\sup_{u,v\in [0,d]}\left|\mathfrak{H}_k^{t,[-1,1]}(u)-\mathfrak{H}_k^{t,[-1,1]}(v) \right|\geq 2^{-1}Kd^{1/2}\right\} . 
\end{align*}
Note that
\begin{align*}
\sup_{u,v\in [0,d]}\left|\mathfrak{H}_k^t(u)-\mathfrak{H}_k^t(v) \right|\leq \sup_{u,v\in [0,d]}\left|\mathfrak{H}_k^{t,[-1,1]}(u)-\mathfrak{H}_k^{t,[-1,1]}(v) \right|+2^{-1}d|\mathfrak{H}_k^t(1)-\mathfrak{H}_k^t(-1)|. 
\end{align*}
When $\mathsf{E}_1$, we have
, $$2^{-1}d|\mathfrak{H}_k^t(1)-\mathfrak{H}_k^t(-1)|\leq rd\leq rd^{1/2} .$$ 
Here we used $d\leq 1$. Because $K\geq 2r$, we deduce that 
$$2^{-1}d|\mathfrak{H}_k^t(1)-\mathfrak{H}_k^t(-1)|\leq 2^{-1}Kd^{1/2}.$$ Therefore we have $\mathsf{G} \cap \mathsf{E}_1\subset 
\mathsf{G}'\cap \mathsf{E}_1$. Applying the Gibbs property (see Definition \ref{def:H-BGP} and Corollary~\ref{cor:KPZGibbs}), we obtain that
\begin{align*}
  \mathbbm{1}_{\mathsf{E}_1\cap\mathsf{E}_2} \cdot \mathbb{E}[\mathbbm{1}_{\mathsf{G}'}\,|\,\mathcal{F}_{ {ext}}([1,k]_{\mathbb{Z}}\times (-1,1))] =  \mathbbm{1}_{\mathsf{E}_1\cap\mathsf{E}_2}\cdot \frac{\mathbb{E}^{1,k,(-1,1),\vec{x},\vec{y} }_{\free }[\mathbbm{1}_{\mathsf{G}'}\cdot W^{1,k,(-1,1),\vec{x},\vec{y},+\infty,\mathfrak{H}_{k+1}^t}_{\mathbf{H}_t} ]}{ Z^{1,k,(-1,1),\vec{x},\vec{y},+\infty,\mathfrak{H}_{k+1}^t}_{\mathbf{H}_t}}. 
\end{align*}
As $\mathsf{E}_2$ occurs, $ Z^{1,k,(-1,1),\vec{x},\vec{y},+\infty,\mathfrak{H}_{k+1}^t}_{\mathbf{H}_t}\geq \delta$. Together with $W^{1,k,(-1,1),\vec{x},\vec{y},+\infty,\mathfrak{H}_{k+1}^t}_{\mathbf{H}_t}\leq 1$, we deduce
\begin{align*}
   \mathbbm{1}_{\mathsf{E}_1\cap\mathsf{E}_2}\cdot \mathbb{E}[\mathbbm{1}_{\mathsf{G}'}\,|\,\mathcal{F}_{ {ext}}([1,k]_{\mathbb{Z}}\times (-1,1))]  & \leq \delta^{-1}\mathbb{P}^{1,k,(-1,1),\vec{x},\vec{y} }_{\free }(  \mathsf{G}' )\cdot  \mathbbm{1}_{\mathsf{E}_1\cap\mathsf{E}_2} \\
    &\leq \delta^{-1}C_0e^{-4^{-1}C_0^{-1}K^2}\cdot  \mathbbm{1}_{\mathsf{E}_1\cap\mathsf{E}_2}.
\end{align*}
In summary, we obtain that
\begin{align*}
\mathbb{P}(\mathsf{G})&\leq \mathbb{P}(\mathsf{G}\cap\mathsf{E}_1\cap \mathsf{E}_2)+\mathbb{P}(\mathsf{E}_1^{\textup{c}})+\mathbb{P}(\mathsf{E}_2^{\textup{c}})\leq \mathbb{P}(\mathsf{G}'\cap\mathsf{E}_1\cap \mathsf{E}_2)+\mathbb{P}(\mathsf{E}_1^{\textup{c}})+\mathbb{P}(\mathsf{E}_2^{\textup{c}})\\
&=\mathbb{E}[\mathbbm{1}_{\mathsf{E}_1\cap\mathsf{E}_2}\cdot \mathbb{E}[\mathbbm{1}_{\mathsf{G}'}\,|\,\mathcal{F}_{ {ext}}([1,k]_{\mathbb{Z}}\times (-1,1))]]+\mathbb{P}(\mathsf{E}_1^{\textup{c}})+\mathbb{P}(\mathsf{E}_2^{\textup{c}})\\
&\leq \delta^{-1}C_0e^{-4^{-1}C_0^{-1}K^2}+2\varepsilon. 
\end{align*}
Recall definitions of $\delta$ in \eqref{eq:delta} and $K$ in \eqref{def:K}, it is easily checked that $\delta^{-1}\leq D_1e^{8^{-1}C_0^{-1}K^2}$. Moreover there exists $D=D(k)$ such that $K\leq D^{2/3}\left(\log \varepsilon^{-1}\right)^{2/3}$. Hence we obtain that
\begin{align*}
\mathbb{P}\left(\mathsf{G} \right)&\leq   C_0D_1 e^{-8^{-1}C_0^{-1}K^2 }+2\varepsilon\\
&\leq   C_0D_1 e^{-8^{-1}C_0^{-1}K^2 }+ 2e^{-D^{-1}K^{3/2}}\\
&\leq   C_0D_1 e^{-8^{-1}C_0^{-1}K^{3/2}}+ 2e^{-D^{-1}K^{3/2}}.
\end{align*}

Picking $D_0(k)=2\max\left\{C_0D_1, 8C_0, D,2\right\}$, we have
\begin{align*}
\mathbb{P}(\mathsf{G} )&\leq   \frac{1}{2}D_0 e^{-D_0^{-1}K^{3/2}}+ \frac{1}{2}D_0e^{-D_0^{-1}K^{ 3/2}}\leq  D_0e^{-D_0^{-1}K^{ 3/2}}.
\end{align*}
The proof is finished.
\end{proof}

Using an union bound argument, we obtain the following corollary. 
\begin{corollary}\label{cor:osc}
Fix $k\in\mathbb{N}$, $L\geq 1$ and $d\in (0,1]$. There exists a constant $ D_2=D_2(k)$ depending only on $k$ such that the following statement holds. Suppose that $K\geq 4Ld^{1/2}$, then for all $t\geq 1$, 
\begin{align*}
\mathbb{P}\left( \sup_{u,v\in [-L,L],\ |u-v|\leq d}\left|\mathfrak{H}_k^t(u)-\mathfrak{H}_k^t(v) \right|\geq K d^{1/2} \right)\leq  d^{-1}L D_2 e^{-D_2^{-1} K^{3/2} }.
\end{align*} 
\end{corollary}

\begin{proof}
Because
\begin{align*}
\left\{\sup_{u,v\in [-L,L],\ |u-v|\leq d}\left|\mathfrak{H}_k^t(u)-\mathfrak{H}_k^t(v) \right|\geq K d^{1/2}\right\}\subset \left\{\sup_{u,v\in [-L,L],\ |u-v|\leq 2^{-1}d }\left|\mathfrak{H}_k^t(u)-\mathfrak{H}_k^t(v) \right|\geq 2^{-1}K d^{1/2}\right\},
\end{align*}
it suffices to prove the assertion for $d\leq 2^{-1}$. From now on we assume $d\in (0,2^{-1}]$. 

Let $m=\lceil d^{-1}L  \rceil$. For $j\in [-m,m]_{\mathbb{Z}}$, define
\begin{align*}
\mathsf{E}_j:=\left\{ \left|\mathfrak{H}_k^t(u)-\mathfrak{H}_k^t(v) \right|\geq K d^{1/2}\,\, \textup{for some}\ u,v\in [jd,(j+2) d]\right\}.
\end{align*}
Note that
\begin{align*}
\left\{ \sup_{u,v\in [-L,L],\ |u-v|\leq d}\left|\mathfrak{H}_k^t(u)-\mathfrak{H}_k^t(v) \right|\geq K d^{1/2} \right\}\subset \bigcup_{j={-m}}^m\mathsf{E}_j.
\end{align*}

Fix $j\in [-m,m]_\mathbb{Z}$. By stationarity of $\mathfrak{H}^t_k(u)+2^{-1}u^2$, $\mathfrak{H}^t_k(u)\overset{(d)}{=} \mathfrak{H}^t_k(u-T)-T(u-T)-2^{-1}T^2$ as a process in $u$ for any $T\in\R$. Therefore, 
\begin{align*}
\mathbb{P}(\mathsf{E}_j)=&\mathbb{P}\left( \sup_{u,v\in [0,2d]}\left| \mathfrak{H}^t_k(u)-\mathfrak{H}^t_k(v)-jd(u-v)  \right|\geq Kd^{1/2} \right)\\
\leq &\mathbb{P}\left( \sup_{u,v\in [0,2d]}\left| \mathfrak{H}^t_k(u)-\mathfrak{H}^t_k(v) \right|\geq Kd^{1/2}-d(L+1) \right)\\
\leq &\mathbb{P}\left( \sup_{u,v\in [0,2d]}\left| \mathfrak{H}^t_k(u)-\mathfrak{H}^t_k(v) \right|\geq Kd^{1/2}-2d L  \right)\\
\leq &\mathbb{P}\left( \sup_{u,v\in [0,2d]}\left| \mathfrak{H}^t_k(u)-\mathfrak{H}^t_k(v) \right|\geq 2^{-1} Kd^{1/2} \right).
\end{align*}
We have used $ |jd(u-v)|  \leq d(L+1)$, $L\geq 1$ and $K\geq  4Ld^{1/2} $. \\[0.1cm]
\indent From Theorem \ref{thm:main1}, there exists $D=D(k)$ such that for all $j\in [-m,m]_{\mathbb{Z}}$, we have $\mathbb{P}(\mathbb{\mathsf{E}}_j)\leq D e^{-D^{-1} K^{3/2} }.$ In conclusion,
\begin{align*}
\mathbb{P}\left( \sup_{u,v\in [-L,L],\ |u-v|\leq d}\left|\mathfrak{H}_k^t(u)-\mathfrak{H}_k^t(v) \right|\geq K d^{1/2} \right)\leq \sum_{j=-m}^m \mathbb{P}(\mathbb{\mathsf{E}}_j)\leq 3d^{-1}LD e^{-D^{-1} K^{3/2} }.
\end{align*}
The proof is finished by taking $D_2=3D$.
\end{proof}

\begin{proof}[Proof of Theorem \ref{thm:main2}(i).]
Let $t_N$, $N\in\mathbb{N}$ be an increasing sequence which diverges to infinity. Without loss of generality, we assume $t_N\geq 1$. Given $\varepsilon>0$, we need to find a compact subset $S\subset C(\mathbb{N}\times \mathbb{R},\mathbb{R})$ (with respect to the topology of uniform convergence on compact subsets in Definition \ref{sec:basics of line ensembles}) such that for all $N$, we have $$\mathbb{P}(\mathfrak{H}^{t_N}\in S)\geq 1-\varepsilon.$$

For $i\in\mathbb{N}$, let $r_i>0$ and $ d_{i,\ell}>0$ with $\ell\in \mathbb{N}$ be positive numbers to be determined soon. Consider closed subsets of $C(\mathbb{N}\times \mathbb{R},\mathbb{R})$ as follows,
\begin{align*}
S_{i,0}\coloneqq&\left\{ \sup_{1\leq j\leq i}|\mathcal{L}_j(0)|\leq r_i \right\},\\
S_{i,\ell}\coloneqq&\left\{  \sup_{1\leq j\leq i}\sup_{u,v\in [-i,i], |u-v|\leq d_{i,\ell}}|\mathcal{L}_j(u)-\mathcal{L}_j(v)|\leq \ell^{-1} \right\},\\
S_{i}\coloneqq&S_{i,0}\cap \bigcap_{\ell=1}^\infty S_{i,\ell}.
\end{align*}

It suffices to pick constants $r_i$ and $d_{i,\ell}$ such that for all $t\geq 1$, we have
\begin{align}\label{h_in_K}
\mathbb{P}(\mathfrak{H}^t\notin S_i )\leq  2^{-i}\varepsilon.
\end{align}
Suppose for a moment \eqref{h_in_K} holds. Let $S\coloneqq \cap_{i=1}^\infty S_i$. From the Arzel\`a–Ascoli theorem and a diagonal argument, any sequence in $S$ has a convergent subsequence. Together with the closeness of $S$, $S$ is sequentially compact. Because the topology of uniform convergence on compact subsets on $C(\mathbb{N}\times\mathbb{R},\mathbb{R})$ is metrizable, $S$ is a compact subset of $C(\mathbb{N}\times \mathbb{R},\mathbb{R})$. Furthermore, \eqref{h_in_K} implies $\mathbb{P}(\mathfrak{H}^t\in S)\geq 1-\varepsilon $ for all $t\geq 1$.\\

It remains to prove \eqref{h_in_K}. We now explain how $r_i$ and $d_{i,\ell}$ are chosen, based on Propositions \ref{pro:lowertail} and \ref{pro:uppertail} and Corollary \ref{cor:osc}. Let $D_9(k)$ and $D_{10}(k)$ be the constants in Propositions \ref{pro:lowertail} and \ref{pro:uppertail}. Then Propositions \ref{pro:lowertail} and \ref{pro:uppertail} implies that for all $t\geq 1$ we have
\begin{align*}
\mathbb{P}(\mathfrak{H}^t\notin S_{i,0})\leq \sum_{j=1}^i \left( D_{10}(j)e^{-D_{10}(j)^{-1}r_i^{3/2}}+D_9(j)e^{-D_9(j)^{-1}r_i^{3/2}} \right).
\end{align*} 
Similarly, let $D_2(k)$ be the constant in Corollary \ref{cor:osc}. Then for all $t\geq 1$ we have
\begin{align*}
\mathbb{P}(\mathfrak{H}^t\notin S_{i,\ell})\leq \sum_{j=1}^i    id_{i,\ell}^{-1}D_2(j)e^{-D_2(j)^{-1}d_{i,\ell}^{-3/4} \ell^{3/2}}
\end{align*} 
provided $\ell^{-1} d_{i,\ell}^{-1/2}\geq 4i d_{i,\ell}^{1/2}$. By choosing $r_i$ large enough and $d_{i,\ell}$ small enough, we have $\mathbb{P}(\mathfrak{H}^t\notin S_{i,0})\leq 2^{-i-1}\varepsilon$ and $\mathbb{P}(\mathfrak{H}^t\notin S_{i,\ell})\leq 2^{-i-\ell-1}\varepsilon$. Then \eqref{h_in_K} follows.\\[0.1cm]

\end{proof}

\subsection{Proof of Theorem \ref{thm:main2}(ii)}
The purpose of this section is to show any subsequential limit of $\mathfrak{H}^t$ is non-intersecting and satisfies the $\mathbf{H}_{\textup{ord}}$-Brownian Gibbs property defined in Definition~\ref{def:todo}. We start by showing that $\mathfrak{H}^t$ is asymptotically strictly ordered.

\begin{proposition}\label{pro:ordered}
For any $k\in\mathbb{N}$ and $\varepsilon\in (0,1]$, there exist $\rho =\rho (\varepsilon,k)>0$ and $t_0=t_0(\varepsilon ,k)\geq 1$ such that for all $t\geq t_0$
\begin{align*}
\mathbb{P}\bigg( \inf_{u\in [-1,1]}(\mathfrak{H}^t_k(u)-\mathfrak{H}^t_{k+1}(u))< \rho \bigg)\leq   \varepsilon.
\end{align*}
\end{proposition}

\begin{proof}
For any $\rho,d>0$ consider the events
\begin{align*}
\mathsf{D}_1:=&\left\{\inf_{u\in [-1,1]}(\mathfrak{H}^t_k(u)-\mathfrak{H}^t_{k+1}(u))\in (-\rho,\rho)\right\},\\
\mathsf{D}_2:=&\left\{\inf_{u\in [-1,1]}(\mathfrak{H}^t_k(u)-\mathfrak{H}^t_{k+1}(u))\leq -\rho \right\},
\end{align*}
and
\begin{align*}
\mathsf{D}_3\coloneqq \left\{ \sup_{u,v\in [-1,1],|u-v|\leq d} |(\mathfrak{H}^t_k(u)-\mathfrak{H}^t_{k+1}(u))-(\mathfrak{H}^t_k(v)-\mathfrak{H}^t_{k+1}(v))| <  2^{-1}{\rho} \right\}.
\end{align*}  

For two positive numbers $M$ and $\delta$ to be determined soon, define the event
\begin{align*}
\mathsf{G}:=\bigg\{ Z^{1,k+1,(-1,1),\vec{x},\vec{y},+\infty,\mathfrak{H}^t_{k+2}}_{\mathbf{H}_t}\geq \delta \bigg\}\cap\bigg\{ |\mathfrak{H}^t_{k}(\pm 1 )|,|\mathfrak{H}^t_{k+1}(\pm 1)| \leq M \bigg\},
\end{align*}
with $\vec{x}=(\mathfrak{H}^{t}_j(-1))_{j=1}^{k+1}$, $\vec{y}=(\mathfrak{H}^{t}_j( 1))_{j=1}^{k+1} $. By Proposition \ref{pro:Z_lowerbound_k}, Proposition \ref{pro:lowertail} and Proposition \ref{pro:uppertail}, there exist $\delta$ and $M$, depending only on $\varepsilon$ and $k$ such that $\mathbb{P}(\mathsf{G})\geq 1-\varepsilon$. 
\begin{claim}
There exists $\rho >0$, depending only on $\varepsilon\delta$ and $M$, such that
\begin{align*}
\mathbbm{1}_{\mathsf{G}}\cdot \mathbb{P}^{1,k+1,(-1,1),\vec{x},\vec{y}}_{\free}(\mathsf{D}_1)<\varepsilon \delta\cdot \mathbbm{1}_{\mathsf{G}} .
\end{align*}
Here $\vec{x}=(\mathfrak{H}^t_i(-1))_{i=1}^{k+1}$ and $\vec{y}=(\mathfrak{H}^t_i(1))_{i=1}^{k+1}$. 
\end{claim}
\begin{proof}
Under the law of $\mathbb{P}^{1,k+1,(-1,1),\vec{x},\vec{y}}_{\free}$, $2^{-1}(B_k(2x-1)-B_{k+1}(2x-1))$ is a standard Brownian bridge defined on $[0,1]$ and has boundary values $2^{-1}(\mathfrak{H}^t_k(-1)-\mathfrak{H}^t_{k+1}(-1))=:2^{-1}a$ and $2^{-1}(\mathfrak{H}^t_k(1)-\mathfrak{H}^t_{k+1})(1)=:2^{-1}b$. From Lemma \ref{lem:BB-max}, we have  
\begin{align*}
\mathbb{P}^{1,k+1,(-1,1),\vec{x},\vec{y}}_{\free}(\mathsf{D}_1)= e^{-2^{-1}(a-\rho_+ )(b-\rho_+ )}-e^{-2^{-1}(a-\rho_- )(b-\rho_- )}.
\end{align*}
Here $\rho_+=\min\{\rho,a,b\}$ and $\rho_-=\min\{-\rho,a,b\}$. $\mathsf{G}$ implies $a,b\in [-2M,2M]$. Therefore for $\rho>0$ small enough depending on $\delta\varepsilon$ and $M$, we have
\begin{align*}
\mathbbm{1}_{\mathsf{G}}\cdot \mathbb{P}^{1,k+1,(-1,1),\vec{x},\vec{y}}_{\free}(\mathsf{D}_1)<\varepsilon \delta\cdot \mathbbm{1}_{\mathsf{G}} .
\end{align*}
\end{proof}
From now on we fix such $\rho>0$. Applying the Gibbs Property (see Definition \ref{def:H-BGP} and Corollary~\ref{cor:KPZGibbs}), we deduce
\begin{align*}
\mathbb{P}( \mathsf{D}_1\cap \mathsf{G} )=&\mathbb{E}[\mathbbm{1}_{\mathsf{G}}\cdot \mathbb{E}[\mathbbm{1}_{\mathsf{D}_1}\,|\, \mathcal{F}_{ {ext}}([1,k+1]_{\mathbb{Z}}\times (-1,1))]]\\
=&\mathbb{E}\bigg[\mathbbm{1}_{\mathsf{G}}\cdot \frac{\mathbb{E}^{1,k+1,(-1,1),\vec{x},\vec{y}}_{\free}[\mathbbm{1}_{\mathsf{D}_1}\cdot W^{1,k+1,(-1,1),\vec{x},\vec{y},+\infty,\mathfrak{H}^t_{k+2}}_{\mathbf{H}_t} ]}{Z^{1,k+1,(-1,1),\vec{x},\vec{y},+\infty,\mathfrak{H}^t_{k+2}}_{\mathbf{H}_t}  } \bigg]\\
\leq & \mathbb{E}\left[\mathbbm{1}_{\mathsf{G}}\cdot \delta^{-1} \mathbb{P}^{1,k+1,(-1,1),\vec{x},\vec{y}}_{\free}( {\mathsf{D}_1}  )  \right]\leq \varepsilon.
\end{align*}
\begin{claim}
There exists $d\in (0,1]$, depending only on $\varepsilon \delta$, $M$ and $\rho$, such that 
\begin{align*}
\mathbbm{1}_{\mathsf{G}}\cdot \mathbb{P}^{1,k+1,(-1,1),\vec{x},\vec{y}}_{\free}  \left(\mathsf{D}_3^{\textup{c}}  \right)<\varepsilon \delta\cdot \mathbbm{1}_{\mathsf{G}}.
\end{align*} 
Here $\vec{x}=(\mathfrak{H}^t_i(-1))_{i=1}^{k+1}$ and $\vec{y}=(\mathfrak{H}^t_i(1))_{i=1}^{k+1}$. 
\end{claim}
\begin{proof}
Under the law of $\mathbb{P}^{1,k+1,(-1,1),\vec{x},\vec{y}}_{\free}$, $2^{-1/2}(B_k(u)-B_{k+1}(u))$ is a standard Brownian bridge defined on $[-1,1]$ and has boundary values $2^{-1/2}(\mathfrak{H}^t_k(-1)-\mathfrak{H}^t_{k+1}(-1))=:2^{-1/2}a$ and $2^{-1/2}(\mathfrak{H}^t_k(1)-\mathfrak{H}^t_{k+1})(1)=:2^{-1}b$. 
\begin{align*}
 \mathbb{P}^{1,k+1,(-1,1),\vec{x},\vec{y}}_{\free}  \left( \mathsf{D}_3^{\textup{c}} \right)&=\mathbb{P}^{1,1,(-1,1),2^{-1/2}a,2^{-1/2}b}_{\free}\left(\sup_{u,v\in [-1,1],|u-v|\leq d} |B(u)-B(v)| \geq  2^{-3/2}{\rho} \right).\\
  &\leq\mathbb{P}^{1,1,(-1,1),0,0}_{\free}\left(\sup_{u,v\in [-1,1],|u-v|\leq d} |B(u)-B(v)| \geq  2^{-3/2}{\rho}-2^{-3/2}d|b-a| \right).
\end{align*} 
As $\mathsf{G}$ occurs, $|b-a|\leq 4M$. Hence by taking $d\leq 8^{-1}\rho M$,
\begin{align*}
 &\mathbbm{1}_{\mathsf{G}}\cdot \mathbb{P}^{1,k+1,(-1,1),\vec{x},\vec{y}}_{\free}  \left( \mathsf{D}_3^{\textup{c}} \right)\\
 \leq &\mathbbm{1}_{\mathsf{G}}\cdot\mathbb{P}^{1,1,(-1,1),0,0}_{\free}\left(\sup_{u,v\in [-1,1],|u-v|\leq d} |B(u)-B(v)| \geq  2^{-3/2}{\rho}-2^{-3/2}d|b-a| \right)\\
  \leq &\mathbbm{1}_{\mathsf{G}}\cdot\mathbb{P}^{1,1,(-1,1),0,0}_{\free}\left(\sup_{u,v\in [-1,1],|u-v|\leq d} |B(u)-B(v)| \geq  2^{-3/2}({\rho}-4dM)  \right)\\
    \leq &\mathbbm{1}_{\mathsf{G}}\cdot\mathbb{P}^{1,1,(-1,1),0,0}_{\free}\left(\sup_{u,v\in [-1,1],|u-v|\leq d} |B(u)-B(v)| \geq  2^{-5/2} {\rho}   \right).
\end{align*}
Applying Lemma \ref{lem:BB-osc}, we conclude that for $d$ small enough,
\begin{align*}
  \mathbbm{1}_{\mathsf{G}}\cdot \mathbb{P}^{1,k+1,(-1,1),\vec{x},\vec{y}}_{\free}  \left(\mathsf{D}_3^{\textup{c}} \right)  \leq &\mathbbm{1}_{\mathsf{G}}\cdot Ce^{-C^{-1}d^{-1}\rho^2}\leq \mathbbm{1}_{\mathsf{G}}\cdot \varepsilon\delta.
\end{align*}
\end{proof}
From now on we fix such $d$. Applying the Gibbs Property, we deduce
\begin{align*}
\mathbb{P}( \mathsf{D}^{\textup{c}}_3\cap \mathsf{G} )=&\mathbb{E}[\mathbbm{1}_{\mathsf{G}}\cdot \mathbb{E}[\mathbbm{1}_{\mathsf{D}^{\textup{c}}_3}\,|\, \mathcal{F}_{ {ext}}([1,k+1]_{\mathbb{Z}}\times (-1,1))]]\\
=&\mathbb{E}\bigg[\mathbbm{1}_{\mathsf{G}}\cdot \frac{\mathbb{E}^{1,k+1,(-1,1),\vec{x},\vec{y}}_{\free}[\mathbbm{1}_{\mathsf{D}^{\textup{c}}_3 }\cdot W^{1,k+1,(-1,1),\vec{x},\vec{y},+\infty,\mathfrak{H}^t_{k+2}}_{\mathbf{H}_t} ]}{Z^{1,k+1,(-1,1),\vec{x},\vec{y},+\infty,\mathfrak{H}^t_{k+2}}_{\mathbf{H}_t}  } \bigg]\\
\leq & \mathbb{E}\left[\mathbbm{1}_{\mathsf{G}}\cdot \delta^{-1} \mathbb{P}^{1,k+1,(-1,1),\vec{x},\vec{y}}_{\free}(  \mathsf{D}^{\textup{c}}_3  )  \right]\leq \varepsilon.
\end{align*}

As $\mathsf{D}_2\cap \mathsf{D}_3$ occurs, there exists an interval of $ [-1,1]$ with length $d$ in which  $B_1(u)-B_2(u)\leq 2^{-1}\rho$. As a consequence
\begin{align*}
W^{1,k+1,(-1,1),\vec{x},\vec{y},+\infty,\mathfrak{H}^t_{k+2}}_{\mathbf{H}_t}\cdot\mathds{1}_{\mathsf{D}_2\cap\mathsf{D}_3}\leq \exp\left( -d e^{ 2^{-1}\rho t^{1/3}}\right)\cdot\mathds{1}_{\mathsf{D}_2\cap\mathsf{D}_3}.
\end{align*}
Now we can take $t$ large such that $\exp\left( -d e^{ 2^{-1}\rho t^{1/3}}\right)\leq \varepsilon \delta$. Applying the Gibbs Property, we deduce
\begin{align*}
\mathbb{P}( \mathsf{D}_2\cap\mathsf{D}_3 \cap \mathsf{G} )=&\mathbb{E}[\mathbbm{1}_{\mathsf{G}}\cdot \mathbb{E}[\mathbbm{1}_{\mathsf{D}_2\cap\mathsf{D}_3 }\,|\, \mathcal{F}_{ {ext}}([1,k+1]_{\mathbb{Z}}\times (-1,1))]]\\
=&\mathbb{E}\bigg[\mathbbm{1}_{\mathsf{G}}\cdot \frac{\mathbb{E}^{1,k+1,(-1,1),\vec{x},\vec{y}}_{\free}[\mathbbm{1}_{\mathsf{D}_2\cap\mathsf{D}_3 }\cdot W^{1,k+1,(-1,1),\vec{x},\vec{y},+\infty,\mathfrak{H}^t_{k+2}}_{\mathbf{H}_t} ]}{Z^{1,k+1,(-1,1),\vec{x},\vec{y},+\infty,\mathfrak{H}^t_{k+2}}_{\mathbf{H}_t}  } \bigg]\\
\leq & \mathbb{E}\left[\mathbbm{1}_{\mathsf{G}}\cdot \varepsilon  \right]\leq \varepsilon.
\end{align*}
In conclusion, we obtain
\begin{align*}
\mathbb{P}(\mathsf{D}_1\cup \mathsf{D}_2)\leq  &\mathbb{P}((\mathsf{D}_1\cup \mathsf{D}_2)\cap \mathsf{G})+\mathbb{P}( \mathsf{G}^{\textup{c}})\\
\leq &\mathbb{P}( \mathsf{D}_1\cap \mathsf{G} )+\mathbb{P}( \mathsf{D}_2\cap\mathsf{D}_3 \cap \mathsf{G} )+\mathbb{P}( \mathsf{D}^{\textup{c}}_3\cap \mathsf{G} )+\mathbb{P}( \mathsf{G}^{\textup{c}})\\
\leq &4\varepsilon.
\end{align*}
The proof is finished.
\end{proof}

From the stationarity of $\mathfrak{H}^t(u)+2^{-1}u^2$, we immediate obtain the following corollary.
\begin{corollary}\label{cor:ordered}
Fix  $k\in\mathbb{N}$ and $\varepsilon\in (0,1]$. Let $\rho =\rho (\varepsilon,k) $ and $t_0=t_0(\varepsilon ,k)\geq 1$ be the constants in Proposition \ref{pro:ordered}. Then for all $x_0\in\mathbb{R}$, we have
\begin{align*}
\mathbb{P}\bigg( \inf_{u\in [x_0,x_0+2]}(\mathfrak{H}^t_k(u)-\mathfrak{H}^t_{k+1}(u))< \rho \bigg)< \varepsilon.
\end{align*}
\end{corollary}

We are ready to demonstrate that Gibbs property survives under weak convergence of line ensembles in Proposition \ref{pro:BGP}. We adopt a coupling argument as used in \cite[Proof of Proposition 5.2]{CH16}.  We also take care of the issue that the interaction Hamiltonian varies in $t$ in our case while it stays the same as in \cite{CH16}. Theorem \ref{thm:main2}(ii) now immediately follows from Proposition \ref{pro:ordered} and Proposition \ref{pro:BGP}.

\begin{proposition}\label{pro:BGP}
Any subsequential limit of $\mathfrak{H}^t$ satisfies the $\mathbf{H}_{\textup{ord}}$-Brownian Gibbs property. 
\end{proposition}

\begin{proof}
Let $ t_N\to \infty$ be an increasing sequence such that $\mathfrak{H}^{t_N}$ converges weakly to $\mathfrak{H}^{\infty}$ as $\mathbb{N}\times\mathbb{R}$-indexed line ensembles, see Definition \ref{def:weakconvergence}. Recall the the interaction Hamiltonian is $\mathbf{H}_t(x)=e^{t^{1/3}x}$. For simplicity, we denote $\mathfrak{H}^N=\mathfrak{H}^{t_N}$ and $\mathbf{H}_N=\mathbf{H}_{t_N}$. 

Fix an index $i\in\mathbb{N}$ and two numbers $a<b$. We will show that the law of $\mathfrak{H}^\infty$ is unchanged when one resamples the trajectory of $\mathfrak{H}^\infty_i$ between $a$ and $b$ according a Brownian bridge which avoids $\mathfrak{H}^\infty_{i-1}$ and $\mathfrak{H}^\infty_{i+1}$. The argument can be easily generalized to take care of multiple curves resampling so we choose to illustrate the argument with the single curve resampling, see Figure \ref{figure:BGP} for an illustration. Note that the Brownian Gibbs property is equivalent to this resampling invariance, hence finishing the proof.

\begin{figure}
\begin{tikzcd}  
\mathfrak{H}^{N,\textup{re}} \arrow[rr,equal ,"(d)","(1)"'] \arrow[dd,"\textup{a.s.}","(2)'"']& & \mathfrak{H}^{N } \arrow[dd,"\textup{a.s.}","(2)"']  \\
 & &\\
\mathfrak{H}^{\infty,\textup{re}} \arrow[rr,equal ,"(d)","(1)'"'] & & \mathfrak{H}^{\infty }
\end{tikzcd}
\caption{(1) is equivalent to the Brownian Gibbs property when resampling a single curve. The {\bf goal} is to prove (1)', which implies the Brownian Gibbs property for the subsequential limit line ensemble.   (1)' follows from the convergence in (2) and (2)'. (2) follows from the Skorohod representation theorem and (2)' is proved in Lemma \ref{lem:ell} and Lemma \ref{lem:ellconvergence}.}\label{figure:BGP}
\end{figure}

Note that $  C(\mathbb{N} \times \mathbb{R},\mathbb{R}) $ (with the topology in Definition \ref{def:line-ensemble}) is separable due to the Stone–Weierstrass theorem. Hence the Skorohod representation theorem \cite[Theorem 6.7]{Bil} applies. There exists a probability space $(\Omega,\mathcal{B},\mathbb{P})$ on which $\mathfrak{H}^N$ for $N\in\mathbb{N}\cup\{\infty\}$ are defined and almost surely $\mathfrak{H}^N(\omega)$ converges to $\mathfrak{H}^\infty(\omega)$ in the topology of $C(\mathbb{N} \times \mathbb{R},\mathbb{R})$.\\[0.1cm]
\indent Let $\{B^{\ell}\}_{\ell\geq 1}$ be a sequence of independent Brownian bridges defined on $[a,b]$ with $B^{\ell}(a)=B^{\ell}(b)=0$. Let $\{U_\ell\}_{\ell\geq 1}$ be a sequence of of independent uniform distributions on $[0,1]$. We augment the probability space $(\Omega,\mathcal{B},\mathbb{P})$ to accommodate all of $\{B^{\ell}\}_{\ell\geq 1}$ and $\{U_\ell\}_{\ell\geq 1}$ in an independent manner.\\[0.1cm]
\indent {\bf Step one}, we define the $\ell$-th candidate of the resampling trajectory. For $u\in [a,b]$, define
\begin{align*}
\mathfrak{H}^{N,\ell}_i(u)\coloneqq \mathfrak{H}^{N}_i(a)+B^\ell(u)+\frac{u-a}{b-a}\cdot (\mathfrak{H}^{N}_i(b)-\mathfrak{H}^{N}_i(a)),
\end{align*}
and $\mathfrak{H}^{N,\ell}_i(u)\coloneqq \mathfrak{H}^{N}_i(u)$ for $u\in (-\infty,a)\cap (b,\infty)$. Similarly,
\begin{align*}
\mathfrak{H}^{\infty,\ell}_i(u)\coloneqq \mathfrak{H}^{\infty}_i(a)+B^\ell(u)+\frac{u-a}{b-a}\cdot (\mathfrak{H}^{\infty}_i(b)-\mathfrak{H}^{\infty}_i(a)),
\end{align*}
and $\mathfrak{H}^{\infty,\ell}_i(u)\coloneqq \mathfrak{H}^{\infty}_i(u)$ for $u\in (-\infty,a)\cap (b,\infty)$.\\[0.1cm]
\indent{\bf Step two}, we check whether
\begin{align*}
U_{\ell}\leq W(N,\ell):=W^{i,i,(a,b),\mathfrak{H}^N_{i}(a),\mathfrak{H}^N_{i}(b),\mathfrak{H}^N_{i-1},\mathfrak{H}^N_{i+1}}_{\mathbf{H}_{N}}(\mathfrak{H}^{N,\ell}_i), 
\end{align*}
and \textbf{accept} the candidate resampling $\mathfrak{H}^{N,\ell}_i$ if this event occurs. Similarly, we \textbf{accept} the candidate resampling $\mathfrak{H}^{\infty,\ell}_i$ if it does not intersect  $\mathfrak{H}^{\infty}_{i-1}$ or $\mathfrak{H}^{\infty}_{i+1}$ in $[a,b]$. For $N\in \mathbb{N}\cup\{\infty\}$, define $\ell(N)$ to be the minimal value of $\ell$ of which we accept $\mathfrak{H}^{N,\ell}_i$. Write $\mathfrak{H}^{N,\textup{re}}$ for the line ensemble with the $i$-th curve replaced by $ \mathfrak{H}^{N,\ell(N)}_i$. The line ensemble $\mathfrak{H}^{N,\textup{re}} $ satisfies the $\mathbf{H}_N$-Brownian Gibbs property on $\{i\}\times [a,b]$. More precisely, we have the following.
\begin{claim}\label{clm:resample}
Let $\mathcal{F}_{ {ext}}^N(\{i\}\times (a,b))$ be the sigma-field generated by $\mathfrak{H}^N$ restricted on $\mathbb{N}\times \mathbb{R}\setminus \{i\}\times (a,b)$. Then for any Borel function $F:C(\{i\}\times [a,b],\R)\to \mathbb{R}$, we have
\begin{align*}
\mathbb{E}\left[F(\mathfrak{H}^{N,\textup{re}}_i|_{[a,b]}) \,|\,\mathcal{F}_{ {ext}}^N(\{i\}\times (a,b))\right]=\mathbb{E}^{i,i,(a,b),\mathfrak{H}^N_{i}(a),\mathfrak{H}^N_{i}(b),\mathfrak{H}^N_{i-1},\mathfrak{H}^N_{i+1}}_{\mathbf{H}_{N}}[F(\mathcal{L})].
\end{align*}
\end{claim}
\begin{proof}Write
\begin{align*}
\mathbb{E}\left[F(\mathfrak{H}^{N,\textup{re}}_i|_{[a,b]}) \,|\,\mathcal{F}_{ {ext}}^N(\{i\}\times (a,b))\right]=\sum_{\ell=1}^\infty  \mathbb{E}\left[F(\mathfrak{H}^{N,\ell }_i|_{[a,b]})\cdot\mathbbm{1}\{\ell=\ell(N)\} \,|\,\mathcal{F}_{ {ext}}^N(\{i\}\times (a,b))\right].
\end{align*}
Because $ B^j $ and $U_j$ are independent and $$\{\ell=\ell(N)\}=\bigcap_{j=1}^{\ell-1}\{ W(N,j)<U_j \}\cap \{W(N,\ell)\geq U_\ell \}, $$
we have
\begin{align*}
&\mathbb{E}\left[F(\mathfrak{H}^{N,\ell }_i|_{[a,b]})\cdot\mathbbm{1}\{\ell=\ell(N)\} \,|\,\mathcal{F}_{ {ext}}^N(\{i\}\times (a,b))\right]\\
=
&\prod_{j=1}^{\ell-1}\mathbb{E}\left[ \mathbbm{1}\{W(N,j)<U_j\} \,|\,\mathcal{F}_{ {ext}}^N(\{i\}\times (a,b))\right] \mathbb{E}\left[F(\mathfrak{H}^{N,\ell }_i|_{[a,b]})\cdot\mathbbm{1}\{W(N,\ell)\geq U_\ell \} \,|\,\mathcal{F}_{ {ext}}^N(\{i\}\times (a,b))\right].
\end{align*}
From the definition of  $B^j$, $U_j$, $\mathfrak{H}^{N,\ell}_i$ and $W(N,j)$, we have
\begin{align*}
\mathbb{E}\left[ \mathbbm{1}\{W(N,j)<U_j\} \,|\,\mathcal{F}_{ {ext}}^N(\{i\}\times (a,b))\right]=&\mathbb{E}\left[  1- W(N,j)  \,|\,\mathcal{F}_{ {ext}}^N(\{i\}\times (a,b))\right]\\
=&1-Z^{i,i,(a,b),\mathfrak{H}^N_{i}(a),\mathfrak{H}^N_{i}(b),\mathfrak{H}^N_{i-1},\mathfrak{H}^N_{i+1}}_{\mathbf{H}_{N}}.
\end{align*}
Similarly,
\begin{align*}
&\mathbb{E}\left[F(\mathfrak{H}^{N,\ell }_i|_{[a,b]})\cdot\mathbbm{1}\{W(N,\ell)\geq U_\ell \} \,|\,\mathcal{F}_{ {ext}}^N(\{i\}\times (a,b))\right]=\mathbb{E}\left[F(\mathfrak{H}^{N,\ell }_i|_{[a,b]})\cdot  W(N,\ell)  \,|\,\mathcal{F}_{ {ext}}^N(\{i\}\times (a,b))\right]\\
=&\mathbb{E}^{i,i,(a,b),\mathfrak{H}^N_{i}(a),\mathfrak{H}^N_{i}(b) }_{\free}\left[F(\mathcal{L})\cdot W^{i,i,(a,b),\mathfrak{H}^N_{i}(a),\mathfrak{H}^N_{i}(b),\mathfrak{H}^N_{i-1},\mathfrak{H}^N_{i+1}}_{\mathbf{H}_{N}}(\mathcal{L})\right].
\end{align*}
As a result,
\begin{align*}
&\mathbb{E}\left[F(\mathfrak{H}^{N,\textup{re}}_i|_{[a,b]}) \,|\,\mathcal{F}_{ {ext}}^N(\{i\}\times (a,b))\right]\\
=&\mathbb{E}^{i,i,(a,b),\mathfrak{H}^N_{i}(a),\mathfrak{H}^N_{i}(b) }_{\free}\left[F(\mathcal{L})\cdot W^{i,i,(a,b),\mathfrak{H}^N_{i}(a),\mathfrak{H}^N_{i}(b),\mathfrak{H}^N_{i-1},\mathfrak{H}^N_{i+1}}_{\mathbf{H}_{N}}(\mathcal{L})\right]\\
&\times \sum_{\ell=1}^{\infty} \left(1-Z^{i,i,(a,b),\mathfrak{H}^N_{i}(a),\mathfrak{H}^N_{i}(b),\mathfrak{H}^N_{i-1},\mathfrak{H}^N_{i+1}}_{\mathbf{H}_{N}}\right)^{\ell-1}\\
=&\frac{\mathbb{E}^{i,i,(a,b),\mathfrak{H}^N_{i}(a),\mathfrak{H}^N_{i}(b) }_{\free}\left[F(\mathcal{L})\cdot W^{i,i,(a,b),\mathfrak{H}^N_{i}(a),\mathfrak{H}^N_{i}(b),\mathfrak{H}^N_{i-1},\mathfrak{H}^N_{i+1}}_{\mathbf{H}_{N}}(\mathcal{L})\right] }{ Z^{i,i,(a,b),\mathfrak{H}^N_{i}(a),\mathfrak{H}^N_{i}(b),\mathfrak{H}^N_{i-1},\mathfrak{H}^N_{i+1}}_{\mathbf{H}_{N}}}\\
=&\mathbb{E}^{i,i,(a,b),\mathfrak{H}^N_{i}(a),\mathfrak{H}^N_{i}(b),\mathfrak{H}^N_{i-1},\mathfrak{H}^N_{i+1} }_{\mathbf{H}_N}[F(\mathcal{L})].
\end{align*}
\end{proof}
By Claim \ref{clm:resample}, the $\mathbf{H}_t$-Brownian Gibbs property implies that for $N\in\mathbb{N}$,
\begin{align}\label{equ:resample-invariance}
\mathfrak{H}^{N,\textup{re}}\overset{(d)}{=\joinrel=}\mathfrak{H}^{N}.
\end{align}
Furthermore, \eqref{equ:resample-invariance} holds for all $i\in \mathbb{N}$ and $a<b$ implies that $\mathfrak{H}^N$ satisfies $\mathbf{H}_N$-Brownian Gibbs property provided only one line is resampled. Our {\bf goal} is to show that \eqref{equ:resample-invariance} holds for $N=\infty$. As a result, $\mathfrak{H}^\infty$ satisfies the Brownian Gibbs property when resampling a single curve.

It suffices to show that almost surely $\ell(\infty)$ is finite  and that almost surely $\ell(N)$ converges to $\ell(\infty)$. Suppose these two hold, then $\mathfrak{H}^{N,\textup{re}}$ converges to $\mathfrak{H}^{\infty,\textup{re}}$ in $   C ( \mathbb{N} \times \mathbb{R},\mathbb{R})  $ almost surely. Hence we have $\mathfrak{H}^{N,\textup{re}}$ converges weakly to $\mathfrak{H}^{\infty,\textup{re}}$ as $\mathbb{N}\times\mathbb{R}$-indexed line ensembles. See Definition \ref{def:weakconvergence} As a consequence, $\mathfrak{H}^{\infty,\textup{re}}$ has the same distribution as $\mathfrak{H}^{\infty}$.
\end{proof}

\begin{lemma}\label{lem:ell}
Almost surely $\ell(\infty)$ is finite.
\end{lemma} 
\begin{proof}
Let $\mathcal{F}_{ {ext}}^{\infty}(\{i\}\times (a,b))$ be the sigma-field generated by $\mathfrak{H}^\infty$ restricted on $\mathbb{N}\times \mathbb{R}\setminus \{i\}\times (a,b)$. Then
\begin{align*}
\mathbb{E}[\mathbbm{1}\{\ell(\infty)\geq \ell\}\,|\, \mathcal{F}_{ {ext}}^{\infty}(\{i\}\times (a,b))]=\left(1-Z^{i,i,(a,b),\mathfrak{H}^{\infty}_{i}(a),\mathfrak{H}^{\infty}_{i}(b),\mathfrak{H}^{\infty}_{i-1},\mathfrak{H}^{\infty}_{i+1}}\right)^{\ell-1}.
\end{align*}
Here
\begin{align*}
Z^{i,i,(a,b),\mathfrak{H}^{\infty}_{i}(a),\mathfrak{H}^{\infty}_{i}(b),\mathfrak{H}^{\infty}_{i-1},\mathfrak{H}^{\infty}_{i+1}}=\mathbb{P}_{\free}^{i,i,(a,b),\mathfrak{H}^{\infty}_{i}(a),\mathfrak{H}^{\infty}_{i}(b)}\left( \mathfrak{H}^{\infty}_{i-1}(u)<\mathcal{L}(u)< \mathfrak{H}^{\infty}_{i+1}(u)\ \textup{in}\ [a,b] \right).
\end{align*}
Define the event 
\begin{align*}
\mathsf{S}\coloneqq \left\{ \mathfrak{H}^{\infty}_{i-1}(u)<\mathfrak{H}^{\infty}_{i+1}(u)\ \textup{in}\  [a,b]  \right\}\cap \{ \mathfrak{H}^{\infty}_{i-1}(u)<\mathfrak{H}^{\infty}_{i} (u)< \mathfrak{H}^{\infty}_{i+1}(u)\ \textup{for}\ u\in \{a,b\}   \}.
\end{align*} 
As $\mathsf{S}$ occurs, $Z^{i,i,(a,b),\mathfrak{H}^{\infty}_{i}(a),\mathfrak{H}^{\infty}_{i}(b),\mathfrak{H}^{\infty}_{i-1},\mathfrak{H}^{\infty}_{i+1}}>0$. From Corollary \ref{cor:ordered}, we have $\mathbb{P}(\mathsf{S})=1$. As a result,
\begin{align*}
\mathbb{P}(\{\ell(\infty)=\infty\})=\mathbb{P}(\{\ell(\infty)=\infty\}\cap\mathsf{S})=\mathbb{E}\left[\mathbbm{1}_{\mathsf{S}}\cdot \mathbb{E}[\mathbbm{1}\{\ell(\infty)=\infty \}\,|\, \mathcal{F}_{ {ext}}^{\infty}(\{i\}\times (a,b))]  \right]=0.
\end{align*}
\end{proof}

\begin{lemma}\label{lem:ellconvergence}
Almost surely $\ell(N)$ converges to $\ell(\infty)$.
\end{lemma}

\begin{proof}
Let $\mathsf{E}$ be the event such that the following five conditions hold \\[-0.3cm]
\begin{enumerate}
\item $\ell(\infty)<\infty$\\[-0.3cm]
\item $U_{\ell}\in (0,1)$ for all $\ell\in\mathbb{N}$\\[-0.3cm]
\item $\displaystyle\sup_{u\in [a,b]} \mathfrak{H}^{\infty,\ell}_{i}(u)-\mathfrak{H}^{\infty}_{i-1}(u) \neq0$ for all $\ell\in\mathbb{N}$\\[0.05cm]
\item $\displaystyle\sup_{u\in [a,b]} \mathfrak{H}^{\infty}_{i+1}(u)-\mathfrak{H}^{\infty,\ell}_{i}(u)\neq0$ for all $\ell\in\mathbb{N}$\\[0.05cm]
\item $\mathfrak{H}^{N}$ converges to $\mathfrak{H}^{\infty}$ in $C(\mathbb{N}\times\mathbb{R},\mathbb{R})$.\\[-0.3cm]
\end{enumerate}

It follows from Lemma \ref{lem:ell} and Corollary \ref{cor:ordered} that $\mathbb{P}(\mathsf{E})=1$. We will show that when $\mathsf{E}$ occurs, $\ell(N) \to \ell(\infty)$. 
From now on we fix a realization $\omega\in \mathsf{E}$ and the constants below may depend on $\omega$. By the definition of $\ell(\infty)$ and $\mathsf{E}$,  there exists a $\delta>0$ such that
\begin{align*}
\displaystyle\sup_{u\in [a,b]} \mathfrak{H}^{\infty,\ell(\infty)}_{i}(u)-\mathfrak{H}^{\infty}_{i-1}(u)\leq -2\delta<0,\\
\displaystyle\sup_{u\in [a,b]} \mathfrak{H}^{\infty}_{i+1}(u)-\mathfrak{H}^{\infty,\ell(\infty)}_{i}(u)\leq -2\delta < 0.  
\end{align*}
Moreover, since $\mathfrak{H}^N$ converges uniformly to $\mathfrak{H}^{\infty}$ on $[1,i+1]_{\mathbb{Z}}\times [a,b]$, we have for $N$ large enough,
\begin{align*}
\displaystyle\sup_{u\in [a,b]} \mathfrak{H}^{N,\ell(\infty)}_{i}(u)-\mathfrak{H}^{N}_{i-1}(u)\leq -\delta<0,\\
\displaystyle\sup_{u\in [a,b]} \mathfrak{H}^{N}_{i+1}(u)-\mathfrak{H}^{N,\ell(\infty)}_{i}(u)\leq -\delta< 0.  
\end{align*}
As a consequence, 
\begin{align*}
W(N,\ell(\infty))\geq \exp(-2(b-a)e^{-t_N^{1/3}\delta})  	
\end{align*}
which converges to $1$. Because $U_{\ell(\infty)}<1$, we have $W(N,\ell(\infty))>U_{\ell(\infty)}$ for $N$ large enough. Hence
\begin{align*}
\limsup_{N\to\infty}\ell(N)\leq \ell(\infty). 
\end{align*}

On the other hand, for all $1\leq\ell<\ell(\infty)$, we have either 
\begin{align*}
\displaystyle\sup_{u\in [a,b]} \mathfrak{H}^{\infty,\ell }_{i}(u)-\mathfrak{H}^{\infty}_{i-1}(u)>0
\end{align*}
or
\begin{align*}
\displaystyle\sup_{u\in [a,b]} \mathfrak{H}^{\infty}_{i+1}(u)-\mathfrak{H}^{\infty,\ell }_{i}(u)> 0.  
\end{align*}
We assume that $\displaystyle\sup_{u\in [a,b]} \mathfrak{H}^{\infty,\ell }_{i}(u)-\mathfrak{H}^{\infty}_{i-1}(u)>0$ occurs and denote $4\delta=\displaystyle\sup_{u\in [a,b]} \mathfrak{H}^{\infty,\ell }_{i}(u)-\mathfrak{H}^{\infty}_{i-1}(u)$. By the continuity of $\mathfrak{H}^{\infty,\ell }_{i}(u)-\mathfrak{H}^{\infty}_{i-1}(u)$, there exists an interval $I\subset [a,b]$ such that $   \displaystyle\inf_{u\in I} \mathfrak{H}^{\infty,\ell }_{i}(u)-\mathfrak{H}^{\infty}_{i-1}(u)\geq 2\delta$. Because $\mathfrak{H}^N$ converges uniformly to $\mathfrak{H}^{\infty}$ on $[1,i+1]_{\mathbb{Z}}\times [a,b]$, we have for $N$ large enough,
\begin{align*}
\inf_{u\in I} \mathfrak{H}^{\infty,\ell }_{i}(u)-\mathfrak{H}^{\infty}_{i-1}(u)\geq \delta.
\end{align*}
As a consequence,
\begin{align*}
W(N,\ell)\leq \exp(-|I|e^{t^{1/3}_N \delta})
\end{align*}
which converges to $0$. Because $U_{\ell}>0$, we obtain that for $N$ large enough,
\begin{align*}
W(N,\ell)<U_{\ell}. 
\end{align*}
Since the above argument holds for all $1\leq \ell< \ell(\infty)$, we deduce
\begin{align*}
\liminf_{N\to\infty}\ell(N)\geq \ell(\infty).
\end{align*}
Hence $\ell(N)$ converges to $\ell(\infty)$  and the proof is finished.
\end{proof} 

\section{Proof of Proposition~\ref{pro:Z_lowerbound_k}}\label{sec:Z}
In this section we seek to prove a uniform and quantitative estimate, Proposition~\ref{pro:Z_lowerbound_k}, on the normalizing constant when resampling multiple curves. Proposition~\ref{pro:Z_lowerbound_k} is the main technical result of this paper and will be further exploited in \cite{Wu21} to study the Brownian regularity for the KPZ line ensemble.  

\begin{figure}
\includegraphics[width =12cm]{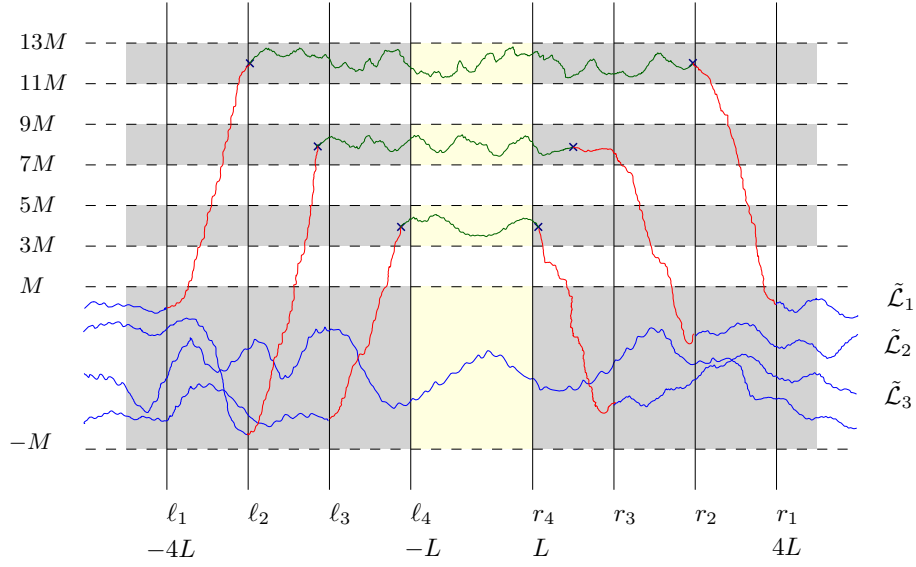}
\caption{Illustration for bounding $Z^{1,3,(-L,L),\vec{x},\vec{y},+\infty,f_{4}}_{\mathbf{H}_t}$ using two-step resampling. The blue trajectories serve as $M$-{\bf Good} boundary data for resampling $\mathcal{L}_1^t, \mathcal{L}_2^t, \mathcal{L}_3^t$. We hope they stay in the preferable yellow region over $[-L,L]$. In that case, curves $\mathcal{L}_1^t, \mathcal{L}_2^t, \mathcal{L}_3^t.$ are separated away from each other by a distance of $2M$. We will do a two-step Gibbsian resampling to show that such configurations happen with a positive (uniform in $t$) probability. The red trajectories are sampled in step one to have them rise up and green trajectories are sampled in step two to have them stay in the yellow region. See details of these two steps respectively in Lemma \ref{lem:lowerbound-k} and Lemma \ref{lem:spacing-k}.}\label{figure:separation}
\end{figure}

The main difficulty in proving Proposition~\ref{pro:Z_lowerbound_k} dates back to estimating the Boltzmann weight when curves go out of order, in particular, when there are multiple curves interacting with their neighbors simultaneously.  We introduce a novel inductive two-step resampling procedure, which first raises the curves one by one starting from the lowest index curve and then lowers them in a reversed order. In doing so, we force curves to stay in the preferable region, see an illustration in Figure \ref{figure:separation}. More precisely, we establish an estimate on the probability of curves being well-separated in Proposition \ref{pro:spacing-k}, which enables us to estimate the Boltzmann weight and furthermore provide a positive lower bound of the normalizing constant in Corollary \ref{cor:Z_lowerbound-k}.

We begin with setting some parameters. Fix $k\in\mathbb{N}$, $L\geq 1$, $t\geq 1$ and take 
\begin{equation}\label{Mbound}
M\geq  L^{1/2}.
\end{equation} 
For $1\leq j\leq k+1$, let
\begin{align*}
\ell_j\coloneqq -(k+2-j)L,\  r_j\coloneqq (k+2-j)L.
\end{align*}
Consider continuous functions $f=(f_1,f_2,\dots, f_k)$ with $f_j\in C([\ell_1,\ell_j]\cup [r_j,r_1],\mathbb{R})$ for $1\leq j\leq k$ and $f_{k+1}\in C\left([\ell_1 ,r_1],\mathbb{R} \right)$. We view $f$ as a boundary condition and call it $M$-{\bf Good} provided the absolute value of each component is bounded by $M$.  

We will consider line ensembles defined on $ {\cup_{j=1}^k\ \{j\}\times [\ell_j,r_j]}$. For any $f$ as above, we let $\mathfrak{X}_f$ be the collection of functions $\mathcal{L}=(\mathcal{L}_1,\mathcal{L}_2,\dots,\mathcal{L}_{k})$ with $\mathcal{L}_j\in C([\ell_j,r_j],\mathbb{R})$, $\mathcal{L}_j(\ell_j)=f(\ell_j)$ and $\mathcal{L}_j(r_j)=f(r_j)$. Given $\mathcal{L}\in\mathfrak{X}_f$, we may extend the domain of $\mathcal{L}$ using $f$. 
\begin{align}\label{equ:extension}
\mathcal{L}_{j,f}(u)\coloneqq \left\{\begin{array}{cc}
\mathcal{L}_j(u) & u\in [\ell_j,r_j],\\
f_j(u) & u\in [\ell_1,\ell_j)\cup (r_j,r_1].
\end{array} \right.
\end{align}
Note that $\mathcal{L}_f\in C([1,k]_{\mathbb{Z}}\times [\ell_1,r_1],\mathbb{R})$. \\[0.1cm] 
\indent For $j\in [1,k]_{\Z} $, let $B_j$ be independent Brownian bridges defined on $[\ell_j,r_j]$ with $B_j(\ell_j  )=f_j(\ell_j )$ and $B_j(r_j  )=f_j(r_j )$. The law of $B_j$ is given by $\mathbb{P}_{\free}^{j,j,(\ell_j,r_j),f_j(\ell_j),f_j(r_j)}$ and we denote by $\mathbb{P}_{\free}$ the joint law of $(B_1,B_2\dots, B_{k})$ and by $\mathbb{E}_{\free}$ the expectation of $\bP_{\free}$. We may view $\mathbb{P}_{\free}$ as a probability measure defined on $\mathfrak{X}_f$.

Next, we define two probability measures, $\bP$ and $\tilde{\bP}$, on  $\mathfrak{X}_f$ through specifying their Radon-Nikodym derivatives with respect to $\mathbb{P}_{\free}$. 
\begin{align}
&\frac{\textup{d}\mathbb{P} }{\textup{d}\mathbb{P}_{\free}}(\cL):=\frac{1}{Z}W^{1,k,(\ell_1,r_1),\vec{x},\vec{y},+\infty,f_{k+1}}_{\mathbf{H}_t}(\cL_f) \label{resample-k},\\
&\frac{\textup{d}\tilde{\mathbb{P}}}{\textup{d}\mathbb{P}_{\free}}(\cL ):=\frac{1}{\tilde{Z} }W^{1,k,(\ell_1,r_1),\vec{x},\vec{y},+\infty,f_{k+1}}_{\mathbf{H}_t,(-L,L)}(\cL_f),\label{resample-kt}
\end{align}
where
\begin{align*}
&Z\coloneqq \mathbb{E}_{\free}\left[W^{1,k,(\ell_1,r_1),\vec{x},\vec{y},+\infty,f_{k+1}}_{\mathbf{H}_t}(\cL_f)\right],\\
&\tilde{Z}\coloneqq \mathbb{E}_{\free}\left[ W^{1,k,(\ell_1,r_1),\vec{x},\vec{y},+\infty,f_{k+1}}_{\mathbf{H}_t,(-L,L )}(\cL_f) \right],
\end{align*}
and $\vec{x}=(f_j(\ell_1))_{j=1}^k$, $\vec{y}=(f_j(r_1))_{j=1}^k$. See Definition~\ref{def:H_Brownian} for the definition of the Boltzmann weights. We denote by $\mathbb{E}$ and $\tilde{\mathbb{E}}$ the expectation of $\bP$ and $\tilde{\bP}$ respectively.

Consider the event that curves are well separated by order $M$ at the endpoints of interval $[-L,L]$.
\begin{align*} 
\mathsf{E}:= \left\{   {\mathcal{L}}_j(\pm L) \in \big[  (4k-4j+3) M,(4k-4j+5)M \big]\ \textup{for all}\ j\in [1,k]_{\Z} \right\}.
\end{align*}

Proposition \ref{pro:spacing-k} below provides a lower bound of $\tilde{\mathbb{P}}(\mathsf{E})$.

\begin{proposition}\label{pro:spacing-k}
Fix $k\in\mathbb{N}$, $L\geq 1$, $t\geq 1$ and $M\geq L^{1/2}$. Let $f$ be a $M$-{\bf Good} boundary condition and $\tilde{\mathbb{P}}$ be the probability measure on $\mathfrak{X}_f$ defined in \eqref{resample-kt}. Then there exists a constant $D_3=D_3(k)$ depending only on $k$ such that  
\begin{align*}
\tilde{\mathbb{P}}(\mathsf{E})\geq D_3^{-1}  e^{-D_3 ( L^{-1}M^2+ L)}.
\end{align*}
\end{proposition}
We postpone the proof of Proposition~\ref{pro:spacing-k} to the end of this section. Next, we prove a lower bound for the  the normalizing constant under the law of  $\tilde{\mathbb{P}}$ (Corollary \ref{cor:Z_lowerbound-k}) and under the law of $\mathbb{P}$ (Proposition \ref{pro:Z_small-k}) respectively.
\begin{corollary}\label{cor:Z_lowerbound-k}
Fix $k\in\mathbb{N}$, $L\geq 1$, $t\geq 1$ and $M\geq L^{1/2}$. Let $f$ be a $M$-{\bf Good} boundary condition and $\tilde{\mathbb{P}}$ be the probability measure on $\mathfrak{X}_f$ defined in \eqref{resample-kt}. Then there exists a constant $D_4=D_4(k)$ depending only on $k$ such that 
\begin{align*}
\tilde{\mathbb{P}}\left(Z^{1,k,(-L,L),\vec{x} ,\vec{y} ,+\infty,f_{k+1}}_{\mathbf{H}_t}\geq  D_4^{-1} e^{-2kL}  \right)\geq D_4^{-1} e^{-D_4 (L^{-1}M^2 + L)}.
\end{align*} 
Here $\vec{x}=( {\mathcal{L}}_j(-L))_{j=1}^k$, $\vec{y}=( {\mathcal{L}}_j(L))_{j=1}^k$. 
\end{corollary}

\begin{proof}
Note that $$2^{-1}L^{-1}\left((u+L) {\mathcal{L}}_j(L)+(L-u) {\mathcal{L}}_j(-L) \right)$$
is the interpolation function connecting $\left(-L,  {\mathcal{L}}_j(-L)\right)$ and $\left(L,  {\mathcal{L}}_j(L)\right)$. Consider the following event $\mathsf{Osc}$ where every layer $ {\mathcal{L}}_j, j=1,2,\cdots,k$ does not deviate from the linear interpolation by $M$,
\begin{align*}
\mathsf{Osc}\coloneqq \left\{ \sup_{1\leq j\leq k}\sup_{u\in [-L,L]}\left| \mathcal{L}_j(u)-2^{-1}L^{-1}\left((u+L) {\mathcal{L}}_j(L)+(L-u) {\mathcal{L}}_j(-L) \right) \right|\leq M  \right\}.
\end{align*}

Recall that event $\mathsf{E}$ says that the boundary values, i.e. $\vec{x}$ and $\vec{y}$ are separated by at least $2M$. Suppose $\mathsf{E}$ and $\mathsf{Osc}$ both occur, then $( {\mathcal{L} }_1, {\mathcal{L} }_2,\dots , {\mathcal{L} }_k,f_{k+1})$ remains ordered on $[-L,L]$. In particular, the Boltzmann weight is bounded below,
\begin{equation}\label{equ:WWW}
W_{\mathbf{H}_t}^{1,k,(-L,L),\vec{x} ,\vec{y} ,+\infty,f_{k+1}}( {\mathcal{L} } )\cdot\mathbbm{1}_{\mathsf{E}}\cdot\mathbbm{1}_{\mathsf{Osc}} \geq  e^{-2kL}\cdot\mathbbm{1}_{\mathsf{E}}\cdot\mathbbm{1}_{\mathsf{Osc}}. 
\end{equation}
Because $M\geq L^{1/2}$, there exists $D=D(k)$ such that
\begin{align}\label{equ:PPP}
\mathbb{P}_{\free}^{1,k,(-L,L),\vec{x} ,\vec{y}}(\mathsf{Osc})=\prod_{j=1}^k\mathbb{P}_{\free}^{j,j,(-L,L),0 ,0}\left(\sup_{u\in [-L,L]}|\cL(u)|\leq M  \right)\geq D^{-1}.
\end{align}

Recall that in \eqref{def:normalcont_Brownian},
\begin{align*}
Z_{\mathbf{H}_t}^{1,k,(-L,L),\vec{x},\vec{y},+\infty,f_{k+1}}:=\mathbb{E}_{\free}^{1,k,(-L,L),\vec{x} ,\vec{y}}\left[W_{\mathbf{H}_t}^{1,k,(-L,L),\vec{x} ,\vec{y} ,+\infty,f_{k+1}}  \right].
\end{align*}
Combining \eqref{def:normalcont_Brownian}, \eqref{equ:WWW} and \eqref{equ:PPP}, we get
\begin{align*}
Z_{\mathbf{H}_t}^{1,k,(-L,L),\vec{x},\vec{y},+\infty,f_{k+1}}\cdot\mathbbm{1}_{\mathsf{E}}\geq D^{-1}e^{-2kL}\cdot\mathbbm{1}_{\mathsf{E}}. 
\end{align*}
This implies that $\mathsf{E}\subset \left\{ Z^{1,k,(-L,L),\vec{x} ,\vec{y},+\infty,f_{k+1}}_{\mathbf{H}_t}\geq  D^{-1}e^{-2kL} \right\}$ and the desired result then follows from Proposition \ref{pro:spacing-k}.
\end{proof}

\begin{proposition}\label{pro:Z_small-k}
Fix $k\in\mathbb{N}$, $L\geq 1$, $t\geq 1$ and $M\geq L^{1/2}$. Let $f$ be a $M$-{\bf Good} boundary condition and $\mathbb{P}$ be the probability measure on $\mathfrak{X}_f$ defined in \eqref{resample-k}. Then there exists a constant $D_5=D_5(k)$ depending only on $k$ such that for all $\varepsilon\in (0,1]$, we have
$$\mathbb{P}\left( Z^{1,k,(-L,L), \vec{x}, \vec{y} ,+\infty,f_{k+1}}_{\mathbf{H}_t}\leq \varepsilon D_5 ^{-1}  e^{-D_5   (L^{-1} M^2+ L)}  \right)\leq \varepsilon.$$
Here $\vec{x}=( {\mathcal{L}}_j(-L))_{j=1}^k$, $\vec{y}=( {\mathcal{L}}_j(L))_{j=1}^k$. 
\end{proposition}
\begin{proof}
Let $\tilde{\mathbb{P}}$ be the probability measure on $\mathfrak{X}_f$ defined in \eqref{resample-kt}. From \eqref{resample-k} and \eqref{resample-kt}, we have
\begin{align}\label{equ:PoverP}
\frac{\textup{d}\mathbb{P}}{\textup{d}\tilde{\mathbb{P}}}(\cL)\propto W^{1,k,(-L,L),\vec{x},\vec{y},+\infty,f_{k+1}}_{\mathbf{H}_t}(\cL|_{[1,k]_{\mathbb{Z}}\times [-L,L]}),
\end{align}
where $\vec{x}=(\cL_j(-L))_{j=1}^k$ and $\vec{y}=(\cL_j(L))_{j=1}^k$. 

let $\mathfrak{X}'$ be the collection of functions $\mathcal{L}'=(\mathcal{L}'_1,\mathcal{L}'_2,\dots,\mathcal{L}'_{k})$ with $\mathcal{L}_j\in C([\ell_j,\ell_1 ]\cup [r_1,r_j],\mathbb{R})$. Notice that there is the restriction map form $\mathfrak{X}_f$ to $\mathfrak{X}'$. Let $\bP'$ and $\tilde{\bP}'$ be the push-forward probability measure of $\bP$ and $\tilde{\bP}$ on $\mathfrak{X}'$ respectively. We write $\mathbb{E}'$ and $\tilde{\mathbb{E}}'$ for the expectations for $\bP'$ and $\tilde{\bP}'$ respectively. From \eqref{equ:PoverP}, we have
$$\frac{\textup{d}\mathbb{P}'}{\textup{d}\tilde{\mathbb{P}}'}(\cL') =\frac{1}{Z'} Z_{\mathbf{H}_t}^{1,k,(-L,L),\vec{x} ,\vec{y},+\infty,f_{k+1}}.$$
Here $\vec{x}=(\cL'_j(-L))_{j=1}^k$, $\vec{y}=(\cL'_j(L))_{j=1}^k$ and $Z'$ is a normalizing constant,
$$ Z'=\tilde{\mathbb{E}}'\left[ Z^{1,k,(-L,L),\vec{x} ,\vec{y} ,+\infty,f_{k+1}}_{\mathbf{H}_t} \right].  $$
From Corollary \ref{cor:Z_lowerbound-k}, 
$$Z'=\tilde{\mathbb{E}}\left[ Z^{1,k,(-L,L),\vec{x} ,\vec{y} ,+\infty,f_{k+1}}_{\mathbf{H}_t} \right]\geq  D_4^{-2}e^{-D_4 (L^{-1}M^2 + L)-2kL}=:\delta.$$ 
Thus
\begin{align*}
\mathbb{P}\left( Z^{1,k,(-L,L), \vec{x}, \vec{y},+\infty,f_{k+1}}_{\mathbf{H}_t}\leq \varepsilon \delta   \right)&=\mathbb{P}'\left( Z^{1,k,(-L,L),\vec{x},\vec{y},+\infty,f_{k+1}}_{\mathbf{H}_t}\leq \varepsilon \delta  \right)\\ 
&=\tilde{\mathbb{E}}' \left[ \frac{1}{Z'} Z_{\mathbf{H}_t}^{1,k,(-L,L),\vec{x} ,\vec{y},+\infty,f_{k+1}}\cdot\mathbbm{1}\left\{Z^{1,k,(-L,L),\vec{x},\vec{y},+\infty,f_{k+1}}_{\mathbf{H}_t}\leq \varepsilon \delta\right\}\right]\\
&\leq\ \frac{\varepsilon\delta}{Z'} \cdot \tilde{\mathbb{E}}'\left[\mathbbm{1}\left\{Z^{1,k,(-L,L),\vec{x},\vec{y},+\infty,f_{k+1}}_{\mathbf{H}_t}\leq \varepsilon \delta\right\}  \right]\leq \varepsilon.
\end{align*}
Thus the assertion follows by picking $D_5(k)=\max\{D^2_4, D_4+2k\}$.
\end{proof}

Now we are ready to prove Proposition \ref{pro:Z_lowerbound_k}.
\begin{proof}[Proof of Proposition \ref{pro:Z_lowerbound_k}]
Fix $K\geq L^2$ and let $M>0$ be a large number to be determined soon. Throughout this proof, we write $\bP_{\mathfrak{H}^t}$ for the law of the scaled KPZ line ensemble. Let $\mathsf{GoodBoundary}$ (shorthanded as $\mathsf{GB}$)  be the the event that the $M$-{\bf Good} boundary conditions holds, i.e.
\begin{align*}
\mathsf{GoodBoundary}:=\left\{ \sup_{u\in [ \ell_1 ,r_1]}|\mathfrak{H}^t_{k+1}(u)|\leq M \right\}\cap\bigcap_{j=1}^k \left\{ \sup_{u\in [\ell_1,\ell_j]\cup [r_j,r_1]  }|\mathfrak{H}^t_{j}(u)|\leq  M \right\}.
\end{align*}
From  Corollary~\ref{cor:tail} , by taking $M=D'K$ with suitably large $D'\geq 1$, we have
$$\mathbb{P}_{\mathfrak{H}^t}(\mathsf{GoodBoundary}^{\textup{c}})\leq 2^{-1}e^{-K^{3/2}}.$$
Moreover, \eqref{Mbound} holds. Let $\mathcal{F}^*_{ {ext}}$ be the sigma-field generated by $\mathfrak{H}^t$ restricted on
\begin{align*}
\mathbb{N}\times \mathbb{R}\setminus \bigcup_{j=1}^k \{j\}\times (\ell_j,r_j).
\end{align*}
By the Gibbs Property (see Definition \ref{def:H-BGP} and Corollary~\ref{cor:KPZGibbs}),
\begin{align*}
   \mathbb{E}_{\mathfrak{H}^t}\left[\mathbbm{1}\{ Z^{1,k,(-L,L), \vec{x},\vec{y} ,+\infty,\mathfrak{H}^t_{k+1}}_{\mathbf{H}_t}\leq 2^{-1}\varepsilon D_5^{-1}  e^{-D_5 (L^{-1}M^2+L)}\} \,|\, \mathcal{F}^*_{ {ext}}\right]
\end{align*}
with $\vec{x}=(\mathfrak{H}^t_j(-L))_{j=1}^k$ and $\vec{y}=(\mathfrak{H}^t_j(L))_{j=1}^k$ equals
\begin{align*}
 \mathbb{P}\left(Z^{1,k,(-L,L), \vec{x},\vec{y} ,+\infty,f_{k+1}}_{\mathbf{H}_t}\leq 2^{-1}\varepsilon D_5^{-1}  e^{-D_5 (L^{-1}M^2+L)}\right),
\end{align*} 
where $\mathbb{P}$ is defined in \eqref{resample-k} with $f_j=\mathfrak{H}^t_j$ in $[\ell_1,\ell_j]\cup [r_j,r_1],\ j\in [1,k]_{\mathbb{Z}} $ and $f_{k+1}=\mathfrak{H}^t_{k+1}$ in $[\ell_1,r_1]$.

Applying Proposition \ref{pro:Z_small-k} with $\varepsilon=2^{-1}e^{-K^{3/2}}$, we deduce that
\begin{align*}
 \mathbbm{1}_{\mathsf{GB}}\cdot\mathbb{E}_{\mathfrak{H}^t}\left[\mathbbm{1}\{ Z^{1,k,(\ell,r), \vec{x},\vec{y} ,+\infty,f_{k+1}}_{\mathbf{H}_t}\leq 2^{-1}\varepsilon D_5^{-1}  e^{-D_5 (L^{-1}M^2+L)}\} \,|\, \mathcal{F}^*_{ {ext}}\right] \leq 2^{-1}e^{-K^{3/2}}\cdot \mathbbm{1}_{\mathsf{GB}}.
\end{align*}
Thus we have
\begin{align*}
&\mathbb{P}_{\mathfrak{H}^t}\left( Z^{1,k,(-L,L), \vec{x},\vec{y} ,+\infty,f_{k+1}}_{\mathbf{H}_t}\leq 2^{-1}\varepsilon D_5^{-1}  e^{-D_5 (L^{-1}M^2+L)}  \right)\\
= & \mathbb{E}_{\mathfrak{H}^t}\left[\mathbbm{1}_{ \mathsf{GB}}\cdot\mathbbm{1}\left\{  Z^{1,k,(-L,L), \vec{x} ,\vec{y} ,+\infty,f_{k+1}}_{\mathbf{H}_t}\leq 2^{-1}  D_5^{-1}  e^{-K^{3/2}-D_5 (L^{-1}M^2+L)}\right\}   \right]+\mathbb{P}_{\mathfrak{H}^t}(\mathsf{GB}^{\textup{c}})\\
\leq &e^{-K^{3/2}}. 
\end{align*}

Because $K\geq L^2$ and $M=D'K$, we have $2^{-1}  D_5^{-1}  e^{-K^{3/2}-D_5 (L^{-1}M^2+L)}\geq D^{-1}e^{-DL^{-1}K^2}$. The proof is finished.
\end{proof}

The rest of this section is devoted to proving Proposition \ref{pro:spacing-k}. We will run a two-step inductive resampling. For $j\in [1,k]_{\mathbb{Z}}$, consider a sequence of events
\begin{align}
\mathsf{A}_j:=  \bigg\{ \inf_{u\in [\ell_{j+1},r_{j+1}]   } {\mathcal{L}}_j(u)\geq (4k-4j+4)M \bigg\}.
\end{align}
The first step is carried out inductively in Lemma \ref{lem:lowerbound-k} to raise curves and hence give a lower bound of $\tilde{\mathbb{P}}\bigg(\bigcap_{j=1}^k \mathsf{A}_j \bigg)$. Denote
\begin{align}
\mathsf{F}_j:=&\bigg\{  {\mathcal{L}}_j \in [(4k-4j+3)M,(4k-4j+5)M]\ \textup{in}\ [\ell_{j+1} , r_{j+1}]  \bigg\}\\
& \cap \bigg\{  {\mathcal{L}}_j \leq (4k-4j+5)M\ \textup{in}\ [\ell_{j} ,r_{j}] \bigg\}.\nonumber
\end{align}
The second step is carried out inductively to lower curve properly in order to separate them in the desired region. See Lemma \ref{lem:spacing-k}, which gives a lower bound of $\tilde{\mathbb{P}}\bigg(\bigcap_{j=1}^k \mathsf{F}_j \bigg)$. Proposition \ref{pro:spacing-k} follows directly from the Lemma \ref{lem:spacing-k} as $\displaystyle\bigcap_{j=1}^k\mathsf{F}_j \subset \mathsf{E}$. 

\begin{lemma}\label{lem:lowerbound-k}
Fix $k\in\mathbb{N}$, $L\geq 1$, $t\geq 1$ and $M\geq L^{1/2}$. Let $f$ be a $M$-{\bf Good} boundary condition and $\tilde{\mathbb{P}}$ be the probability measure on $\mathfrak{X}_f$ defined in \eqref{resample-kt}. Then there exists a constant $D_6= D_6(k)$ depending only on $k$ such that   	
\begin{align*}
\tilde{\mathbb{P}}\bigg(\bigcap_{j=1}^k \mathsf{A}_j \bigg)\geq D_6^{-1} e^{-D_6(L^{-1}M^2 + L)}. 
\end{align*}
\end{lemma}

\begin{figure}[h!]
\includegraphics[scale=0.32]{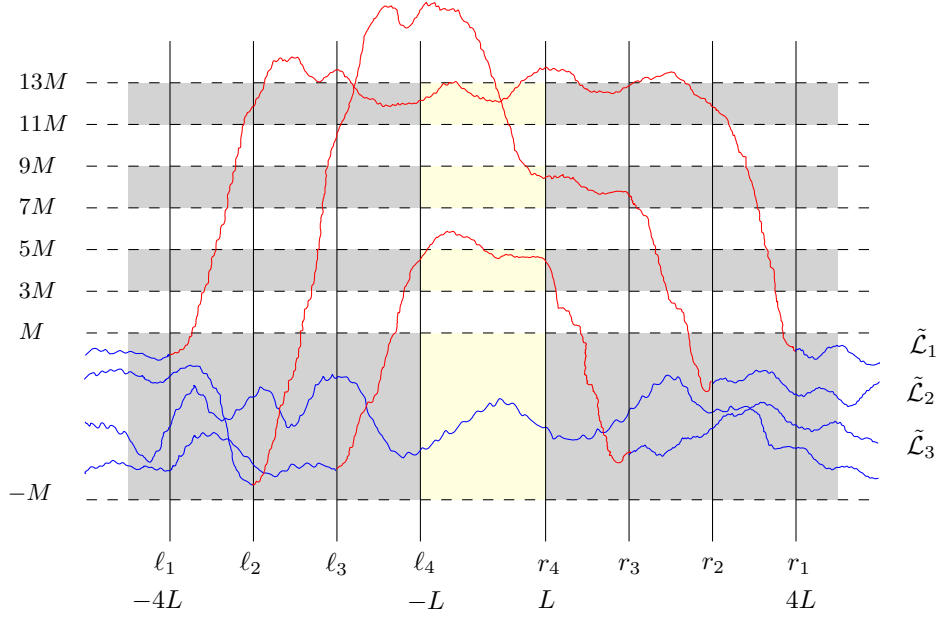}
\caption{A realization of $\bigcap_{j=1}^k \mathsf{A}_j$ with $k=3$. This is an illustration of the step one resampling. We first resample $ {\mathcal{L}}_1$ to have it rise above $12M$ within $[\ell_1, \ell_2]$ and $[r_1, r_2]$ and also have it stay above $12M$ after they reach it. Then resample $ {\mathcal{L}}_2$ to have it rise above $8M$ within $[\ell_2, \ell_3]$ and $[r_2, r_3]$ and also have it stay above $8M$ after they reach it. Furthermore resample $ {\mathcal{L}}_3$ to have it rise above $4M$ within $[\ell_3, \ell_4]$ and $[r_3, r_4]$ and also have it stay above $4M$ after they reach it.}
\end{figure}

\begin{proof}
We start by showing the lower bound for $\tilde{\mathbb{P}}(\mathsf{A}_1)$. From \eqref{resample-kt} and the Gibbs property (see Definition \ref{def:H-BGP}),  
\begin{align*}
\tilde{\mathbb{E}} [\mathbbm{1}_{\mathsf{A}_1} \,|\,\mathcal{F}_{ {ext}}(\{1\}\times (\ell_1,r_1))]= {\mathbb{P}_{\mathbf{H}_t,(-L,L)}^{1,1,(\ell_1,r_1),f_1(\ell_1),f_1(r_1),+\infty, \mathcal{{L}}_{2,f}} ( {\mathsf{A}_1})}
\end{align*}  
with $\mathcal{{L}}_{2,f}$ defined by \eqref{equ:extension}. Note that equivalently we have $${\mathbb{P}_{\mathbf{H}_t,(-L,L)}^{1,1,(\ell_1,r_1),f_1(\ell_1),f_1(r_1),+\infty, \mathcal{{L}}_{2,f}}}={\mathbb{P}_{\mathbf{H}_t }^{1,1,(\ell_1,r_1),f_1(\ell_1),f_1(r_1),+\infty, g_2}}$$
with
\begin{align*}
g_2(u)\coloneqq\left\{\begin{array}{cc}
\mathcal{{L}}_{2,f}(u) & u\in [\ell_1,-L)\cup(L,r_1],\\
-\infty & u\in [-L,L].
\end{array} \right.
\end{align*}
Using the stochastic monotonicity, Lemma~\ref{monotonicity}, we have
\begin{align*}
{\mathbb{P}_{\mathbf{H}_t}^{1,1,(\ell_1,r_1),f_1(\ell_1),f_1(r_1),+\infty, g_2} ( {\mathsf{A}_1})}&\geq{\mathbb{P}_{\mathbf{H}_t}^{1,1,(\ell_1,r_1),f_1(\ell_1),f_1(r_1),+\infty, -\infty} ( {\mathsf{A}_1})}\\
&=\mathbb{P}_{\free}^{1,1,(\ell_1,r_1),f_1(\ell_1),f_1(r_1) }\left( {\mathsf{A}_1} \right). 
\end{align*}

Let
\begin{align*}
a_1\coloneqq \inf\left\{\mathbb{P}_{\free}^{1,1,(\ell_1,r_1),x,y }\left( {\mathsf{A}_1} \right)\,\Big|\, |x|\leq M,\ |y|\leq M \right\}.
\end{align*}
Under $M$-{\bf Good} boundary conditions, $M\geq L^{1/2}$, we seek for a lower bound of $a_1$. To realize $\mathsf{A}_1$, it suffices to have the Brownian bridge $B_1$ first jump over a height $(4k+2)M$ within the interval $[\ell_1  ,\ell_2]$,  secondly remain above $4kM$ within the interval $[\ell_2 ,r_2]$, and thirdly drop down a height $(4k+2)M$ within $[r_2, r_1]$. $a_1$ is bounded below by the Brownian kinetic cost of such trajectory. Note that $|\ell_2-\ell_1|=|r_2-r_1|=L$ and $|\ell_2-r_2|=2kL$. Hence there exists $D=D(k)$ such that $a_1\geq D^{-1}e^{-DL^{-1}M^2}.$ Thus we deduce
\begin{align*}
\tilde{\mathbb{P}}(\mathsf{A}_1)&=  \tilde{\mathbb{E}}\left[ \tilde{\mathbb{E}} [\mathbbm{1}_{\mathsf{A}_1} \,|\,\mathcal{F}_{ {ext}}(\{1\}\times (\ell_1,r_1))] \right]\\
& \geq \mathbb{P}_{\free}^{1,1,(\ell_1,r_1),f_1(\ell_1),f_1(r_1) }\left( {\mathsf{A}_1} \right)  \geq D^{-1}e^{-DL^{-1}M^2}.
\end{align*}

Now we proceed by induction on $j$. Assume for $j\in [2,k]_{\Z}$, we have
$$\tilde{\mathbb{P}}\bigg(\bigcap_{i=1}^{j-1} \mathsf{A}_{i} \bigg)\geq D^{-1} e^{-D (L^{-1}M^2 +L)}.$$ 
We aim to show 
$$\tilde{\mathbb{P}}\bigg(\bigcap_{i=1}^{j}\mathsf{A}_{i} \bigg)\geq D^{-1}e^{-D (L^{-1}M^2 +L)}.$$ 

From \eqref{resample-kt} and the Gibbs property (see Definition \ref{def:H-BGP}), 
\begin{align*}
\tilde{\mathbb{E}} [\mathbbm{1}_{\mathsf{A}_j} \,|\,\mathcal{F}_{ {ext}}(\{j\}\times (\ell_j,r_j))]=      \mathbb{P}^{j,j,(\ell_j,r_j),f_j(\ell_j),f_j(r_j), {\mathcal{L}}_{j-1 }, {\mathcal{L}}_{j+1,f} }_{\mathbf{H}_t,(-L,L)}(\mathsf{A}_j),    
\end{align*}
with $ {\mathcal{L}}_{j+1,f} $ defined in \eqref{equ:extension} and we adopt the convention that $ {\mathcal{L}}_{k+1,f}=f_{k+1}$. Note that
\begin{align*}
\mathbb{P}^{j,j,(\ell_j,r_j),f_j(\ell_j),f_j(r_j), {\mathcal{L}}_{j-1 }, {\mathcal{L}}_{j+1,f} }_{\mathbf{H}_t,(-L,L)}=\mathbb{P}^{j,j,(\ell_j,r_j),f_j(\ell_j),f_j(r_j),g_{j-1 },g_{j+1} }_{\mathbf{H}_t }
\end{align*}
with
\begin{align*}
g_{j-1}(u)\coloneqq\left\{\begin{array}{cc}
\mathcal{ {L}}_{j-1}(u) & u\in [\ell_j,-L)\cup(L,r_j],\\
+\infty & u\in [-L,L],\end{array} \right.
\end{align*}
and
\begin{align*}
g_{j+1}(u)\coloneqq\left\{\begin{array}{cc}
\mathcal{ {L}}_{j+1,f }(u) & u\in [\ell_j,-L)\cup(L,r_j],\\
-\infty & u\in [-L,L].
\end{array} \right.
\end{align*}
We claim that 
\begin{equation}\label{equ:01181110}
\tilde{\mathbb{E}}[\mathbbm{1}_{\mathsf{A}_j}\,|\,\mathcal{F}_{ {ext}}(\{j\}\times (\ell_j,r_j)) ]\geq  \mathbb{E}_{\free}^{j,j,(\ell_j,r_j),f_j(\ell_j),f_j(r_j) }\left[\mathbbm{1}_{\mathsf{A}_j} \cdot W^{j,j,(\ell_j,r_j),f_j(\ell_j),f_j(r_j),g_{j-1 },-\infty }_{\mathbf{H}_t } \right]. 
\end{equation}
This can be derived as
\begin{align*}
\tilde{\mathbb{E}}[\mathbbm{1}_{\mathsf{A}_j}\,|\,\mathcal{F}_{ {ext}}(\{j\}\times (\ell_j,r_j)) ]&=\mathbb{P}^{j,j,(\ell_j,r_j),f_j(\ell_j),f_j(r_j),g_{j-1 },g_{j+1} }_{\mathbf{H}_t }(\mathsf{A}_j)\\
&\geq \mathbb{P}^{j,j,(\ell_j,r_j),f_j(\ell_j),f_j(r_j),g_{j-1 },-\infty }_{\mathbf{H}_t } (\mathsf{A}_j)\\
&=\frac{1}{Z_{\mathbf{H}_t }}\mathbb{E}_{\free}^{j,j,(\ell_j,r_j),f_j(\ell_j),f_j(r_j) }\left[\mathbbm{1}_{\mathsf{A}_j} \cdot W^{j,j,(\ell_j,r_j),f_j(\ell_j),f_j(r_j),g_{j-1 },-\infty }_{\mathbf{H}_t } \right].\\
&\geq \mathbb{E}_{\free}^{j,j,(\ell_j,r_j),f_j(\ell_j),f_j(r_j) }\left[\mathbbm{1}_{\mathsf{A}_j} \cdot W^{j,j,(\ell_j,r_j),f_j(\ell_j),f_j(r_j),g_{j-1 },-\infty }_{\mathbf{H}_t } \right]. 
\end{align*}
Here the in the second equality we use $Z_{\mathbf{H}_t }$ to abbreviate 
$$Z_{\mathbf{H}_t }:=\mathbb{E}_{\free}^{j,j,(\ell_j,r_j),f_j(\ell_j),f_j(r_j) }\left[  W^{j,j,(\ell_j,r_j),f_j(\ell_j),f_j(r_j),g_{j-1 },-\infty }_{\mathbf{H}_t } \right].$$ 
In the first inequality we apply stochastic monotonicity, Lemma~\ref{monotonicity}, and in the second inequality we use the fact that normalizing constant is bounded from above by 1.
 
Now we proceed to find a lower bound for
$$\mathbbm{1}_{\mathsf{A}_{j-1}}\cdot \mathbb{E}_{\free}^{j,j,(\ell_j,r_j),f_j(\ell_j),f_j(r_j) }\left[\mathbbm{1}_{\mathsf{A}_j} \cdot W^{j,j,(\ell_j,r_j),f_j(\ell_j),f_j(r_j),g_{j-1 },-\infty }_{\mathbf{H}_t } \right].$$
Consider the event
\begin{align*}
\mathsf{D}_j:= &\bigg\{ \inf_{u\in [\ell_{j+1},  r_{j+1}]  } {\mathcal{L}}_j(u)\geq (4k-4j+4)M \bigg\} \cap \bigg\{ \sup_{u\in [\ell_{j} ,r_{j}]  }  {\mathcal{L}}_j(u)\leq (4k-4j+8)M \bigg\}.
\end{align*}
Note that $\mathsf{D}_j\subset \mathsf{A}_j$. As $\mathsf{D}_j$ and $\mathsf{A}_{j-1}$ occur, $ {\mathcal{L}}_{j-1}\geq  {\mathcal{L}}_{j}$ in $[\ell_{j} ,r_{j}]$ and hence
\begin{align*}
W^{j,j,(\ell_j,r_j),f_j(\ell_j),f_j(r_j),g_{j-1 },-\infty }_{\mathbf{H}_t }\cdot\mathbbm{1}_{\mathsf{D}_j}\cdot\mathbbm{1}_{\mathsf{A}_{j-1}}  \geq e^{-2(k+1-j)L}\cdot\mathbbm{1}_{\mathsf{D}_j}\cdot\mathbbm{1}_{\mathsf{A}_{j-1}}.
\end{align*}
Consequently
\begin{equation*}
\begin{split}
 & \mathbbm{1}_{\mathsf{A}_{j-1}} \mathbb{E}_{\free}^{j,j,(\ell_j,r_j),f_j(\ell_j),f_j(r_j) }\left[\mathbbm{1}_{\mathsf{A}_j} \cdot W^{j,j,(\ell_j,r_j),f_j(\ell_j),f_j(r_j),g_{j-1 },-\infty }_{\mathbf{H}_t } \right]\\
&\geq \mathbbm{1}_{\mathsf{A}_{j-1}} \mathbb{E}_{\free}^{j,j,(\ell_j,r_j),f_j(\ell_j),f_j(r_j) }\left[\mathbbm{1}_{\mathsf{D}_j} \cdot W^{j,j,(\ell_j,r_j),f_j(\ell_j),f_j(r_j),g_{j-1 },-\infty }_{\mathbf{H}_t } \right]\\
&\geq \mathbbm{1}_{\mathsf{A}_{j-1}}\mathbb{E}_{\free}^{j,j,(\ell_j,r_j),f_j(\ell_j),f_j(r_j) }\left[\mathbbm{1}_{\mathsf{D}_j} \cdot e^{-2(k+1-j)L}  \right]\\
&\geq \mathbbm{1}_{\mathsf{A}_{j-1}}\cdot e^{-2(k+1-j)L} \mathbb{P}_{\free}^{j,j,(\ell_j,r_j),x,y } (\mathsf{D}_j).
\end{split}
\end{equation*}
In short, we derived
\begin{equation}\label{equ:01191118}
\begin{split}
 & \mathbbm{1}_{\mathsf{A}_{j-1}} \mathbb{E}_{\free}^{j,j,(\ell_j,r_j),f_j(\ell_j),f_j(r_j) }\left[\mathbbm{1}_{\mathsf{A}_j} \cdot W^{j,j,(\ell_j,r_j),f_j(\ell_j),f_j(r_j),g_{j-1 },-\infty }_{\mathbf{H}_t } \right]\\
\geq & \mathbbm{1}_{\mathsf{A}_{j-1}}\cdot e^{-2(k+1-j)L} \mathbb{P}_{\free}^{j,j,(\ell_j,r_j),x,y } (\mathsf{D}_j).
\end{split}
\end{equation}
From Lemma~\ref{lem:BBjump}, the $M$-{\bf Good} boundary condition implies that there exists $D=D(k)$ such that
\begin{align}\label{def:a_j}
a_j&:=  \inf \left\{\mathbb{P}_{\free}^{j,j,(\ell_j,r_j),x,y} (\mathsf{D}_j)\,|\, |x|\leq M,\ |y|\leq M\right\} \geq D^{-1} e^{-DL^{-1} M^2 }.
\end{align}
From \eqref{equ:01181110}, \eqref{equ:01191118}, \eqref{def:a_j} and the induction hypothesis, we have
\begin{align*}
\tilde{\mathbb{P}}\bigg( \bigcap_{i=1}^{j } \mathsf{A}_{i} \bigg)=&\tilde{\mathbb{E}}\left[\prod_{i=1}^{j-1}\mathbbm{1}_{\mathsf{A}_i}\cdot \tilde{\mathbb{E}}[\mathbbm{1}_{\mathsf{A}_j}\,|\,\mathcal{F}_{ {ext}}(\{j\}\times (\ell_j,r_j)) ] \right]\\
&\geq a_j\cdot e^{-2(k+1-j)L}\cdot \tilde{\mathbb{P}}\bigg( \bigcap_{i=1}^{j-1 } \mathsf{A}_{i} \bigg)\\
&\geq D^{-1}e^{-D(L^{-1}M^2+L)}.
\end{align*}
This completes the induction argument and hence proves the desired result.
\end{proof}

For $j\in [1,k]_{\mathbb{Z}}$, recall the events
\begin{align*}
\mathsf{F}_j:=&\bigg\{  {\mathcal{L}}_j \in [(4k-4j+3)M,(4k-4j+5)M]\ \textup{in}\ [\ell_{j+1} , r_{j+1}]  \bigg\}\\
& \cap \bigg\{  {\mathcal{L}}_j \leq (4k-4j+5)M\ \textup{in}\ [\ell_{j} ,r_{j}] \bigg\}.
\end{align*}

In the following Lemma \ref{lem:spacing-k}, we give a lower bound for $\tilde{\mathbb{P}}\bigg(\bigcap_{j=1}^k \mathsf{F}_j \bigg)$.

\begin{lemma}\label{lem:spacing-k}
Fix $k\in\mathbb{N}$, $L\geq 1$, $t\geq 1$ and $M\geq L^{1/2}$. Let $f$ be a $M$-{\bf Good} boundary condition and $\tilde{\mathbb{P}}$ be the probability measure on $\mathfrak{X}_f$ defined in \eqref{resample-kt}. Then there exists a constant  $D_7=D_7(k)$ such that  
\begin{align*}
\tilde{\mathbb{P}}\bigg(\bigcap_{j=1}^k \mathsf{F}_j \bigg)\geq D_7^{-1}e^{-D_7 (L^{-1}M^2 +L)}.
\end{align*}
\end{lemma}

\begin{figure}[h!]
\includegraphics[width=12cm]{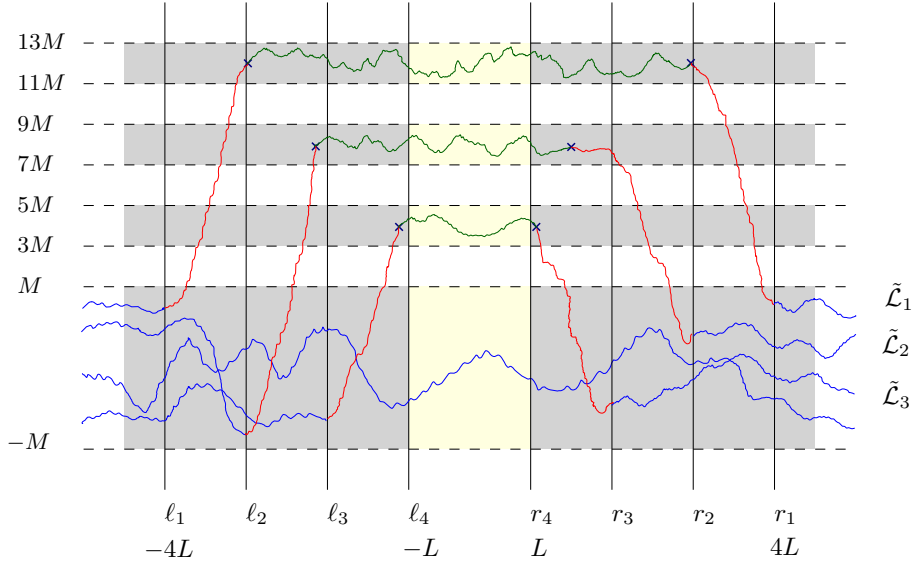}
\caption{A realization of $\bigcap_{j=1}^k \mathsf{F}_j$ with $k=3$. This is an illustration of the step two resampling. First for $ {\mathcal{L}}_3$, find the $\{3\}$-stopping domain $(\mathfrak{l}_3, \mathfrak{r}_3)$ where it first hits $4M$ from left and right respectively. Then resample $ {\mathcal{L}}_3$ on $[\mathfrak{l}_3, \mathfrak{r}_3]$ to have it stay within $[3M, 5M]$. Then for $ {\mathcal{L}}_2$, find the $\{2\}$-stopping domain $(\mathfrak{l}_2,\mathfrak{r}_2)$ where it first hits $8M$ from left and right respectively. Then resample $ {\mathcal{L}}_2$ on $[\mathfrak{l}_2, \mathfrak{r}_2]$ to have it stay within $[7M, 9M]$. Lastly for $ {\mathcal{L}}_1$, find the $\{1\}$-stopping domain $(\mathfrak{l}_1,\mathfrak{r}_1)$ where it first hits $12M$ from left and right respectively. Then resample $ {\mathcal{L}}_1$ on $[\mathfrak{l}_1, \mathfrak{r}_1]$ to have it stay within $[11M, 13M]$.}
\end{figure}

\begin{proof}
We will run a resampling in a reversed order starting from the $k$-layer and argue inductively. More precisely, we start by showing a lower bound for  $\tilde{\mathbb{P}}\bigg(\bigcap_{j=1}^{k-1} \mathsf{A}_j\cap \mathsf{F}_k \bigg)$. 

Let $[\mathfrak{l}_k, \mathfrak{r}_k]$ be a $\{k\}$-stopping domain such that
\begin{align*}
\mathfrak{l}_k:=&\sup \big\{u\in [\ell_k ,\ell_{k+1}] \ |\  {\mathcal{L}}_{k}(u')<  4M\ \textup{for all}\ u'\in [\ell_k,u]\big\},\\
\mathfrak{r}_k:=&\inf\big\{u\in [r_{k+1},r_k] \ |\  {\mathcal{L}}_{k}(u')< 4M\ \textup{for all}\ u'\in [u,r_k]\big\},
\end{align*}
where we set $\mathfrak{l}_k=\ell_k$ or $\mathfrak{r}_k=r_k$ if the sets are empty respectively. Consider the event
\begin{align*}
\mathsf{F}'_k:=\bigg\{  {\mathcal{L}}_k(\mathfrak{l}_k)= {\mathcal{L}}_k(\mathfrak{r}_k)=4M \bigg\}.
\end{align*}

Because $ {\mathcal{L}}_k(\ell_k), {\mathcal{L}}_k(r_k)\in [-M,M]$,  we have $ \mathsf{A}_k\subset\mathsf{F}'_k$. In the view of Lemma \ref{lem:lowerbound-k},
\begin{equation}\label{equ:inducttion}
\tilde{\mathbb{P}}\bigg(\bigcap_{j=1}^{k-1} \mathsf{A}_j \cap \mathsf{F}'_k \bigg)\geq D^{-1}e^{-D(L^{-1}M^2+L)}.  
\end{equation}
We would like to have $ {\mathcal{L}}_k$ stay in the preferable region $[3M,5M]$ over $[\mathfrak{l}_k,\mathfrak{r}_k]$. Let
\begin{align*}
\mathsf{F}_k''\coloneqq\{ {\mathcal{L}}_k \in [3M,5M]\ \textup{in}\ [\mathfrak{l}_k,\mathfrak{r}_k ]\}.
\end{align*}
Note that the occurrence of $\mathsf{A}_{k-1}\cap\mathsf{F}''_{k}$ implies ordering between $ {\mathcal{L}}_{k-1}$ and $ {\mathcal{L}}_{k}$, i.e. 
\begin{align*}
\mathsf{A}_{k-1}\cap\mathsf{F}''_{k}\subset \{  {\mathcal{L}}_{k-1}\geq  {\mathcal{L} }_k\geq f_{k+1}\ \textup{in}\  [\mathfrak{l}_k, \mathfrak{r}_k ]\}.
\end{align*}
Hence
\begin{align*}
W^{k,k,(\mathfrak{l}_k ,\mathfrak{r}_k),\mathcal{\tilde{L}}_k(\mathfrak{l}_k ),\mathcal{\tilde{L}}_k(\mathfrak{r}_k ), {\mathcal{L}}_{k-1},f_{k+1} }_{\mathbf{H}_t,(-L,L)}\cdot \mathbbm{1}_{\mathsf{A}_{k-1} }\cdot \mathbbm{1}_{\mathsf{F}_{k}'' } \geq e^{-4L}\cdot \mathbbm{1}_{\mathsf{A}_{k-1} }\cdot \mathbbm{1}_{\mathsf{F}_{k}'' } . 
\end{align*} 
And
\begin{equation}\label{equ:01181127}
\begin{split}
&\quad\ \mathbbm{1}_{\mathsf{A}_{k-1}}\cdot \mathbbm{1}_{\mathsf{F}'_{k}}\cdot \tilde{\mathbb{E}}[\mathbbm{1}_{\mathsf{F}''_k}\,|\,\mathcal{F}_{ {ext}}(\{k\}\times (\mathfrak{l}_k,\mathfrak{r}_k))]\\
&=\mathbbm{1}_{\mathsf{A}_{k-1}}\cdot \mathbbm{1}_{\mathsf{F}'_{k}}\cdot\mathbb{P}^{k,k,(\mathfrak{l}_k ,\mathfrak{r}_k),\mathcal{\tilde{L}}_k(\mathfrak{l}_k ),\mathcal{\tilde{L}}_k(\mathfrak{r}_k ), {\mathcal{L}}_{k-1},f_{k+1} }_{\mathbf{H}_t,(-L,L)}\left( \mathsf{F}_k''  \right) \\
& \geq \mathbbm{1}_{\mathsf{A}_{k-1}}\cdot \mathbbm{1}_{\mathsf{F}'_{k}}\cdot{\mathbb{E}^{k,k,(\mathfrak{l}_k ,\mathfrak{r}_k),\mathcal{\tilde{L}}_k(\mathfrak{l}_k ),\mathcal{\tilde{L}}_k(\mathfrak{r}_k ) }_{\free}\left[\mathbbm{1}_{\mathsf{F}_k''}\cdot W^{k,k,(\mathfrak{l}_k ,\mathfrak{r}_k),\mathcal{\tilde{L}}_k(\mathfrak{l}_k ),\mathcal{\tilde{L}}_k(\mathfrak{r}_k ), {\mathcal{L}}_{k-1},f_{k+1} }_{\mathbf{H}_t,(-L,L)} \right] } \\
&\geq \mathbbm{1}_{\mathsf{A}_{k-1}}\cdot \mathbbm{1}_{\mathsf{F}'_{k}} \cdot e^{-4L} \cdot \mathbb{P}^{k,k,(\mathfrak{l}_k ,\mathfrak{r}_k),4M,4M   }_{\free} (\mathsf{F}''_k).
\end{split}
\end{equation}
Under $M$-{\bf Good} boundary conditions, there exists $D=D(k)$ such that
\begin{align}\label{def:b_k}
b_k\coloneqq \inf\left\{ \mathbb{P}^{k,k,(\ell ,r),4M,4M   }_{\free} \left(\cL_k \in [3M,5M] \right) \right\}\geq D^{-1}.
\end{align}
The infimum is taken over all $\ell\in [\ell_k,\ell_{k+1}]$ and $ r\in [r_{k+1},r_k]$. Moreover
\begin{align*}
\bigcap_{j=1}^{k-1} \mathsf{A}_j \cap \mathsf{F}'_k\cap\mathsf{F}_k''\subset \bigcap_{j=1}^{k-1} \mathsf{A}_j \cap \mathsf{F}_k.
\end{align*}
Combining \eqref{equ:01181127}, \eqref{def:b_k} and \eqref{equ:inducttion}, we derive
\begin{align*}
\tilde{\mathbb{P}}  \bigg(  \bigcap_{j=1}^{k-1} \mathsf{A}_j \cap   \mathsf{F}_k\bigg)\geq & \tilde{\mathbb{P}}  \bigg(  \bigcap_{j=1}^{k-1} \mathsf{A}_j \cap \mathsf{F}'_k\cap \mathsf{F}''_k\bigg)= \tilde{\mathbb{E}}\left[\prod_{j=1}^{k-1}\mathbbm{1}_{\mathsf{A}_j}\cdot\mathbbm{1}_{\mathsf{F}'_{k }}\cdot \tilde{\mathbb{E}}[\mathbbm{1}_{\mathsf{F}''_k}\,|\,\mathcal{F}_{ {ext}}(\{k\}\times (\mathfrak{l}_k,\mathfrak{r}_k))] \right]\\
\geq & b_k\cdot e^{-4L}\cdot \tilde{\mathbb{P}}  \bigg(  \bigcap_{j=1}^{k-1} \mathsf{A}_j \cap \mathsf{F}'_k \bigg) \geq  D^{-1}e^{-D(L^{-1}M^2+L)}.
\end{align*}

We now proceed by a reversed induction. Assume for some $i\in [1,k-1]_{\Z}$, we have 
\begin{align*}
\tilde{\mathbb{P}}\bigg(\bigcap_{j=1}^{i} \mathsf{A}_j\cap \bigcap_{j=i+1}^{k} \mathsf{F}_j \bigg)\geq D^{-1}e^{-D(L^{-1}M^2+L)}.
\end{align*}
We aim to show
\begin{align*}
\tilde{\mathbb{P}}\bigg(\bigcap_{j=1}^{i-1} \mathsf{A}_j\cap \bigcap_{j=i}^{k} \mathsf{F}_j \bigg)\geq  D^{-1}e^{-D(L^{-1}M^2+L)}
\end{align*}
and we adopt the convention that $\bigcap_{j=1}^{0} \mathsf{A}_j$ means the total probability space. 

Let $[\mathfrak{l}_{i}, \mathfrak{r}_{i}]$ be a $\{i\}$-stopping domain such that
\begin{align*}
\mathfrak{l}_{i}:=&\sup\big\{u\in [\ell_i ,\ell_{i+1}] \ |\  {\mathcal{L}}_{i}(u')< (4k-4i+4) M\ \textup{for all}\ u'\in [\ell_i,u]\big\},\\
\mathfrak{r}_{i}:=&\inf\big\{u\in [r_{i+1},r_{i }] \ |\  {\mathcal{L}}_{i}(u')< (4k-4i+4) M\ \textup{for all}\ u'\in [u,r_i]\big\},
\end{align*}
where we set $\mathfrak{l}_i=\ell_i$ or $\mathfrak{r}_i=r_i$ if the set is empty respectively. Consider the event
\begin{align*}
\mathsf{F}'_i:=\bigg\{  {\mathcal{L}}_i(\mathfrak{l}_i)= {\mathcal{L}}_i(\mathfrak{r}_i)=(4k-4i+4)M \bigg\}.
\end{align*}
Because $ {\mathcal{L}}_i(\ell_i), {\mathcal{L}}_i(r_i)\in [-M,M]$, we have $ \mathsf{A}_i\subset\mathsf{F}'_i$. We would like to have $ {\mathcal{L}}_i$ stay in the preferable region $[(4k-4j+3)M,(4k-4j+5)M]$ over $[\mathfrak{l}_i,\mathfrak{r}_i]$.
Let
\begin{align*}
\mathsf{F}_i''\coloneqq\{ \cL_i \in [(4k-4i+3)M,(4k-4i+5)M]\ \textup{in}\ [\mathfrak{l}_i,\mathfrak{r}_i ]\}.
\end{align*}
Note that the occurrence of $\mathsf{A}_{i-1}\cap\mathsf{F}''_{i}\cap \mathsf{F}_{i+1}$ implies ordering between $ {\mathcal{L}}_{i-1}$ ,$ {\mathcal{L}}_{i}$ and $ {\mathcal{L}}_{i+1,f}$, i.e. 
\begin{align*}
\mathsf{A}_{i-1}\cap\mathsf{F}''_{i}\cap \mathsf{F}_{i+1} \subset \{  {\mathcal{L}}_{i-1}\geq  {\mathcal{L} }_i\geq  {\mathcal{L}}_{i+1,f}\ \textup{in}\  [\mathfrak{l}_i, \mathfrak{r}_i ]\}.
\end{align*}
Hence
\begin{align*}
W^{i,i,(\mathfrak{l}_i ,\mathfrak{r}_i),\mathcal{\tilde{L}}_i(\mathfrak{l}_i ),\mathcal{\tilde{L}}_i(\mathfrak{r}_i ), {\mathcal{L}}_{i-1}, {\mathcal{L}}_{i+1,f}  }_{\mathbf{H}_t,(-L,L)}\cdot \mathbbm{1}_{\mathsf{A}_{i-1} }\cdot \mathbbm{1}_{\mathsf{F}_{i}'' }\cdot \mathbbm{1}_{\mathsf{F}_{i+1}  } \geq e^{-4L}\cdot \mathbbm{1}_{\mathsf{A}_{i-1} }\cdot \mathbbm{1}_{\mathsf{F}_{i}'' }\cdot \mathbbm{1}_{\mathsf{F}_{i+1}  } . 
\end{align*} 
And
\begin{equation}\label{equ:01181134}
\begin{split}
&\quad\ \mathbbm{1}_{\mathsf{A}_{i-1}}\cdot \mathbbm{1}_{\mathsf{F}'_{i}}\cdot \mathbbm{1}_{\mathsf{F}_{i+1}  } \cdot \tilde{\mathbb{E}}[\mathbbm{1}_{\mathsf{F}''_i}\,|\,\mathcal{F}_{ {ext}}(\{i\}\times (\mathfrak{l}_i,\mathfrak{r}_i))]\\
&=\mathbbm{1}_{\mathsf{A}_{i-1}}\cdot \mathbbm{1}_{\mathsf{F}'_{i}}\cdot \mathbbm{1}_{\mathsf{F}_{i+1}  }\cdot\mathbb{P}^{i,i,(\mathfrak{l}_i ,\mathfrak{r}_i),\mathcal{\tilde{L}}_i(\mathfrak{l}_i ),\mathcal{\tilde{L}}_i(\mathfrak{r}_i ), {\mathcal{L}}_{i-1}, {\mathcal{L}}_{i+1,f}  }_{\mathbf{H}_t,(-L,L)}\left( \mathsf{F}_i''  \right) \\
& \geq \mathbbm{1}_{\mathsf{A}_{i-1}}\cdot \mathbbm{1}_{\mathsf{F}'_{i}}\cdot \mathbbm{1}_{\mathsf{F}_{i+1}  }\cdot{\mathbb{E}^{i,i,(\mathfrak{l}_i ,\mathfrak{r}_i),\mathcal{\tilde{L}}_i(\mathfrak{l}_i ),\mathcal{\tilde{L}}_i(\mathfrak{r}_i ) }_{\free}\left[\mathbbm{1}_{\mathsf{F}_i''}\cdot W^{i,i,(\mathfrak{l}_i ,\mathfrak{r}_i),\mathcal{\tilde{L}}_i(\mathfrak{l}_i ),\mathcal{\tilde{L}}_i(\mathfrak{r}_i ), {\mathcal{L}}_{i-1}, {\mathcal{L}}_{i+1,f}  }_{\mathbf{H}_t,(-L,L)} \right] } \\
&\geq \mathbbm{1}_{\mathsf{A}_{i-1}}\cdot \mathbbm{1}_{\mathsf{F}'_{i}}\cdot \mathbbm{1}_{\mathsf{F}_{i+1}  } \cdot e^{-4L} \cdot \mathbb{P}^{i,i,(\mathfrak{l}_i ,\mathfrak{r}_i),(4k-4i+4)M,(4k-4i+4)M   }_{\free} (\mathsf{F}''_i).
\end{split}
\end{equation}
It is straightforward to check that there exists $D=D(k)$ such that
\begin{align}\label{def:b_i}
b_i\coloneqq \inf\left\{ \mathbb{P}^{i,i,(\ell ,r),(4k-4i+4)M,(4k-4i+4)M   }_{\free} \left(B_i \in [(4k-4i+3)M,(4k-4i+5)M] \right)\right\}\geq D^{-1}.
\end{align}
The infimum is taken over all $\ell\in [\ell_i,\ell_{i+1}]$ and $r\in [r_{i+1},r_i]$. Moreover $  \mathsf{F}'_i\cap\mathsf{F}_i''\subset  \mathsf{F}_i.$ 

Combining \eqref{equ:01181134}, \eqref{def:b_i} and the induction hypothesis, we get
\begin{align*}
\tilde{\mathbb{P}}\bigg(\bigcap_{j=1}^{i-1} \mathsf{A}_j\cap \bigcap_{j=i}^{k} \mathsf{F}_j \bigg)\geq & \tilde{\mathbb{P}}\bigg(\bigcap_{j=1}^{i-1} \mathsf{A}_j  \cap \bigcap_{j=i}^{k} \mathsf{F}_j \cap\mathsf{F}_i'\cap\mathsf{F}_i''\bigg)\\
=&  \tilde{\mathbb{E}}\left[\prod_{j=1}^{i-1}\mathbbm{1}_{\mathsf{A}_j}\cdot  \prod_{j=i+1}^k \mathbbm{1}_{\mathsf{F}_j}\cdot\mathbbm{1}_{\mathsf{F}'_{i }} \cdot \tilde{\mathbb{E}}[\mathbbm{1}_{\mathsf{F}''_i}\,|\,\mathcal{F}_{ {ext}}(\{i\}\times (\mathfrak{l}_i,\mathfrak{r}_i))] \right]\\
\geq & b_i\cdot e^{-4(k+1-i)L}\cdot \tilde{\mathbb{P}}\bigg(\bigcap_{j=1}^{i-1} \mathsf{A}_j \cap \bigcap_{j=i+1}^{k} \mathsf{F}_j \cap\mathsf{F}_i' \bigg)\\
\geq & b_i\cdot e^{-4(k+1-i)L}\cdot \tilde{\mathbb{P}}\bigg(\bigcap_{j=1}^{i} \mathsf{A}_j \cap \bigcap_{j=i+1}^{k} \mathsf{F}_j  \bigg)\\
\geq &D^{-1}e^{-D(L^{-1}M^2+L)}.
\end{align*}
We used $\mathsf{A}_i\subset \mathsf{F}'_i$ in the second to last inequality and used the induction hypothesis to the last inequality. The induction argument is finished.
\end{proof}

\begin{appendix}
\section{Results on Brownian Bridges}\label{sec:Brownian}
We record in this section some properties of Brownian bridges that we need. Recall that we denote $\mathbb{P}^{1,1, (a,b),x,y}_{\free}$ for the law of a Brownian bridge $B$ defined on $[a,b]$ with $B(a)=u$ and $B(b)=v$. The next lemma is an analogue of \cite[Chapter 4 (3.40)]{KS}, where the supremum is considered.
\begin{lemma}\label{lem:BB-max}
\begin{align*}
\mathbb{P}^{1,1, (a,b),x,y}_{\free} \left( \inf_{u\in [a,b]}B(u)\leq \beta \right)= e^{-2(b-a)^{-1}(x-\min\{ \beta, x, y\})(y-\min\{ \beta, x, y\})}.
\end{align*}
\end{lemma}
\begin{proof}
By \cite[Chapter 4 (3.40)]{KS}, we have that
\begin{align*}
\mathbb{P}^{1,1, (a,b),x,y}_{\free} \left( \sup_{u\in [a,b] }B(u)\geq \beta \right)= e^{-2(b-a)^{-1}(\max\{ \beta, x, y\}-x)(\max\{ \beta, x, y\}-y)}.
\end{align*}
Assume $B$ has the law law $\mathbb{P}^{1,1, (a,b),x,y}_{\free} $, the law of $\tilde{B}\coloneqq -B$ is $\mathbb{P}^{1,1,(a,b),-x,-y}_{\free}$. Hence
\begin{align*}
\mathbb{P}^{1,1, (a,b),x,y}_{\free} \left( \inf_{u\in [a,b]}B(u)\leq \beta \right)=&\mathbb{P}^{1,1,(a,b),-x,-y}_{\free}\left( \sup_{u\in [a,b]}\tilde{B}(u)\geq -\beta \right)\\
=&e^{-2(b-a)^{-1}(\max\{-\beta,-x,-y\}+x)(\max\{-\beta,-x,-y\}+y) }\\
= & e^{-2(b-a)^{-1}(x-\min\{ \beta, x, y\})(y-\min\{ \beta, x, y\})}.
\end{align*}
\end{proof}
The following two lemmas can be found in \cite[Lemma 5.13]{Ham1}. We present the proofs for completeness. 
\begin{lemma}\label{lem:BB-osc}
There exists a constant $C_0$ such that for all $d\in (0,1]$ and $K\geq 0$,
\begin{align*}
\mathbb{P}^{1,1,(0,1),0,0}_{\free} \left( \sup_{u,v\in [0,1],\ |u-v|\leq d}|B(u)-B(v)|>Kd^{1/2} \right)\leq d^{-1} C_0e^{-C_0^{-1}K^2}.
\end{align*}
\end{lemma}
\begin{lemma}\label{lem:BB-pre}
There exists a constant $C_1$ such that for all $\ell<r$ in  $[0,1]$ and $K\geq 0$,
\begin{align*}
\mathbb{P}^{1,1,(0,1),0,0}_{\free}  \left( \sup_{u,v\in [\ell,r]  }|B(u)-B(v)|>K(r-\ell)^{1/2}  \right)<C_1e^{-C_1^{-1}K^2}.
\end{align*}
\end{lemma}
\begin{proof}[Proof of Lemma \ref{lem:BB-osc}]
Suppose $d\in [4^{-1},1]$, the assertion of Lemma \ref{lem:BB-osc} follows easily by applying Lemma \ref{lem:BB-pre} with $\ell=0$ and $r=1$. From now on we assume $d\in (0,4^{-1})$. Let $m=\lfloor d^{-1}\rfloor$. For $j\in [1,m]_{\mathbb{Z}}$, define
\begin{align*}
\mathsf{E}_j\coloneqq\left\{ \sup_{u,v\in [(j-1)d,(j+1)d]\cap [0,1]  }|B(u)-B(v)|>Kd^{1/2} \right\}.
\end{align*}
As $$\left\{ \sup_{u,v\in [\ell,r]  }|B(u)-B(v)|>K(r-\ell)^{1/2} \right\}\subset \bigcup_{i=1}^{m}\, \mathsf{E}_j,$$
by Lemma \ref{lem:BB-pre}, there exists a universal constant $C_0$ such that
\begin{align*}
\mathbb{P}^{1,1,(0,1),0,0}_{\free}  \left( \sup_{u,v\in [0,1],\ |u-v|\leq d}|B(u)-B(v)|>Kd^{1/2} \right)\leq \sum_{j=1}^m \mathbb{P}^{1,1,(0,1),0,0}_{\free}   (\mathsf{E}_j)\leq d^{-1}C_0e^{-C_0^{-1}K^2}.
\end{align*}
The proof is finished.
\end{proof}
Given $\sigma>0$ and $a<b$ with $a,b \in [-\infty,\infty]$, let
\begin{align*}
\nu_{\sigma^2}( a,b )\coloneqq \int_a^b (2\pi\sigma^2)^{-1/2}e^{-2^{-1}\sigma^{-2} u^2}\, dx.
\end{align*}
In other words, $\nu_{\sigma^2}( a,b )$ is the probability that a Gaussian random variable with mean $0$ and variance $\sigma^2$ is contained in $[a,b]$. In particular, $\nu_{\sigma^2}( a,b )=\nu_{1}( \sigma^{-1} a,\sigma^{-1} b ).$ From \cite[Section 14.8]{Wil}, for all $a>0$,
\begin{equation}\label{normal}
\nu_1(a,\infty)\leq (2\pi)^{-1/2}\cdot a^{-1}e^{-2^{-1}a^2}.
\end{equation}
\begin{proof}[Proof of Lemma \ref{lem:BB-pre}]
For notational simplicity, we denote $\mathbb{P}^{1,1,(0,1),0,0}_{\free}$ by $\mathbb{P}_{\free}$. Suppose $K\leq 1$, by choosing $C_1$ large, we can get and $C_1e^{-C_1^{-1}K^2}\geq 1$ and the assertion holds trivially. From now on we assume $K\geq 1$. Let $X(t)$ be a standard Brownian motion. The argument is based on the following property,
\begin{align*}
B(u)\overset{(d)}{=}(1-u)X\left( \frac{u}{1-u} \right). 
\end{align*}

Suppose first $r-\ell\geq 2^{-1}.$ Then
\begin{align*}
&\mathbb{P}_{\free} \left( \sup_{u,v\in [\ell,r]  }|B(u)-B(v)|>K(r-\ell)^{1/2}  \right)\leq \mathbb{P}_{\free} \left( \sup_{u,v\in [0,1]  }|B(u)-B(v)|>2^{-1/2}K  \right)\\
\leq & \mathbb{P}_{\free} \left( \sup_{u \in [0,1 ]  }|B(u) |>2^{-3/2}K  \right) 
\leq   2\mathbb{P}_{\free} \left( \sup_{u \in [0,2^{-1} ]  }|B(u) |>2^{-3/2}K  \right)\\
= & 2\mathbb{P} \left( \sup_{u \in [0,2^{-1}  ]  }(1-u) \left|X\left(\frac{u}{1-u} \right)\right|>2^{-3/2}K  \right)\\
\leq &2\mathbb{P} \left( \sup_{u' \in [0,1  ]  } \left|X\left(u'\right)\right|>2^{-3/2}K  \right)= 4\nu_1(2^{-3/2}K,\infty).
\end{align*}
Here we used the reflection principle for the last equality. Because $K\geq 1$, from \eqref{normal} there exists a constant $C$ such that
\begin{align*}
\mathbb{P}_{\free} \left( \sup_{u,v\in [\ell,r]  }|B(u)-B(v)|>K(r-\ell)^{1/2}  \right)\leq Ce^{-C^{-1}K^2}.
\end{align*}

From now on, we assume $r-\ell<2^{-1}$. Because $B(u)\overset{(d)}{=} B(1-u)$, we can assume without loss of generality that $[r,\ell]\subset [0,\frac{3}{4}]$. Given $u,v\in [\ell,r]$
\begin{align*}
|B(u)-B(v)|\overset{(d)}{=}&\left|(1-u)X\left(\frac{u}{1-u}\right)-(1-v)X\left(\frac{v}{1-v}\right) \right|\\
\leq &(1-u)\left|X\left(\frac{u}{1-u}\right)-X\left(\frac{v}{1-v}\right) \right|+|u-v|\left|X\left(\frac{v}{1-v}\right)\right|\\
\leq & \left|X\left(\frac{u}{1-u}\right)-X\left(\frac{v}{1-v}\right) \right|+(r-\ell)\left|X\left(\frac{v}{1-v}\right)\right|.
\end{align*}    
Thus, by setting $\ell'=\frac{\ell}{1-\ell}$ and $r'=\frac{r}{1-r}$,
\begin{align*}
&\mathbb{P}_{\free}  \left( \sup_{u,v\in [\ell,r]  }|B(u)-B(v)|>K(r-\ell)^{1/2}  \right)\\
\leq &\mathbb{P}\left( \sup_{u',v'\in [\ell',r']  }|X(u')-X(v')|>2^{-1}K(r-\ell)^{1/2}  \right) +\mathbb{P}\left( \sup_{u' \in [\ell',r']  }|X(u') |>2^{-1}K(r-\ell)^{-1/2}  \right).  
\end{align*}
As $\frac{d}{dt}\left(\frac{t}{1-t} \right)=\frac{1}{(1-t)^2}\leq 16$ for $t\in [0,\frac{3}{4}]$, we have $r'-\ell'\leq 16(r-\ell)$. Together with $r'\leq 3$ and $r-\ell\leq 2^{-1}$, we deduce
\begin{align*}
&\mathbb{P}_{\free}  \left( \sup_{u,v\in [\ell,r]  }|B(u)-B(v)|>K(r-\ell)^{1/2}  \right)\\
\leq &\mathbb{P}\left( \sup_{u' \in [\ell',r']  }|X(u')|>4^{-1}K(r-\ell)^{1/2}  \right) +\mathbb{P}\left( \sup_{u' \in [0,3]  }|X(u') |>2^{-1}K(r-\ell)^{-1/2}  \right) \\
=&2\nu_1(4^{-1}(r-\ell)^{1/2}(r'-\ell')^{-1/2}K,\infty)+2\nu_1(2^{-1}(r-\ell)^{-1/2}3^{-1/2}K,\infty)\\
\leq &2\nu_1( K/16,\infty)+2\nu_1(6^{-1/2}K,\infty)
\end{align*}
Because $K\geq 1$, in the view of \eqref{normal}, we conclude that there exists a universal constant $C$ such that
\begin{align*}
\mathbb{P}_{\free} \left( \sup_{u,v\in [\ell,r]  }|B(u)-B(v)|>K(r-\ell)^{1/2}  \right)\leq Ce^{-C^{-1}K^2}.
\end{align*}
The proof is finished.
\end{proof}
\begin{lemma}\label{lem:BBjump}
Fix $k\in\mathbb{N}$, $L\geq 1$ and $M\geq L^{1/2}$. There exists a constant $D_8=D_8(k)$ depending only on $k$ such that the following statement holds. Let $\ell,r,\lambda,x$ and $y$ be numbers that satisfy $4L\leq r-\ell\leq (2k+2)L$, $4\leq \lambda\leq 4k$ and $|x|,|y|\leq M$. Let
\begin{align*}
\mathsf{J}\coloneqq \left\{  \inf_{u\in [\ell+L,r-L]} B(u)\geq \lambda M, \sup_{u\in [\ell ,r ]} B(u)\leq (\lambda+4) M    \right\}.
\end{align*} 
Then we have
\begin{align*}
\mathbb{P}_{\free}^{1,1,(\ell,r),x,y}\left(\mathsf{J} \right)\geq D_8e^{-D_8^{-1}L^{-1}M^2}.
\end{align*}
\end{lemma}
\begin{proof}
Consider the following events
\begin{align*}
\mathsf{J}_1\coloneqq&\left\{ B(\ell+L),B(r-L)\in [(\lambda+1)M,(\lambda+3)M] \right\},\\
\mathsf{J}_2\coloneqq&\left\{ \sup_{u\in [\ell,\ell+L]}\left| B(u)-L^{-1}((u-\ell)B(\ell+L)+(\ell+L-u)u) \right| \leq M \right\},\\
\mathsf{J}_3\coloneqq&\left\{ \sup_{u\in [\ell+L,r-L]}\left| B(u)-(r-\ell-2L)^{-1}((u-\ell-L)B(r- L)+(r-L-u)B(\ell+L)) \right| \leq M \right\},\\
\mathsf{J}_4\coloneqq&\left\{ \sup_{u\in [r-L,r]}\left| B(u)-L^{-1}((u-r+L)y+(r-u)B(r-L)) \right| \leq M \right\}. 
\end{align*}
The event $\mathsf{J}_1$ implies that $B(\ell+L)$ and $B(r-L)$ are contained in $[(\lambda+1)M,(\lambda+3)M]$. The events $\mathsf{J}_2$, $\mathsf{J}_3$ and $\mathsf{J}_4$ imply respectively that $B(u)$ is deviated from the linear interpolation on the intervals $[\ell,\ell+L]$, $[\ell+L,r-L]$ and $[r-L,r]$ at most by $M$. It is straightforward to check that $\cap_{i=1}^4\mathsf{J}_i \subset \mathsf{J}$. Furthermore, $\mathsf{J}_i$ are independent. Hence it suffices to bound $\mathbb{P}_{\free}^{1,1,(\ell,r),x,y}(\mathsf{J}_i)$ from below. We start with $\mathsf{J}_2$. 
\begin{align*}
\mathbb{P}_{\free}^{1,1,(\ell,r),x,y}(\mathsf{J}_2)=&\mathbb{P}_{\free}^{1,1,(0,L),0,0}(\sup_{u\in[0,L]}|B(u)|\leq M ) = \mathbb{P}_{\free}^{1,1,(0,1),0,0}(\sup_{u\in[0,1]}|B(u)|\leq L^{-1/2}M )\\
\geq &\mathbb{P}_{\free}^{1,1,(0,1),0,0}(\sup_{u\in[0,1]}|B(u)|\leq 1 ) 
\geq   C^{-1}.
\end{align*}
We have used $M\geq L^{1/2}$. Similarly, $\mathbb{P}_{\free}^{1,1,(\ell,r),x,y}(\mathsf{J}_4)\geq C^{-1}$ and 
\begin{align*}
\mathbb{P}_{\free}^{1,1,(\ell,r),x,y}(\mathsf{J}_3)=&\mathbb{P}_{\free}^{1,1,(0,1),0,0}(\sup_{u\in[0,1]}|B(u)|\leq (r-\ell-2L)^{-1/2}M )\\
\geq &\mathbb{P}_{\free}^{1,1,(0,1),0,0}(\sup_{u\in[0,1]}|B(u)|\leq (2k)^{-1/2} )\geq D^{-1}.
\end{align*}

It remains to deal with $\mathsf{J}_1$. Define
\begin{align*}
\mathsf{J}_{1,-}\coloneqq&\left\{ B(\ell+L)\in [(\lambda+1)M,(\lambda+3)M] \right\}\\
\mathsf{J}_{1,+}\coloneqq&\left\{  B(r-L)\in [(\lambda+1)M,(\lambda+3)M] \right\}.
\end{align*}
Under the law $\mathbb{P}_{\free}^{1,1,(\ell,r),x,y}$, $B(\ell+L)$ is a Gaussian random variable with mean and variance $$m_-=(\ell-r)^{-1}( Ly+(r-\ell-L)x ),\ \sigma_-^2=(r-\ell)^{-1}L(r-\ell-L).$$ Therefore
\begin{align*}
\mathbb{P}_{\free}^{1,1,(\ell,r),x,y}(\mathsf{J}_{1,-})=\nu_1(\sigma_-^{-1}((\lambda+1)M-m_-),\sigma_-^{-1}((\lambda+3)M-m_-)).
\end{align*}
From the assumption, $|m_-|\leq M$ and $\frac{3}{4}L\leq \sigma_-^2\leq L $. Hence
\begin{align*}
\sigma_-^{-1}((\lambda+1)M-m_-)\leq  2(\lambda+2)\cdot 3^{-1/2}L^{-1/2}M\leq 3^{-1/2}(8k+4)L^{-1/2}M,
\end{align*}
and
\begin{align*}
\sigma_-^{-1}((\lambda+3)M-m_-)-\sigma_-^{-1}((\lambda+1)M-m_-)=2\sigma^{-1}M\geq 2L^{-1/2}M\geq 2.
\end{align*}
Therefore,
\begin{align*}
\mathbb{P}_{\free}^{1,1,(\ell,r),x,y}(\mathsf{J}_{1,-})=&\nu_1(\sigma_-^{-1}((\lambda+1)M-m_-),\sigma_-^{-1}((\lambda+3)M-m_-))\\
\geq& \nu_1(  3^{-1/2}(8k+4)L^{-1/2}M ,3^{-1/2}(8k+4)L^{-1/2}M+2 )\\
\geq& D^{-1}e^{-DL^{-1}M^2}.
\end{align*}
Conditioned on the event $\mathsf{J}_{1,-}$, $B(r-L)$ is a Gaussian distribution with mean and variance
\begin{align*}
m_+=(r-\ell-L)^{-1}((r-\ell-2L)y+LB(\ell+L)),\ \sigma_+^{2}=(r-\ell-L)^{-1}L(r-\ell-2L).
\end{align*}
A similar argument yields
\begin{align*}
\mathbb{P}_{\free}^{1,1,(\ell,r),x,y}(\mathsf{J}_{1,+}\, |\, \mathsf{J}_{1,-})\geq  D^{-1}e^{-DL^{-1}M^2}.
\end{align*}
This implies $$\mathbb{P}_{\free}^{1,1,(\ell,r),x,y}(\mathsf{J}_{1})=\mathbb{P}_{\free}^{1,1,(\ell,r),x,y}(\mathsf{J}_{1,-}) \mathbb{P}_{\free}^{1,1,(\ell,r),x,y}(\mathsf{J}_{1,+}\, |\, \mathsf{J}_{1,-})\geq D^{-1}e^{-DL^{-1}M^2}.$$
The assertion then follows by combining the bounds on $\mathbb{P}_{\free}^{1,1,(\ell,r),x,y}(\mathsf{J}_i)$. 
\end{proof}
\section{Tail Bounds}\label{sec:tail}
In this section we prove quantitative tail estimates, Propositions \ref{pro:lowertail} and \ref{pro:uppertail} for the scaled KPZ line ensemble $\mathfrak{H}^t$ (defined in Definition~\ref{def:scaledKPZLE}). These two propositions are used as key inputs for Proposition \ref{pro:Z_lowerbound_k}, a quantitative estimate on the normalizing constant.
\begin{proposition}\label{pro:lowertail}
Fix $k\in \mathbb{N}$. There exists a constant  $D_9=D_9(k)$ depending only on $k$ such that for all $t\geq 1$ and $r\geq 0$,
\begin{align}\label{lowertail_k}
\mathbb{P}\left( \inf_{u\in [-1,1]}\left(\mathfrak{H}^t_k(u)+2^{-1} {u^2}   \right)\leq -r \right)\leq D_9 e^{-D_9^{-1}r^{3/2}}.
\end{align}
\end{proposition}
\begin{proposition}\label{pro:uppertail}
Fix $k\in \mathbb{N}$. There exists a constant  $D_{10}=D_{10}(k)$ depending only on $k$ such that for all for all $t\geq 1$ and $r\geq 0$,
\begin{align}\label{uppertail_k}
\mathbb{P}\left( \sup_{u\in [-1,1]}\left(\mathfrak{H}^t_k(u)+2^{-1} {u^2}   \right)\geq r \right)\leq D_{10}e^{-D_{10}^{-1}r^{3/2}} .
\end{align}
\end{proposition}
Since $\mathfrak{H}^t_k(u)+2^{-1}u^2 $ is stationary in $u$, we have the following corollary.
\begin{corollary}\label{cor:tail}
Let $I\subset \mathbb{R}$ be an interval and $|I|$ be the length of $I$. Then for all $k\geq 1$, $t\geq 1$ and $r\geq 0$,
\begin{align}\label{uppertail_k_I}
\mathbb{P}\left( \inf_{u\in I}\left(\mathfrak{H}^t_k(u)+2^{-1} {u^2}   \right)\leq -r \right)\leq (|I|/2+1) D_9 e^{-D_9^{-1}r^{3/2}}.
\end{align}
\begin{align}\label{lowertail_k_I}
\mathbb{P}\left( \sup_{u\in I}\left(\mathfrak{H}^t_k(u)+2^{-1} {u^2}   \right)\geq r \right)\leq  (|I|/2+1)D_{10} e^{-D_{10}^{-1}r^{3/2}}.
\end{align}
\end{corollary}

We will run induction on $k$. The case $k=1$ follows the tail bounds for the solution to KPZ equation with the narrow wedge initial condition \cite{CG1, CG2}. In the rest of this section, we consider $k\geq 2$ and assume Propositions \ref{pro:lowertail} and \ref{pro:uppertail} hold for $1,2,\dots k-1$. In particular, Corollary \ref{cor:tail} holds for for $1,2,\dots k-1$.
\subsection{Proof of Proposition \ref{pro:lowertail}}
For $\tau>0$, define the events
\begin{align*}
\mathsf{Low}^{[-4\tau,-2\tau]}_k:=\left\{ \sup_{u\in [-4\tau,-2\tau]}\left(\mathfrak{H}^t_k(u)+ 2^{-1} u^2 \right) \leq -4^{-1} \tau^2-\tau^{1/2} \right\},\\
\mathsf{Low}^{[2\tau,4\tau]}_k:=\left\{ \sup_{u\in [2\tau,4\tau]}\left(\mathfrak{H}^t_k(u)+2^{-1} u^2\right) \leq -4^{-1} \tau^2-\tau^{1/2} \right\}.
\end{align*}
\begin{lemma}\label{lem:Lowk}
There exist $\tau_0>0$ and $D_{11}=D_{11}(k)$ depending only on $k$ such that for all $t\geq 1$ and $\tau\geq \tau_0$,
\begin{align*}
\mathbb{P}(\mathsf{Low}^{[-4\tau,-2\tau]}_k),\ \mathbb{P}(\mathsf{Low}^{[2\tau,4\tau]}_k)\leq D_{11} e^{-D_{11}^{-1} \tau^{3}}.
\end{align*}	
\end{lemma}
\begin{proof}
Because of the stationarity of $ \mathfrak{H}^t_k(u)+2^{-1}u^2$, it suffices to prove 
$$\mathbb{P}(\mathsf{Low}^{[-\tau, \tau]}_k)\leq \bar{E}_k e^{-\bar{c}_k\tau^{3 }} $$
with
\begin{align*}
\mathsf{Low}^{[-\tau,\tau]}_k:=\left\{ \sup_{u\in [-\tau,\tau]}\left(\mathfrak{H}^t_k(u)+2^{-1} {u^2} \right) \leq -4^{-1} \tau^2-\tau^{1/2} \right\}.
\end{align*}
Let
\begin{align*}
\mathsf{G}_{k-1}:=\left\{ \mathfrak{H}_{k-1}(\tau)+2^{-1} \tau^2 \leq 4^{-1}\tau^2 \right\}\cap\left\{ \mathfrak{H}_{k-1}(-\tau)+2^{-1} \tau^2  \leq 4^{-1}\tau^2 \right\}
\end{align*}
and 
$$\mathsf{A}_k:=\mathsf{Low}^{[-\tau,\tau]}_k\cap \mathsf{G}_{k-1}.$$
Applying Proposition~\ref{pro:uppertail} for $k-1$, there exists $D=D(k)$ such that
\begin{align}\label{equ:01230305}
\bP(\mathsf{G}^{\textup{c}}_{k-1})\leq De^{-D^{-1}\tau^3}.
\end{align}

Next, we bound $\bP(\{\mathfrak{H}^t_{k-1}(0)\geq -8^{-1}\tau^2\}\cap\mathsf{A}_{k})$ from above.
\begin{align}\label{equ:01230334-1}
\bP(\{\mathfrak{H}^t_{k-1}(0)\geq -8^{-1}\tau^2\}\cap\mathsf{A}_{k})=\mathbb{E}\left[\mathbbm{1}_{\mathsf{A}_k}\mathbb{E}[\mathbbm{1}\{\mathfrak{H}^t_{k-1}(0)\geq -8^{-1}\tau^2\}  |  \mathcal{F}_{ext}(\{k-1\}\times (-\tau,\tau)) ]\right].
\end{align}
By the Gibbs property of the scaled KPZ line ensemble (see Corollary~\ref{cor:KPZGibbs} and Definition~\ref{def:H_Brownian}),
\begin{equation}\label{equ:01230334-2}
\begin{split}
&\mathbb{E}[\mathbbm{1}\{\mathfrak{H}^t_{k-1}(0)\geq -8^{-1}\tau^2\}  |  \mathcal{F}_{ext}(\{k-1\}\times (-\tau,\tau)) ]\\
=&\mathbb{P}^{k-1,k-1,(-\tau,\tau),x,y,\mathfrak{H}^t_{k-2},\mathfrak{H}^t_{k}}_{\mathbf{H}_t}(\{\cL(0)\geq -8^{-1}\tau^2\}).
\end{split}
\end{equation}
Here $x=\mathfrak{H}^t_{k-1}(-\tau)$ and $y=\mathfrak{H}^t_{k-1}(\tau)$. We use the stochastic monotonicity to simply the boundary condition. From the definition of $\mathsf{A}_k$ and Lemma~\ref{monotonicity}, we have
\begin{equation}\label{equ:01230334-3}
\begin{split}
&\mathbbm{1}_{\mathsf{A}_k} \mathbb{P}^{k-1,k-1,(-\tau,\tau),x,y,\mathfrak{H}^t_{k-2},\mathfrak{H}^t_{k}}_{\mathbf{H}_t}(\{\cL(0)\geq 8^{-1}\tau^2\})\\
\leq &\mathbbm{1}_{\mathsf{A}_k} \mathbb{P}^{k-1,k-1,(-\tau,\tau),-\tau^2/4,-\tau^2/4,\infty,-u^2/2-\tau^2/4-\tau^{1/2}}_{\mathbf{H}_t}(\{\cL(0)\geq 8^{-1}\tau^2\}).
\end{split}
\end{equation}

We then look for an upper bound for $\mathbb{P}^{k-1,k-1,(-\tau,\tau),-\tau^2/4,-\tau^2/4,\infty,-u^2/2-\tau^2/4-\tau^{1/2}}_{\mathbf{H}_t}(\{\cL(0)\geq 8^{-1}\tau^2\})$. By a direct calculation,
\begin{equation}\label{equ:01230331}
\begin{split}
&\mathbb{P}^{k-1,k-1,(-\tau,\tau),-\tau^2/4,-\tau^2/4,\infty,-u^2/2-\tau^2/4-\tau^{1/2}}_{\free}(\{\cL(0)\geq -8^{-1}\tau^2\})\\
=&\mathbb{P}_{\free}^{1,1,(-1,1),0,0}\left(  \cL(0)\geq   -8^{-1}\tau^{3/2} \right)\leq Ce^{-C^{-1}\tau^3}.
\end{split}
\end{equation}
To bound the normalizing constant, consider the event 
\begin{align*}
\mathsf{H}:=\left\{\inf_{u\in [-\tau,\tau]} \cL(u)\geq -\tau^2/4-\tau^{1/2}\right\}.
\end{align*}
As $\mathsf{H}$ occurs, $\cL(u)\geq -u^2/2-\tau^2/4-\tau^{1/2}$ on $[-\tau,\tau]$. Hence
\begin{align*}
\mathbbm{1}_{\mathsf{H}} W^{k-1,k-1,(-\tau,\tau),-\tau^2/4,-\tau^2/4,\infty,-u^2/2-\tau^2/4-\tau^{1/2}}_{\mathbf{H}_t}(B)\geq e^{-2\tau}\mathbbm{1}_{\mathsf{H}}.
\end{align*}
Therefore,
\begin{equation}\label{equ:01230332}
\begin{split}
&Z^{k-1,k-1,(-\tau,\tau),-\tau^2/4,-\tau^2/4,\infty,-u^2/2-\tau^2/4-\tau^{1/2}}_{\mathbf{H}_t}\\
\geq & e^{-2\tau}\mathbb{P}_{\free}^{k-1,k-1,(-\tau,\tau),-\tau^2/4,-\tau^2/4} (\mathsf{H})\\
=&e^{2\tau} \mathbb{P}_{\free}^{1,1,(-1,1),0,0}\left( \inf_{u\in [-1,1]}\cL (u)\geq -1 \right)\geq C^{-1} e^{-2\tau}.
\end{split}
\end{equation}
Combining \eqref{equ:01230331} and \eqref{equ:01230332} and setting $\tau_0$ large enough, we have
\begin{equation}\label{equ:01230334}
\begin{split}
&\mathbb{P}^{k-1,k-1,(-\tau,\tau),-\tau^2/4,-\tau^2/4,\infty,-u^2/2-\tau^2/4-\tau^{1/2}}_{\mathbf{H}_t}(\{\cL(0)\geq -8^{-1}\tau^2\}) \leq  Ce^{2\tau-C^{-1}\tau^3}\leq 2^{-1}.
\end{split} 
\end{equation}
Combining \eqref{equ:01230334-1}, \eqref{equ:01230334-2}, \eqref{equ:01230334-3} and \eqref{equ:01230334}, we obtain
\begin{align*}
\bP(\{\mathfrak{H}^t_{k-1}(0)\geq -8^{-1}\tau^2\}\cap\mathsf{A}_{k})\leq 2^{-1}\bP(\mathsf{A}_k).
\end{align*}
Applying Proposition~\ref{pro:lowertail} for $k-1$, we get
\begin{align}\label{equ:01230347}
\bP(\mathsf{A}_k)\leq 2\bP(\{\mathfrak{H}^t_{k-1}(0)< -8^{-1}\tau^2\})\leq De^{-D^{-1}\tau^3}.
\end{align}
Combining \eqref{equ:01230305} and \eqref{equ:01230347}, we conclude
\begin{align*}
\mathbb{P}(\mathsf{Low}^{[-\tau,\tau]}_k)&\leq \mathbb{P}(\mathsf{A}_k)+\mathbb{P}(\mathsf{G}_{k-1}^{\textup{c}}) \leq De^{ -D^{-1}\tau^3} .
\end{align*}
\end{proof}
Now we are ready to prove Proposition \ref{pro:lowertail}. 
\begin{proof}[Proof of Proposition \ref{pro:lowertail}]
Let $\tau_0$ be the number in Lemma \ref{lem:Lowk} and  $\tau\geq \max\{2,\tau_0\}$. Define the events
\begin{align*}
\mathsf{Up}^{[-4\tau,4\tau]}_{k-1} =&\left\{ \inf_{u\in [-4\tau,4\tau]}\left( \mathfrak{H}^t_{k-1}(u)+2^{-1}u^2\right)\geq -2^{-1}  \tau^2   \right\},\\
\mathsf{E}=&\mathsf{Up}^{[-4\tau,4\tau]}_{k-1}\cap \left(\mathsf{Low}^{[-4\tau,-2\tau]}_k\right)^{\textup{c}}\cap \left(\mathsf{Low}^{[2\tau,4\tau]}_k\right)^{\textup{c}}. 
\end{align*}
By Lemma \ref{lem:Lowk} and Corollary \ref{cor:tail}, 
\begin{align}\label{equ:01230456}
\mathbb{P}(\mathsf{E}^{\textup{c}})\leq (4\tau+1)D e^{ -D^{-1}\tau^3}.
\end{align}
Define
\begin{align*}
\sigma_1&=\inf\left\{ u\in [-4\tau,-2\tau]\ |\mathfrak{H}^t_k(u)+2^{-1}u^2\geq - 4^{-1}\tau^2-\tau^{1/2} \right\},\\
 \sigma_2&=\sup\left\{ u\in [ 2\tau,4\tau]\ |\mathfrak{H}^t_k(u)+2^{-1}u^2\geq - 4^{-1}\tau^2-\tau^{1/2} \right\}.
\end{align*}
In the above sets are empty, we set $\sigma_1=-2\tau$ or $\sigma_2=2\tau$ respectively.

Next, we bound $\bP(\{\inf_{u\in [-2\tau,2\tau]}  \mathfrak{H}^t_{k}(u) \leq -10\tau^2 \}\cap\mathsf{E})$ from above. 
\begin{equation}\label{equ:01230435-1}
\begin{split}
&\bP\left(\left\{\inf_{u\in [-2\tau,2\tau]}  \mathfrak{H}^t_{k}(u) \leq -10\tau^2 \right\}\cap\mathsf{E}\right)\\
=&\mathbb{E}\left[\mathbbm{1}_{\mathsf{E}}\mathbb{E}[\mathbbm{1}\{\inf_{u\in [-2\tau,2\tau]}  \mathfrak{H}^t_{k}(u) \leq -10\tau^2 \}\,|\,\mathcal{F}_{ext}(\{k\}\times (\sigma_1,\sigma_2))]  \right].
\end{split} 
\end{equation}
By the strong Gibbs property of the scaled KPZ line ensemble (see Corollary~\ref{cor:KPZGibbs} and Lemma~\ref{lem:stronggibbs}), we have
\begin{equation}\label{equ:01230435-2}
\begin{split}
&\mathbb{E}[\mathbbm{1}\{\inf_{u\in [-2\tau,2\tau]}  \mathfrak{H}^t_{k}(u) \leq -10\tau^2 \}\,|\,\mathcal{F}_{ext}(\{k\}\times (\sigma_1,\sigma_2))]  \\
=&\mathbb{P}^{k,k,(\sigma_1,\sigma_2),x,y,\mathfrak{H}^t_{k-1},\mathfrak{H}^t_{k+1}}_{\mathbf{H}_t}\left(\inf_{u\in [-2\tau,2\tau]}  \cL (u) \leq -10\tau^2\right).
\end{split}
\end{equation}
Here $x=\mathfrak{H}^t_k(\sigma_1)$ and $y=\mathfrak{H}^t_k(\sigma_2)$. As $\mathsf{E}$ occurs, 
\begin{align*}
\mathfrak{H}^t_k(\sigma_1),\ \mathfrak{H}^t_k(\sigma_2) \geq -8\tau^2- 4^{-1}\tau^2-\tau^{1/2}\geq -9\tau^2.
\end{align*}
Here we used $\tau\geq 2$. By the stochastic monotonicity, Lemma~\ref{monotonicity}, we have
\begin{equation}\label{equ:01230435-3}
\begin{split}
&\mathbbm{1}_{\mathsf{E}} \mathbb{P}^{k,k,(\sigma_1,\sigma_2),x,y,\mathfrak{H}^t_{k-1},\mathfrak{H}^t_{k+1}}_{\mathbf{H}_t}\left(\inf_{u\in [-2\tau,2\tau]}  \cL (u) \leq -10\tau^2\right)\\
\leq &\mathbbm{1}_{\mathsf{E}} \mathbb{P}^{k,k,(\sigma_1,\sigma_2),-9\tau^2,-9\tau^2,-u^2-\tau^2 ,-\infty }_{\mathbf{H}_t}\left(\inf_{u\in [-2\tau,2\tau]}  \cL (u) \leq -10\tau^2\right).
\end{split}
\end{equation}
We then look for an upper bound of $\mathbb{P}^{k,k,(\sigma_1,\sigma_2),-9\tau^2,-9\tau^2,-u^2-\tau^2 ,-\infty }_{\mathbf{H}_t}(\inf_{u\in [-2\tau,2\tau]}  \cL (u) \leq -10\tau^2).$ From a direct computation,
\begin{equation}\label{equ:01230428}
\begin{split}
&\mathbb{P}^{k,k,(\sigma_1,\sigma_2),-9\tau^2,-9\tau^2}_{\free}\left(\inf_{u\in [-2\tau,2\tau]}  \cL (u) \leq -10\tau^2\right)\\
\leq &\mathbb{P}^{k,k,(\sigma_1,\sigma_2),-9\tau^2,-9\tau^2}_{\free}\left(\inf_{u\in [\sigma_1,\sigma_2]}  \cL (u) \leq -10\tau^2\right)\\
=&\mathbb{P}^{1,1,(-1,1),0,0}_{\free}\left(\inf_{u\in [-1,1]}  \cL (u) \leq -\tau^2/(2^{-1}(\sigma_2-\sigma_1))^{1/2}\right)\\
\leq &\mathbb{P}^{1,1,(-1,1),0,0}_{\free}\left(\inf_{u\in [-1,1]}  \cL (u) \leq -2^{-1}\tau^{3/2} \right)\leq C e^{-C^{-1}\tau^3}.
\end{split}
\end{equation}
To bound the normalizing constant, consider the event
\begin{align*}
\mathsf{H}:=\left\{\sup_{u\in [\sigma_1,\sigma_2]}\cL(u)\leq -(17/2)\tau^2\right\}.
\end{align*}
The event $\mathsf{H}$ implies $\cL$ is smaller than $-u^2/2-\tau^2/2$ in $[\sigma_1,\sigma_2]$. Therefore,
\begin{align*}
\mathbbm{1}_{\mathsf{H}}W^{k,k,(\sigma_1,\sigma_2),-9\tau^2,-9\tau^2,-u^2-\tau^2 ,-\infty }_{\mathbf{H}_t} \geq e^{-8\tau}\mathbbm{1}_{\mathsf{H}}.
\end{align*}
Taking expectation, we obtain
\begin{equation}\label{equ:01230429}
\begin{split}
&Z^{k,k,(\sigma_1,\sigma_2),-9\tau^2,-9\tau^2,-u^2-\tau^2 ,-\infty }_{\mathbf{H}_t}\\
\geq &e^{-8\tau}\mathbb{P}^{k,k,(\sigma_1,\sigma_2),-9\tau^2,-9\tau^2 }_{\free}(\mathsf{H})\\
=  &e^{-8\tau}\mathbb{P}^{1,1,(-1,1),0,0}_{\free}\left(\sup_{u\in [-1,1]}B(u)\leq	 2^{-1}  \tau^2 / (2^{-1}( \sigma_2-\sigma_1  ) )^{1/2} \right)\\
 \geq&e^{-8\tau}\mathbb{P}^{1,1,(-1,1),0,0}_{\free}\left(\sup_{u\in [-1,1]}B(u)\leq 4^{-1}{\tau^{3/2}}   \right)\geq C^{-1} e^{-8\tau}
\end{split}
\end{equation}
Combining \eqref{equ:01230428} and \eqref{equ:01230429}, we get
\begin{equation}\label{equ:01230435-4}
\mathbb{P}^{k,k,(\sigma_1,\sigma_2),-9\tau^2,-9\tau^2,-u^2-\tau^2 ,-\infty }_{\mathbf{H}_t}\left(\inf_{u\in [-2\tau,2\tau]}  \cL (u) \leq -10\tau^2\right)\leq Ce^{8\tau-C^{-1}\tau^3}.
\end{equation}
Combining \eqref{equ:01230435-1}, \eqref{equ:01230435-2}, \eqref{equ:01230435-3} and \eqref{equ:01230435-4}, we have
\begin{align}
\bP(\{\inf_{u\in [-2\tau,2\tau]}  \mathfrak{H}^t_{k}(u) \leq -10\}\cap\mathsf{E})\leq Ce^{8\tau-C^{-1}\tau^3}\bP(\mathsf{E}).
\end{align}
Together with \eqref{equ:01230456}, we get
\begin{align*}
\mathbb{P}\left(\inf_{u\in [-2\tau,2\tau]}  \mathfrak{H}^t_{k}(u)  \leq    -10\tau^2  \right)\leq&         
\bP(\{\inf_{u\in [-2\tau,2\tau]}  \mathfrak{H}^t_{k}(u) \leq -10\}\cap\mathsf{E})+\mathbb{P}\left(\mathsf{E}^{\textup{c}} \right)\\
\leq&  De^{  - D^{-1} \tau^3}.
\end{align*}
As $[-1,1]\subset [-2\tau,2\tau]$, we conclude for all $\tau\geq\max\{2,\tau_0\}$
\begin{align*}
\mathbb{P}\left(\inf_{u\in [-1,1]}  \mathfrak{H}^t_{k}(u)  \leq    -10\tau^2  \right)\leq E_k'\exp (- {c}_{k }' \tau^3).
\end{align*}
by setting $r=10\tau^2$, we have for $r\geq 10\max\{4,\tau_0^2\}$
$$\mathbb{P}\left(\inf_{u\in [-1,1]}  \mathfrak{H}^t_{k}(u)  \leq    -r  \right)\leq De^{  - D^{-1} r^{3/2}}.$$
Thus \eqref{uppertail_k} follows. 
\end{proof}
\subsection{Proof of Proposition \ref{pro:uppertail}}
For any real number  $\hat{R} $,  define the event 
\begin{align*}
\mathsf{E}_k(\hat{R}):=\left\{ \sup_{u\in [0,1/2]}\left( \mathfrak{H}_k^t(u)+2^{-1}u^2\right)\geq \hat{R} \right\}
\end{align*}
Let 
$$\chi(\hat{R})=\inf\left\{ u\in [0,1/2] |\  \mathfrak{H}_k^t(u)+2^{-1}u^2\geq \hat{R} \right\}$$
For simplicity, we denote
$$\chi=\chi(\hat{R}).$$
As $\mathsf{E}_k(\hat{R})$ occurs, $\chi\in [0,1/2]$.
For any real numbers $K$, $R$ and $R'$, define the events
\begin{align*}
\mathsf{Q}_{k-2}(K):=&\left\{\inf_{u\in [\chi,2]} \mathfrak{H}^t_{k-2}(u)+2^{-1}u^2\geq -K  \right\},\\
\mathsf{A}_{k-1,k}(R):=&\left\{  \mathfrak{H}^t_{k-1}(\chi)+2^{-1} {\chi^2} \geq -R  \right\}\cap\left\{  \mathfrak{H}^t_{j}(2)+2^{-1} {2^2} \geq -R,\ j=k-1\ \text{and}\ j=k  \right\},\\
\mathsf{B}_{k-1}(R'):=&\left\{ \sup_{u\in [\chi,2]}\left( \mathfrak{H}_{k-1}^t(u)+2^{-1}u^2 \right)\geq R' \right\}.
\end{align*}
By the stochastic monotonicity, Lemma~\ref{monotonicity}, 
\begin{align*}
\mathbb{P}(\mathsf{B}_{k-1}(R')|\mathsf{E}_k(\hat{R})\cap\mathsf{Q}_{k-2}(K)\cap\mathsf{A}_{k-1,k}(R))\geq \inf_{\chi\in [0,1/2]} p_{\chi,t}(  R,K,\hat{R},R'),
\end{align*}
where
\begin{align*}
p_{\chi,t}(  R,K,\hat{R},R')=\mathbb{P}_{\mathbf{H}_t}^{k-1,k,(\chi,2),(-R-\chi^2/2,-\hat{R}-\chi^2/2),(-R-2^2/2,-R-2^2/2),-K-u^2/2,-\infty}(\mathsf{B}_{k-1}(R')).
\end{align*}
The following proposition is a simplified version of \cite[Proposition 7.6]{CH16}
\begin{proposition}\label{pro:IvanHammond}
There exists functions $K^0(R) $, $\hat{R}^0(R,K)$ and $R'^{0}(\hat{R})$ such that the following holds. For all $R>1$, $K>K^0(R)$, $\hat{R}>\hat{R}^0(R,K)$, all $t\geq 1$ and $\chi\in [0,1/2]$,
\begin{align*}
p_{\chi,t}\left(  R,K,\hat{R},R'^{0}(\hat{R})  \right)\geq 2^{-1}.
\end{align*}
Furthermore, the functions $K^0(R) $, $\hat{R}^0(R,K)$ and $R'^{0}(\hat{R})$ are of the form
\begin{align*}
K^0(R)=\max\{R,C_1\},\ \hat{R}^0(R,K)=\max\{1600R,1600K,C_2\},\ R'^{0}(\hat{R})=6400^{-1}\hat{R}-1
\end{align*}
with some universal constants $C_1$ and $C_2$
\end{proposition}
\begin{remark}
Compared to \cite[Proposition 7.6]{CH16}, we made the following simplifications. $\mu$ is chosen to be $2^{-1}$. The extra minus one in $R'^{0}(\hat{R})$ comes from \cite[Lemma 7.7]{CH16}. On page 75 of \cite{CH16}, the choice of $R^0$ is arbitrary and we set it to be $1$. On the same page, $\delta$ can be $1600^{-1}$. The form of $K^0(R)$ and $\hat{R}^0(R,K)$ can also be found on page 75 of \cite{CH16}.
\end{remark}
In view of Proposition \ref{pro:IvanHammond}, we assume
\begin{align}\label{assumption}
R\geq 2, K\geq 2\max\{R,C_1\}, \hat{R}\geq 2\max\{ 1600R,1600K,C_2,6400 \}.
\end{align}
Then
\begin{align*}
\mathbb{P}\left(\mathsf{B}_{k-1}\left( R'^{0}(\hat{R}) \right)| \mathsf{E}_k(\hat{R})\cap\mathsf{Q}_{k-2}(K)\cap\mathsf{A}_{k-1,k}(R)\right)\geq \frac{1}{2}. 
\end{align*}
Hence
\begin{align*}
\mathbb{P}\left(  \mathsf{E}_k(\hat{R})\cap\mathsf{Q}_{k-2}(K)\cap\mathsf{A}_{k-1,k}(R)\right)\leq & 2 \mathbb{P}\left(\mathsf{B}_{k-1}\left( R'^{0}(\hat{R}) \right)\right)\\
\leq &D\exp\left(-D^{-1}\left( R'^{0}(\hat{R}) \right)^{3/2}\right)
\end{align*}
Here we used \eqref{uppertail_k} for $k-1$. Similarly, from \eqref{lowertail_k} for $k-2,k-1,k$ 
\begin{align*}
\mathbb{P}(\mathsf{Q}_{k-2}(K)^{\textup{c}})\leq&  De^{ -D^{-1}K^{3/2}},\\
 \mathbb{P}(\mathsf{A}_{k-1,k}(R)^{\textup{c}})\leq & De^{-D^{-1} R^{3/2}}.
 \end{align*}  
Thus, provided \eqref{assumption} holds, $ \mathbb{P}\left(  \mathsf{E}_k(\hat{R}) \right) $ is bounded from above by
\begin{align*}
D\exp\left(-D^{-1}\left( R'^{0}(\hat{R}) \right)^{3/2}\right)+ De^{-D^{-1} R^{3/2}}.
\end{align*}
For any $\hat{R}\geq \max\{2C_0',3200C_0,12800 \}$, take $2R=K= 3200^{-1} {\hat{R}}.$ Then \eqref{assumption} holds and $R'^{0}(\hat{R})\geq 12800^{-1}  \hat{R}$. Thus  
\begin{align*}
\mathbb{P}\left( \sup_{u\in [0,1/2]}\left( \mathfrak{H}_k^t(u)+2^{-1}u^2\right)\geq \hat{R} \right)= \mathbb{P}\left(  \mathsf{E}_k(\hat{R}) \right)\leq D\exp(-D^{-1} \hat{R}^{-3/2}).
\end{align*}
Together with the stationarity of $ \mathfrak{H}_k^t(u)+2^{-1}u^2$, we have
\begin{align*}
\mathbb{P}\left( \sup_{u\in [-1,1]}\left( \mathfrak{H}_k^t(u)+2^{-1}u^2\right)\geq \hat{R} \right)\leq D\exp(-D^{-1} \hat{R}^{-3/2}).
\end{align*}  
provided $\hat{R}\geq \max\{2C_2,3200C_1,12800 \}$. Thus \eqref{uppertail_k} follows and this finishes the proof.
\end{appendix}

\end{document}